\documentclass[a4paper,leqno,12pt]{amsart}
\usepackage[all]{xy}
\usepackage{xspace}
\usepackage{bm}
\usepackage{amsmath}
\usepackage{amstext}
\usepackage{amsfonts}
\usepackage[mathscr]{euscript}
\usepackage{amscd}
\usepackage{latexsym}
\usepackage{amssymb}
\usepackage{color}
\DeclareMathOperator{\indec}{ind}

\begin{document}
%
%
\theoremstyle{plain}
\swapnumbers
    \newtheorem{thm}{Theorem}[section]
    \newtheorem{prop}[thm]{Proposition}
    \newtheorem{lemma}[thm]{Lemma}
    \newtheorem{slemma}[thm]{Separation Lemma}
    \newtheorem{cor}[thm]{Corollary}
    \newtheorem{subsec}[thm]{}
    \newtheorem*{thma}{Theorem A}
    \newtheorem*{thmb}{Theorem B}
    \newtheorem*{propc}{Proposition C}
    \newtheorem{conj}{Conjecture}
\theoremstyle{definition}
    \newtheorem{defn}[thm]{Definition}
    \newtheorem{assume}[thm]{Assumption}
    \newtheorem{example}[thm]{Example}
    \newtheorem{examples}[thm]{Examples}
    \newtheorem{notn}[thm]{Notation}
\theoremstyle{remark}
    \newtheorem{remark}[thm]{Remark}
    \newtheorem{aside}[thm]{Aside}
    \newtheorem{ack}[thm]{Acknowledgements}

\newdir{ >}{{}*!/-5pt/@{>}}

%
%
\newenvironment{myeq}[1][]
{\stepcounter{thm}\begin{equation}\tag{\thethm}{#1}}
{\end{equation}}
\newcommand{\myeqn}[2][]
{\stepcounter{thm}\begin{equation}
     \tag{\thethm}{#1}\vcenter{#2}\end{equation}}
\newcommand{\diagr}[1]{\begin{equation*}\xymatrix{#1}\end{equation*}}
\newcommand{\diags}[1]{\begin{equation*}\xymatrix@R=10pt@C=25pt{#1}\end{equation*}}
\newcommand{\diagt}[1]{\begin{equation*}\xymatrix@R=15pt@C=25pt{#1}\end{equation*}}
\newcommand{\mydiag}[2][]{\myeq[#1]\xymatrix{#2}}
\newcommand{\mydiagram}[2][]
{\stepcounter{thm}\begin{equation}
     \tag{\thethm}{#1}\vcenter{\xymatrix{#2}}\end{equation}}
\newcommand{\myodiag}[2][]
{\stepcounter{thm}\begin{equation}
     \tag{\thethm}{#1}\vcenter{\xymatrix@R=25pt@C=15pt{#2}}\end{equation}}
\newcommand{\mypdiag}[2][]
{\stepcounter{thm}\begin{equation}
     \tag{\thethm}{#1}\vcenter{\xymatrix@R=39pt@C=5pt{#2}}\end{equation}}
\newcommand{\myqdiag}[2][]
{\stepcounter{thm}\begin{equation}
     \tag{\thethm}{#1}\vcenter{\xymatrix@R=25pt@C=21pt{#2}}\end{equation}}
\newcommand{\myrdiag}[2][]
{\stepcounter{thm}\begin{equation}
     \tag{\thethm}{#1}\vcenter{\xymatrix@R=10pt@C=25pt{#2}}\end{equation}}
\newcommand{\mysdiag}[2][]
{\stepcounter{thm}\begin{equation}
     \tag{\thethm}{#1}\vcenter{\xymatrix@R=18pt@C=25pt{#2}}\end{equation}}
\newcommand{\mytdiag}[2][]
{\stepcounter{thm}\begin{equation}
     \tag{\thethm}{#1}\vcenter{\xymatrix@R=30pt@C=25pt{#2}}\end{equation}}
\newcommand{\myudiag}[2][]
{\stepcounter{thm}\begin{equation}
     \tag{\thethm}{#1}\vcenter{\xymatrix@R=25pt@C=30pt{#2}}\end{equation}}
\newcommand{\myvdiag}[2][]
{\stepcounter{thm}\begin{equation}
     \tag{\thethm}{#1}\vcenter{\xymatrix@R=25pt@C=45pt{#2}}\end{equation}}
\newcommand{\mywdiag}[2][]
{\stepcounter{thm}\begin{equation}
     \tag{\thethm}{#1}\vcenter{\xymatrix@R=35pt@C=20pt{#2}}\end{equation}}
%
\newcommand{\supsect}[2]
{\vspace*{-5mm}\quad\\\begin{center}\textbf{{#1}}\vsm.~~~~\textbf{{#2}}\end{center}}
\newenvironment{mysubsection}[2][]
{\begin{subsec}\begin{upshape}\begin{bfseries}{#2.}
\end{bfseries}{#1}}
{\end{upshape}\end{subsec}}
\newenvironment{mysubsect}[2][]
{\begin{subsec}\begin{upshape}\begin{bfseries}{#2\vsn.}
\end{bfseries}{#1}}
{\end{upshape}\end{subsec}}
\newcommand{\sect}{\setcounter{thm}{0}\section}
\newcommand{\wh}{\ -- \ }
\newcommand{\wwh}{-- \ }
\newcommand{\w}[2][ ]{\ \ensuremath{#2}{#1}\ }
\newcommand{\ww}[1]{\ \ensuremath{#1}}
\newcommand{\wwb}[1]{\ \ensuremath{(#1)}-}
\newcommand{\wb}[2][ ]{\ (\ensuremath{#2}){#1}\ }
\newcommand{\wref}[2][ ]{\ \eqref{#2}{#1}\ }
\newcommand{\wwref}[2][ ]{\ \eqref{#2}{#1}}
\newcommand{\wbref}[2][ ]{\eqref{#2}{#1}}
%
%
\newcommand{\xra}[1]{\xrightarrow{#1}}
\newcommand{\xla}[1]{\xleftarrow{#1}}
\newcommand{\xsim}{\xrightarrow{\sim}}
\newcommand{\hra}{\hookrightarrow}
\newcommand{\epic}{\to\hspace{-5 mm}\to}
\newcommand{\xepic}[1]{\stackrel{#1}{\epic}}
\newcommand{\adj}[2]{\substack{{#1}\\ \rightleftharpoons \\ {#2}}}
\newcommand{\ccsub}[1]{\circ_{#1}}
\newcommand{\DEF}{:=}
\newcommand{\EQUIV}{\Leftrightarrow}
\newcommand{\lngra}{\longrightarrow}
\newcommand{\hsp}{\hspace{10 mm}}
\newcommand{\hs}{\hspace{5 mm}}
\newcommand{\hsm}{\hspace{3 mm}}
\newcommand{\vsm}{\vspace{3 mm}}
\newcommand{\vsn}{\vspace{1 mm}}
\newcommand{\vs}{\vspace{4 mm}}
\newcommand{\rest}[1]{\lvert_{#1}}
\newcommand{\lra}[1]{\langle{#1}\rangle}
\newcommand{\lrau}[2]{\langle\underline{#1}{#2}\rangle}
\newcommand{\li}[1]{_{({#1})}}
%
%
\newcommand{\A}{{\EuScript A}}
\newcommand{\C}{{\mathcal C}}
\newcommand{\D}{{\mathcal D}}
\newcommand{\E}{{\EuScript E}}
\newcommand{\F}{{\EuScript F}}
\newcommand{\G}{{\mathcal G}}

\newcommand{\cE}{{\mathcal E}}
\newcommand{\cI}{I}
\newcommand{\cJ}{{\mathcal J}}
\newcommand{\cK}{{\mathcal K}}
\newcommand{\LL}{{\mathcal L}}
\newcommand{\LQ}{\LL_{\bQ}}
\newcommand{\Jbar}{\overline\cJ}
\newcommand{\Jhat}{\widehat\cJ}
\newcommand{\rE}{{\mathcal E_{\ast}}}

\newcommand{\K}{{\mathcal K}}
\newcommand{\M}{{\mathcal M}}
\newcommand{\MA}{\M_{\A}}
\newcommand{\OA}{\OO^{\A}}
\newcommand{\OAp}{\OA_{+}}
\newcommand{\OAh}{\widehat{\OO}^{\A}}
\newcommand{\N}{{\EuScript N}}
\newcommand{\OO}{{\EuScript O}}
\newcommand{\Op}{\OO_{+}}
\newcommand{\Oh}{\widehat{\OO}}
\newcommand{\MO}{\M_{\OO}}
\newcommand{\MOp}{\M_{\Op}}
\newcommand{\PP}{{\mathcal P}}
\newcommand{\Pe}[1]{{\EuScript Pe}\sp{#1}}
\newcommand{\PO}{{\mathcal PO}}
\newcommand{\QQ}{{\mathcal Q}}
\newcommand{\QQp}{{\mathcal Q}_{+}}
\newcommand{\Ss}{{\mathcal S}}
\newcommand{\Sa}{\Ss_{\ast}}
\newcommand{\Sr}{\Ss_{\ast}^{\red}}
\newcommand{\TT}{{\mathcal T}}
\newcommand{\Ta}{\TT_{\ast}}
\newcommand{\U}{{\mathcal U}}
\newcommand{\V}{{\mathcal V}}
\newcommand{\eV}{{\EuScript V}}
\newcommand{\eVn}[1]{\eV\lra{#1}}
\newcommand{\W}{{\mathcal W}}
\newcommand{\X}{{\mathcal X}}
\newcommand{\eX}{{\EuScript X}}
\newcommand{\Y}{{\mathcal Y}}
\newcommand{\eY}{{\EuScript Y}}
\newcommand{\Z}{{\mathcal Z}}
\newcommand{\eZ}{{\EuScript Z}}
%
%
\newcommand{\tJ}{\widetilde{\mathbf{J}}}
\newcommand{\bJ}{\overline{\cJ}}
\newcommand{\bK}{\mathbf{K}}
\newcommand{\rM}{\overline{\mathbf{M}}}
\newcommand{\rN}{\overline{\mathbf{N}}}
\newcommand{\rP}{\overline{\mathbf{P}}}
\newcommand{\rQ}{\overline{\mathbf{Q}}}
%
%
\newcommand{\uM}{\mathbf{M}}
\newcommand{\uN}{\mathbf{N}}
\newcommand{\uP}{\mathbf{P}}
\newcommand{\uQ}{\mathbf{Q}}
\newcommand{\uR}{\mathbf{R}}
%
%
\newcommand{\bmo}{\mathbf{-1}}
\newcommand{\bze}{\mathbf{0}}
\newcommand{\bo}{\mathbf{1}}
\newcommand{\bt}{\mathbf{2}}
\newcommand{\bth}{\mathbf{3}}
\newcommand{\bfo}{\mathbf{4}}
\newcommand{\bn}{\mathbf{n}}
%
%
\newcommand{\hy}[2]{{#1}\text{-}{#2}}
%
%
\newcommand{\Alg}[1]{{#1}\text{-}{\EuScript Alg}}
\newcommand{\Ab}{{\EuScript Ab}}
\newcommand{\Abgp}{{\Ab\Grp}}
\newcommand{\AbL}{\Ab_{\Lambda}}
\newcommand{\Cat}{{\EuScript Cat}}
\newcommand{\Ch}{{\EuScript Ch}}
\newcommand{\DiGa}{{\EuScript D}i{\EuScript G}\sb{\ast}}
\newcommand{\Grp}{{\EuScript Gp}}
\newcommand{\Gpd}{{\EuScript Gpd}}
\newcommand{\RM}[1]{\hy{{#1}}{\EuScript Mod}}
\newcommand{\RL}{\RM{\Lambda}}
\newcommand{\Set}{{\EuScript Set}}
\newcommand{\Sets}{\Set_{\ast}}
\newcommand{\SC}{\hy{\Ss}{\Cat}}
\newcommand{\VC}{\hy{\V}{\Cat}}
%
%
\newcommand{\Tal}[1][ ]{$\Theta$-algebra{#1}}
\newcommand{\TAlg}{\Alg{\Theta}}
\newcommand{\Pa}[1][ ]{$\Pi$-algebra{#1}}
\newcommand{\PAa}[1][ ]{$\PiA$-algebra{#1}}
\newcommand{\PAMa}[1][ ]{$\PiAM$-algebra{#1}}
\newcommand{\PAMas}[1][ ]{$\PiAM$-algebras{#1}}
\newcommand{\PiA}{\Pi\sb{\A}}
\newcommand{\PiAM}{\Pi\sb{\A}\sp{\M}}
\newcommand{\PAlg}{\Alg{\Pi}}
\newcommand{\PAAlg}{\Alg{\PiA}}
\newcommand{\PAMAlg}{\Alg{\PiAM}}
%
%
\newcommand{\COC}{\hy{\CO}{\Cat}}
\newcommand{\GO}{(\Gpd,\OO)}
\newcommand{\GOC}{\hy{\GO}{\Cat}}
\newcommand{\GOp}{(\Gpd,\Op)}
\newcommand{\GOpC}{\hy{\GOp}{\Cat}}
\newcommand{\OC}{\hy{\OO}{\Cat}}
\newcommand{\SO}{(\Ss,\OO)}
\newcommand{\SaO}{(\Sa,\OO)}
\newcommand{\SOC}{\hy{\SO}{\Cat}}
\newcommand{\SaOC}{\hy{\SaO}{\Cat}}
\newcommand{\SOp}{(\Ss,\Op)}
\newcommand{\SOpC}{\hy{\SOp}{\Cat}}
\newcommand{\SaOp}{(\Sa,\Op)}
\newcommand{\SaOpC}{\hy{\SaOp}{\Cat}}
\newcommand{\HSO}[1]{H^{#1}\sb{\operatorname{SO}}}
\newcommand{\VO}{(\V,\OO)}
\newcommand{\VOC}{\hy{\VO}{\Cat}}
%
%
\newcommand{\co}[1]{c({#1})}
\newcommand{\Ad}{A_{\bullet}}
\newcommand{\Bd}{B_{\bullet}}
\newcommand{\tBL}{B\Lambda}
\newcommand{\Ed}{\E_{\bullet}}
\newcommand{\Fd}{F_{\bullet}}
\newcommand{\Gd}{G_{\bullet}}
\newcommand{\Kd}{\K_{\bullet}}
\newcommand{\Vd}{V_{\bullet}}
\newcommand{\whV}{\widehat{V}}
\newcommand{\hVd}[1]{\whV\sp{\lra{#1}}\sb{\bullet}}
\newcommand{\cVd}[1]{\breve{V}^{\lra{#1}}\sb{\bullet}}
\newcommand{\tVd}[1]{\tilde{V}^{\lra{#1}}\sb{\bullet}}
\newcommand{\qV}[2]{V^{\lra{#1}}_{#2}}
\newcommand{\qVd}[1]{\qV{#1}{\bullet}}
\newcommand{\qW}[2]{W^{({#1})}_{#2}}
\newcommand{\qWd}[1]{\qW{#1}{\bullet}}
\newcommand{\Wd}{W_{\bullet}}
\newcommand{\tWd}{\tilde{W}_{\bullet}}
\newcommand{\whW}{\widehat{W}}
\newcommand{\whWd}{\whW\sb{\bullet}}
\newcommand{\Xd}{X_{\bullet}}
\newcommand{\Yd}{Y_{\bullet}}
\newcommand{\sps}{semi-Postnikov section}
\newcommand{\qps}{quasi-Postnikov section}
%
%
\newcommand{\bF}{\mathbb F}
\newcommand{\Fp}{\bF_{p}}
\newcommand{\bN}{{\mathbb N}}
\newcommand{\bQ}{{\mathbb Q}}
\newcommand{\bR}{{\mathbb R}}
\newcommand{\bZ}{{\mathbb Z}}
%
%
\newcommand{\Fs}{F_{s}}
\newcommand{\Gp}{\Gamma_{+}}
\newcommand{\GpG}{\Gp,\Gamma}
\newcommand{\cGp}{\check{\Gamma}_{+}}
\newcommand{\baK}{\overline{\K}}
\newcommand{\bM}{\overline{M}}
\newcommand{\bX}{\bar{X}}
\newcommand{\bY}[1]{\overline{Y}\sb{#1}}
\newcommand{\hX}[1]{\widehat{X}\sp{#1}}
\newcommand{\cX}{\overline{X}}
\newcommand{\tX}{\widetilde{X}}
%
%
\newcommand{\pis}{\pi_{\ast}}
\newcommand{\piul}[2]{\pi\sp{#1}\sb{#2}}
\newcommand{\piAn}[1]{\piul{\A}{#1}}
\newcommand{\piA}{\piAn{\ast}}
\newcommand{\hpi}{\piul{\A}{0}}
\newcommand{\pinat}[1]{\operatorname{\pi^{\natural}_{#1}}}
\newcommand{\Po}[1]{\mathbf{P}^{#1}}
\newcommand{\bE}[2]{\mathbf{E}({#1},{#2})}
\newcommand{\tE}[2]{E({#1},{#2})}
\newcommand{\tEL}[2]{E\sb{\Lambda}({#1},{#2})}
\newcommand{\tPo}[1]{\Po{#1}}
%
%
\newcommand{\ab}{\sb{\operatorname{ab}}}
\newcommand{\arr}{\sb{\operatorname{arr}}}
\newcommand{\Arr}{\operatorname{Arr}}
\newcommand{\HAQ}[1]{H^{#1}}
\newcommand{\HL}[1]{H^{#1}\sb{\Lambda}}
\newcommand{\csk}[1]{\operatorname{csk}_{#1}}
\newcommand{\Coef}{\operatorname{Coef}}
\newcommand{\Cok}{\operatorname{Cok}}
\newcommand{\colim}{\operatorname{colim}}
\newcommand{\hocolim}{\operatorname{hocolim}}
\newcommand{\cone}[1]{\operatorname{Co}(#1)}
\newcommand{\sk}[1]{\operatorname{sk}_{#1}}
\newcommand{\cskc}[1]{\operatorname{cosk}^{c}_{#1}}
\newcommand{\fib}{\operatorname{Fib}}
\newcommand{\fiber}[1]{\operatorname{fiber}(#1)}
\newcommand{\fin}{\operatorname{fin}}
\newcommand{\hc}[1]{\operatorname{hc}_{#1}}
\newcommand{\ho}{\operatorname{ho}}
\newcommand{\holim}{\operatorname{holim}}
\newcommand{\Hom}{\operatorname{Hom}}
\newcommand{\uHom}{\underline{\Hom}}
\newcommand{\Id}{\operatorname{Id}}
\newcommand{\Image}{\operatorname{Im}}
\newcommand{\inc}{\operatorname{inc}}
\newcommand{\init}{\operatorname{init}}
\newcommand{\Ker}{\operatorname{Ker}}
\newcommand{\vf}{v_{\fin}}
\newcommand{\vfi}[1]{v{#1}_{\fin}}
\newcommand{\vi}{v_{\init}}
\newcommand{\vin}[1]{v{#1}_{\init}}
\newcommand{\bstar}{\mbox{\large $\star$}}
\newcommand{\Mor}{\operatorname{Mor}}
\newcommand{\Obj}{\operatorname{Obj}\,}
\newcommand{\op}{\sp{\operatorname{op}}}
\newcommand{\pt}{\operatorname{pt}}
\newcommand{\Ran}{\operatorname{Ran}}
\newcommand{\red}{\operatorname{red}}
\newcommand{\wPh}[1]{\widetilde{\Phi}\sb{#1}}
%
%
\newcommand{\map}{\operatorname{map}}
\newcommand{\mapp}{\map\,}
\newcommand{\mapa}{\map_{\ast}}
%
%
\newcommand{\tg}[1]{\widetilde{\gamma}_{#1}}
\newcommand{\bdz}[1]{\bar{d}^{#1}_{0}}
\newcommand{\bbd}{\mathbf{d}}
\newcommand{\bd}[1]{\bbd^{#1}_{0}}
\newcommand{\tbd}[1]{\widetilde{\bbd}^{#1}_{0}}
\newcommand{\odel}{\overline{\delta}}
\newcommand{\bDelta}{\mathbf{\Delta}}
\newcommand{\hD}[1]{\hat{\Delta}^{#1}}
\newcommand{\tD}[1]{\tilde{\Delta}^{#1}}
\newcommand{\tDl}[2]{\tilde{\Delta}^{#1}_{#2}}
\newcommand{\tDp}[1]{\tD{#1}_{+}}
\newcommand{\hDp}[1]{\hD{#1}_{+}}
\newcommand{\bG}[1]{\overline{G}\sb{#1}}
\newcommand{\bS}[1]{\mathbf{S}^{#1}}
\newcommand{\bSp}[2]{\bS{#1}_{({#2})}}
\newcommand{\bV}[1]{\overline{V}_{#1}}
\newcommand{\bW}{\overline{W}}
%
%
\newcommand{\fG}{\mathfrak{G}}
\newcommand{\cF}[1]{{\mathcal K}\sb{#1}}
\newcommand{\parz}{\partial_{0}}
\newcommand{\cfbase}[1]{\parz\cF{#1}}
\newcommand{\cftop}[1]{\widetilde{\partial}\cF{#1}}
\newcommand{\mC}[1]{\mathbf{C}\sb{#1}}
\newcommand{\baC}[1]{\overline{C}\sb{#1}}
\newcommand{\bbE}{\overline{E}'\sb{\ast}}
\newcommand{\buD}[1]{\overline{D}\sp{#1}}
\newcommand{\baD}[1]{\overline{D}\sb{#1}}
\newcommand{\baE}[1]{\overline{E}\sb{#1}}
\newcommand{\baP}[1]{\overline{P}\sb{#1}}
\newcommand{\mZ}[1]{\mathbf{Z}\sb{#1}}
\newcommand{\mZu}[1]{\mathbf{Z}\sp{#1}}
\newcommand{\latch}[1]{\mathbf{L}_{#1}}
\newcommand{\match}[1]{\mathbf{M}_{#1}}
\newcommand{\norm}[1]{\mathbf{N}_{#1}}
%
%
\newcommand{\dprod}[2][]{\displaystyle \prod_{\begin{subarray}{c}{#2} \\ {#1}\end{subarray}}}
\newcommand{\Pbar}{\overline P}

%
%
\newcommand{\smx}[1][k]{\Jul[x]{#1}}
\newcommand{\snx}[1][k]{\mathfrak n^{#1}_x}
\newcommand{\smnox}[1][k]{\mathfrak m^{#1}_{-x}}
\newcommand{\ind}[1][k]{\indec_{#1}}
\newcommand{\xcomma}[1][k]{(x \downarrow \Jul{#1})}

\newcommand{\Iul}[2][]{\cI^{#1}_{#2}}
\newcommand{\Jul}[2][]{\cJ^{#1}_{#2}}
\newcommand{\Mul}[2][]{\uM^{#1}_{#2}}
\newcommand{\Nul}[2][]{\uN^{#1}_{#2}}
\newcommand{\Pul}[2][]{\uP^{#1}_{#2}}
\newcommand{\Psul}[2]{\Psi^{#1}_{#2}}
\newcommand{\uPsi}{\overline{\Psi}}
\newcommand{\uPsul}[2]{\uPsi^{#1}_{#2}}
\newcommand{\Qul}[2][]{\uQ^{#1}_{#2}}
\newcommand{\Rul}[2][]{\uR^{#1}_{#2}}
\newcommand{\Yul}[2][]{Y^{#1}_{#2}}
\newcommand{\Jx}{\Jul[x]{}}
\newcommand{\Jxk}[1][k]{\Jul[x]{#1}}
\newcommand{\Mxk}[1][k]{\Mul[x]{#1}}
\newcommand{\Nxk}[1][k]{\Nul[x]{#1}}
\newcommand{\Pxk}[1][k]{\Pul[x]{#1}}
\newcommand{\Qxk}[1][k]{\Qul[x]{#1}}
\newcommand{\Rxk}[1][k]{\Rul[x]{#1}}
\newcommand{\rIul}[2][]{\overline\cI^{#1}_{#2}}
\newcommand{\rJul}[2][]{\overline\cJ^{#1}_{#2}}
\newcommand{\rMul}[2][]{\overline\uM^{#1}_{#2}}
\newcommand{\rNul}[2][]{\overline\uN^{#1}_{#2}}
\newcommand{\rPul}[2][]{\overline\uP^{#1}_{#2}}
\newcommand{\rQul}[2][]{\overline\uQ^{#1}_{#2}}
\newcommand{\rRul}[2][]{\overline\uR^{#1}_{#2}}
\newcommand{\rYul}[2][]{\overlineY^{#1}_{#2}}
\newcommand{\rJx}{\rJul[x]{}}
\newcommand{\rJxk}[1][k]{\rJul[x]{#1}}
\newcommand{\rMxk}[1][k]{\rMul[x]{#1}}
\newcommand{\rNxk}[1][k]{\rNul[x]{#1}}
\newcommand{\rPxk}[1][k]{\rPul[x]{#1}}
\newcommand{\rQxk}[1][k]{\rQul[x]{#1}}
\newcommand{\rRxk}[1][k]{\rRul[x]{#1}}

\newcommand{\tY}{\widetilde{Y}}
\newcommand{\hY}{\widehat{Y}}
\newcommand{\Yof}[2][]{Y_{#1}(#2)}
\newcommand{\Ybarof}[2][]{\overline{Y}_{#1}(#2)}
\newcommand{\Yk}[1][k]{\Yul[ ]{#1}}
\newcommand{\Yxk}[1][k]{\Yul[x]{#1}}
\newcommand{\tYxk}[2][x]{\tY^{#1}_{#2}} 
\newcommand{\tYk}[1]{\tY^{}_{#1}} 
\newcommand{\mxk}[2][x]{\operatorname m^{#1}_{#2}}
\newcommand{\hmxk}[2][x]{\widehat{\operatorname m}^{#1}_{#2}}
\newcommand{\rmxk}[2][x]{\overline{\operatorname m}^{#1}_{#2}}

\newcommand{\Cof}{\operatorname C} 
\newcommand{\Cofib}{\operatorname {cof}}
\newcommand{\cofib}[1]{\Cofib(#1)}
\newcommand{\Sigmap}{\Sigma'} 
\newcommand{\mapcone}[1]{\operatorname M_{#1}}

\newcommand{\sigmaul}[2]{\sigma^{#1}_{#2}}

\newcommand{\forget}{\operatorname{forget}}
\newcommand{\rforget}{\overline{\forget}}

\newcommand{\cprod}{\displaystyle \prod}
%
%
\newcommand{\Jass}[1][]{weak lattice} 

\newcommand{\Cyl}{\operatorname{Cyl}}
\newcommand{\Lift}{\operatorname{Lift}}
\newcommand{\Path}{\operatorname{Path}}
\newcommand{\proj}{\operatorname{pr}}
\newcommand{\ovp}{\overline{p}}
\newcommand{\rel}{\operatorname{rel}}
\newcommand{\Var}{\operatorname{Var}}
\newcommand{\varu}[1]{\Var_u(#1)}
\newcommand{\ovaru}[1]{\overline{\Var}_u(#1)}
\newcommand{\lift}[2]{\Lift_{#1}(#2)}
\newcommand{\Yu}{Y\langle u \rangle}
\newcommand{\Wp}{W\langle p \rangle}
\newcommand{\oYu}{\overline Y \langle u \rangle}
\newcommand{\oWpu}{\overline W \langle p,u \rangle}
\newcommand{\Prel}{P^{\rel}}
\newcommand{\siml}{\sim\sp{l}}
\newcommand{\simr}{\sim\sp{r}}
\newcommand{\rcol}[1]{\textcolor{red}{#1}}
\newcommand{\rgrn}[1]{\textcolor{green}{#1}}
\newcommand{\rmag}[1]{\textcolor{magenta}{#1}}
\newcommand{\rblue}[1]{\textcolor{blue}{#1}}

%
%
\title{A Constructive Approach to Higher Homotopy Operations}
%
%
\author[D.~Blanc]{David Blanc}
\address{Department of Mathematics\\ University of Haifa\\ 34988 Haifa\\ Israel}
\email{blanc@math.haifa.ac.il}
\author [M.W.~Johnson]{Mark W.~Johnson}
\address{Department of Mathematics\\ Penn State Altoona\\
                  Altoona, PA 16601\\ USA}
\email{mwj3@psu.edu}
\author[J.M.~Turner]{James M.~Turner}
\address{Department of Mathematics\\ Calvin College\\
         Grand Rapids, MI 49546\\ USA}
\email{jturner@calvin.edu}
\date{\today}

\begin{abstract}
  In this paper we provide an explicit general construction of higher homotopy operations
  in model categories, which include classical examples such as (long) Toda brackets
  and (iterated) Massey products, but also cover unpointed operations not usually considered
  in this context. We show how such operations, thought of as obstructions to rectifying a
  homotopy-commutative diagram, can be defined in terms of a double induction, yielding
  intermediate obstructions as well.
\end{abstract}
\maketitle

\setcounter{section}{0}

%
%
\section*{Introduction}
\label{cint}

Secondary homotopy and cohomology operations have always played an important role
in classical homotopy theory (see, e.g., \cite{AdHI,BJMahT,MPetS,PSteS} and later
\cite{GPorW,GPorH,AlldR,MOdaC,SnaiM,CWaneG}), as well as other
areas of mathematics (see \cite{AlSabtFT,FGMorgC,GLevinT,GrantT,SSterQ}).

Toda's construction of what we now call Toda brackets in \cite{TodG}
(cf.\ \cite[Ch.\ I]{TodC}) was the first example of
a secondary homotopy operation \emph{stricto sensu}, although
Adem's secondary cohomology operations and Massey's
triple products in cohomology appeared at about the same time (see \cite{AdemI,MassN}).

In \cite[Ch.\ 3]{AdHI}, Adams first tried to give a general definition of
secondary stable cohomology operations (see also \cite{HarpSC}). Kristensen gave a
description of such operations in terms of chain complexes (cf.\ \cite{KristS,KKrisS}),
which was extended by Maunder and others to $n$-th order
cohomology operations (see \cite{MaunC,HoltH,KlauC,KlauT}).

Higher operations have also figured over the years in rational homotopy theory,
where they are more accessible to computation (see, e.g., \cite{AlldR,SBasuS,RetaL,TanrH}).
In more recent years there has been a certain revival of interest in the subject,
notably in algebraic contexts (see for example, \cite{BaskT,GartH,SagaU,EfraZ,CFranH,HWickS}).

In \cite{SpanH}, Spanier gave a general theory of higher order homotopy
operations (extending the definition of secondary operations given in \cite{SpanS}).
Special cases of higher order homotopy operations appeared in
\cite{GWalkL,KraiM,MoriHT,BBGondH}, and other general definitions may be found in
\cite{BMarkH,BJTurnHH}.

The last two approaches cited present higher order operations as
the (last) obstruction to rectifying certain homotopy-commutative
diagrams (in spaces or other model categories). In particular,
they highlight the special role played by null maps in almost all
examples occurring in practice. Implicitly, they both assume an
inductive approach to rectifying such diagrams. However, in
earlier work no attempt was made to describe a useable inductive
procedure, which should (inter alia) explain precisely which
lower-order operations are required to vanish in order for a
higher order operation to be even \emph{defined}.

The goal of the present note is to make explicit the inductive process underlying
our earlier definitions of higher order operations, in as general a framework
as possible. We hope the explicit nature of this approach will help in future work
both to clarify the question of indeterminacy of the higher operations, and possibly
to produce an ``algebra of higher operations,'' in the spirit of Toda's original ``juggling
lemmas'' (see \cite[Ch.\ I]{TodC}).

An important feature of the current approach is that we assume that our indexing category
  is directed, and we consistently proceed in one direction in rectifying the given homotopy-commutative
  diagram (say, from right to left, in the ``right justified'' version).
  As a result, when we come to define the operation
  associated to an indexing category of length $n$, we use as initial data a specific choice of rectification
  for the right segment of length \w[.]{n-1} This sequence of earlier choices will appear only implicitly
  in our description and general notation for higher operations, but will be made explicit for our
    (long) Toda brackets (see \S \ref{rrjtodabr}-\ref{cohvan}).

Since our higher operations appear as obstructions to rectification, they fit into the usual
  framework of obstruction theory: when they do not vanish, one must go back along the
  thread of earlier choices until reaching a point from which one can proceed along a new branch.
  From the point of view of the obstruction theory, the important fact is their vanishing or
  non-vanishing (see Remark \ref{cohvan} for the relation to coherent vanishing).
  Nevertheless, since our higher operations are always described as a certain set of
  homotopy classes of maps into a suitable pullback, at least in some cases it is possible to describe
  the indeterminacy more explicitly. However, this would only be a part
  of the total indeterminacy, since the most general obstruction to rectification consists of the union
  of these sets, taken over all possible choices of initial data of length \w[.]{n-1}

After a brief discussion of the classical Toda bracket from our point of
view in Section \ref{cctb}, in Section \ref{cgrms}.A we describe the basic
constructions we need, associated to the type of Reedy indexing categories
for the diagrams we consider. The changes needed for pointed diagrams are discussed
in Section \ref{cgrms}.B. We give our general definition of higher order
operations in Section \ref{cgdhho}: it is hard to relate this construction to more
familiar examples, because it is intended to cover a number of different situations,
and in particular the less common unpointed version. In all cases the ``total
higher operation'' serves as an obstruction to extending a partial rectification of a
homotopy-commutative diagram one further stage in the induction.

In Section \ref{csto} we provide a refinement of this obstruction to a sequence of
intermediate steps (in an inner induction), culminating in the total operation for
the given stage in the induction.  Section \ref{crsd} is devoted to a commonly occurring problem:
rigidifying a (reduced) simplicial object in a model category, for which the simplicial identities hold only
up to homotopy. This serves to illustrate how the general (unpointed) theory works in low dimensions.

In Section \ref{cgho} we define pointed higher operations, which arise when the indexing
category has designated null maps, and we want to rectify our diagram while simultaneously
sending these to the strict zero map in the model category. This involves certain simplifications of the
general definition, as illustrated in the motivating examples of (long) Toda brackets
and Massey products, described in Section \ref{cltbmp}.

Finally, in Section \ref{cfrd} we make a tentative first step towards a possible
``algebra of higher operations,'' by showing how we can decompose our pointed higher operations into
ordinary (long) Toda brackets for a certain class of \emph{fully reduced diagrams}.

In Appendix \ref{abm} we review some basic facts in model categories needed in the paper;
Appendix \ref{abind} contains some preliminary remarks on the indeterminacy of the operations.

\begin{ack}
 We wish to thank the referee and editor for their detailed and pertinent comments.
The research of the first author was supported by Israel Science
Foundation grants 74/11 and 770/16, and the third author by National Science Foundation
grant DMS-1207746.
\end{ack}

%
%
\section{The classical Toda Bracket}
\label{cctb}

We start with a review of the classical Toda bracket, the primary example of a
pointed secondary homotopy operation. In keeping with tradition we give a
left justified description, in terms of pushouts, although for technical reasons our
general approach will be right justified, in terms of pullbacks.

\begin{mysubsection}{Left Justified Toda Brackets}

A classical \emph{Toda diagram} in any pointed model category consists of
three composable maps:
\mydiagram[\label{eqtodadiag}]{
\Yof{3} \ar[r]^h & \Yof{2} \ar[r]^g & \Yof{1} \ar[r]^f & \Yof{0}
}
\noindent with each adjacent composite left null-homotopic.
We shall assume that all objects in \wref[,]{eqtodadiag} and the analogous diagrams
  throughout the paper, are both fibrant and cofibrant, so we may disregard the distinction between
  left and right homotopy classes). To define the associated
Toda bracket, we first change $h$ into a cofibration (to avoid excessive notation, we do
not change the names of $h$ or its target).
By Lemma \ref{dlhpp} we can alter $g$ within its homotopy class to a \w{g'}
to produce a factorization:
\mydiagram[\label{eqlefttodadi}]{
  \Yof{3} \ar@{}[dr] |>{\mbox{\large{$\ulcorner$}}}  \ar@{ >->}[r]^h \ar[d] &
  \Yof{2} \ar[d] \ar@/^2.5em/[dd]^{g'} \\
\ast \ar@{ >->}[r] \ar@/_1em/[dr]^{0} & \cofib{h} \ar@{ >->}[d]^{g_2} \\
& \Yof{1}
}
\noindent so \w{g'\circ h} is the zero map (not just null-homotopic).

We use \w{i:\Yof{2}\hra\Cof\Yof{2}} (an inclusion into a reduced cone) to
extend \wref{eqlefttodadi} to the solid diagram:
\diagr{
  \Yof{3} \ar@{}[dr] |>{\mbox{\large{$\ulcorner$}}}   \ar@{ >->}[r]^h \ar[d] &
  \Yof{2} \ar@{}[dr] |>{\mbox{\large{$\ulcorner$}}}  \ar[d] \ar@/_2em/[dd]_(0.7){g'} \ar@{ >->}[r] &
\Cof\Yof{2} \ar@{}[dr] |>{\mbox{\large{$\ulcorner$}}}  \ar[d] \ar@{-->}@/^3em/[dddrr]^{\phi} \\
\ast \ar@{ >->}[r] &
\cofib{h} \ar@{}[dr] |>{\mbox{\large{$\ulcorner$}}}  \ar@{ >->}[d]^{g_2} \ar@{ >->}[r] &
\Sigmap \Yof{3} \ar@{}[dr] |>{\mbox{\large{$\ulcorner$}}}  \ar@/_2em/@{-->}[ddrr]^{\psi_{\phi}}
    \ar[r] \ar@{ >->}[d]_{j} & \ast \ar@{ >->}[d] \\
& \Yof{1} \ar@{ >->}[r]^{i} \ar@/_1em/[drrr]_{f} & \mapcone{g'} \ar[r] \ar@{-->}@/_1em/[drr]_(0.25){\kappa} &
\cofib{g_2} \ar@{.>}[dr] \\
& & & & \Yof{0}
}
\noindent where all squares (and thus all rectangles) are pushouts, with
cofibrations as indicated.

In particular, \w{\Sigmap \Yof{3}} is a model for the reduced
suspension of \w[,]{\Yof{3}} \w{\mapcone{g'}} is a mapping cone on \w[,]{g'}
and $\phi$ is a  nullhomotopy for \w[.]{f\circ g'} Note that any choice of such a
nullhomotopy $\phi$ induces maps \w{\psi_\phi:\Sigmap \Yof{3} \to \Yof{0}}
and \w[,]{\kappa:\mapcone{g'} \to \Yof{0}} with \w[.]{\kappa \circ j = \psi_\phi}
Suppose that for some choice of $\phi$, the map \w{\psi_\phi} is null-homotopic,
so \w[.]{\kappa \circ j = \psi_\phi \sim 0}  Then by Lemma \ref{dlhpp}, we could alter $\kappa$
within its homotopy class to \w{\kappa'} such that \w[,]{\kappa' \circ j = 0} whence
the pushout property for the lower right square would induce the dotted map \w[.]{\cofib{g_2} \to \Yof{0}}
As a consequence, choosing \w{f' = \kappa' \circ i \sim \kappa \circ i=f}
provides a replacement for $f$ in the same homotopy class satisfying
\w[,]{f' \circ g' = \kappa' \circ i \circ g' = 0} rather than only agreeing up to homotopy.
\end{mysubsection}

\begin{defn}\label{dljtodabr}
Given \wref[,]{eqtodadiag} the subset of the homotopy classes of maps
\w{[\Sigmap\Yof{3}, \Yof{0}]} consisting
of all classes \w{\psi_{\phi}} (for all choices of $\phi$ and \w{g_{2}} as above)
forms the \emph{Toda bracket} \w[.]{\lra{f,g,h}} Each such \w{\psi_{\phi}} is called a
\emph{value} of \w[,]{\lra{f,g,h}} and we say that the Toda bracket \emph{vanishes}
(at \w{\psi_{\phi}:\Sigmap\Yof{3}\to\Yof{0}} as above) if \w{\psi_{\phi}\sim\ast} \wwh
that is, if \w{\lra{f,g,h}} includes the null map.
\end{defn}

\begin{remark}\label{rljtodabr}
By what we have shown, \w{\lra{f,g,h}} vanishes if and only if we can vary
the spaces \w{\Yof{0},\dotsc,\Yof{3}} and the maps \w{f,g,h} within their homotopy
classes so as to make the adjacent composites in \wref{eqtodadiag} (strictly)
\emph{zero}, rather than just null-homotopic.

In fact, by considering the cofiber sequence
$$
\Yof{3} \to \Yof{2} \to \cofib{h} \to \Sigmap \Yof{3}
$$
\noindent one can show that \w{\lra{f,g,h}} is a double coset in the group
\w[:]{[\Sigmap\Yof{3},\,\Yof{0}]} In fact, the choices for homotopy classes of a
nullhomotopy for any fixed pointed map \w{\varphi:A\to B} are in one-to-one
correspondence with classes \w{[\Sigma A,\,B]} (see \cite[\S 1]{SpanS}), and
thus the contribution of the choices for $\phi$ and \w{g_{2}} respectively to
the value of \w{\lra{f,g,h}} are given by
\w{(\Sigmap h)^{\#}[\Sigmap \Yof{2}, \Yof{0}]} and
\w[,]{f_{\#}[\Sigmap \Yof{3},\Yof{1}]} respectively.

The two subgroups
\begin{myeq}\label{eqindet}
(\Sigmap h)^{\#}[\Sigmap \Yof{2}, \Yof{0}]\hsp \text{and}\hsp
 f_{\#}[\Sigmap \Yof{3}, \Yof{1}],
\end{myeq}
\noindent of \w{[\Sigmap\Yof{3}, \Yof{0}]} are referred to as the
\emph{indeterminacy} of \w[;]{\lra{f,g,h}} when \w{\Yof{3}} is a homotopy cogroup
object or \w{\Yof{1}} is a homotopy group object, the sum of \wref{eqindet} is
a subgroup of the abelian group \w[.]{[\Sigmap\Yof{3}, \Yof{0}]}

In any case, \emph{vanishing} means precisely that the (well-defined) class of
\w{\lra{f,g,h}} in the double quotient
$$
[(\Sigmap h)^{\#}\Sigmap \Yof{2}, \Yof{0}] \backslash [\Sigma\Yof{3}, \Yof{0}]
/ f_{\#}[\Sigmap \Yof{3}, \Yof{1}]
$$
\noindent is the trivial element in the quotient set.
\end{remark}

\begin{remark}\label{rrjtodabr}
The `right justified' definition of our ordinary Toda bracket is given in Step (c)
of Section \ref{cltbmp}.A below.  This will depend on a specific initial choice
of maps $f$ and $g$ with \w{f \circ g=\ast} (rather than \w[),]{f \circ g\sim\ast}
and will be denoted by \w[,]{\lrau{f,g}{,h}} so
\[
\lra{f,g,h}= \bigcup_{f \circ g=\ast} \lrau{f,g}{,h}
\]
where the union is indexed over those pairs with $f$ and $g$ in the specified homotopy classes.

The reader is advised to refer to that section
for examples of all constructions in Sections \ref{cgdhho}-\ref{csto} below,
since the example of our long Toda bracket \w{\lrau{f,g,h}{,k}} in Section \ref{cltbmp}
was the template for our more general setup.
\end{remark}

%
%
\section{Graded Reedy Matching Spaces}
\label{cgrms}

Our goal is now to extend the notions recalled in Section \ref{cctb} \wh
of Toda diagrams, and Toda brackets as obstructions to their (pointed)
realization \wh to more general diagrams \w[,]{Y:\cJ\to\E}  where $\E$ is
some complete category (eventually, a pointed model category).

\supsect{\protect{\ref{cgrms}}.A}{Reedy indexing categories}

Since our approach will
be inductive, we need to be able to filter our indexing category $\cJ$, for which
purpose we need the following notions.  Recall that a category is said to be locally finite
if each \ww{\Hom}-set is finite.

\begin{defn}\label{areedy}
We define a \emph{\Jass} to be a locally finite  Reedy indexing category $\cJ$
(see \cite[15.1]{PHirM}), equipped with a degree function \w[,]{\deg:\Obj\cJ\to\bN}
written  \w[,]{|x|=\deg(x)} such that:
\begin{itemize}
\item $\cJ$ is connected,
\item there are only finitely many objects in each degree,
\item all non-identity morphisms strictly decrease degree, and
\item every object maps to (at least) one of degree zero.
\end{itemize}
\end{defn}

\begin{remark}
Note that a \Jass[ ] $\cJ$ has no directed loops or non-trivial endomorphisms, and
\w{x\in\Obj\cJ} has only \w{\Id\sb{x}} mapping out of it if and only if \w[.]{|x|=0}
Moreover, each object is the source of only finitely many morphisms, although there may be
  elements of arbitrarily large degree.
\end{remark}

\begin{notn}\label{nreedy}
For a \Jass[ ] $\cJ$ as above:
\begin{enumerate}
\renewcommand{\labelenumi}{(\alph{enumi})~}
\item We denote by \w{\Jul{k}} the full subcategory of $\cJ$ consisting of the
objects of degree $\leq k$, with \w{\Iul{k}:\Jul{k} \to \cJ} the inclusion.
\item  For any  \w{x \in\Obj\cJ} in a positive degree, \w{\Jx} will denote the full
subcategory of $\cJ$ whose objects are those \w{t \in \cJ} with \w{\cJ(x,t)}
 non-empty. Thus \w{x\in\Jx}  and \w{\Jx\cap\Jul{0}\neq\emptyset} (by \S \ref{areedy}).
\item We denote by \w{\Jxk} the full subcategory of \w{\Jx} containing $x$ and all
objects (under $x$) of degree at most $k$, with \w{\Iul[x]{k}:\Jxk \to \Jx} the
inclusion. We implicitly assume that \w{|x|>k} when we use this notation.
Similarly, \w{\partial \Jxk} is the full subcategory of \w{\Jxk} containing
all objects other than $x$.
\item Given \w{|x|\geq k>0} and a functor \w{Y:\smx[k-1]\to\E} we have maps
$$
\sigmaul{x}{k-1}:\Yof{x} \to \dprod[|t|=k-1]{\cJ(x,t)} \Yof{t}
\hspace{.5in}\text{and}\hspace{.5in}
\sigmaul{x}{<k}:\Yof{x} \to \dprod[|t|<k]{\cJ(x,t)} \Yof{t}
$$
\noindent given by \w{\Yof{f}:\Yof{x}\to\Yof{t}} into the factor \w{\Yof{t}}
indexed by \w[.]{f:x \to t}
\item Given \w{Y:\Jxk[k-1]\to\cE} as above, there is a natural
\emph{generalized diagonal} map:
\begin{myeq}\label{eqgendiag}
\Psi=\Psul{x}{k}~:~\dprod[|v|<k]{\cJ(x,v)} \Yof{v} ~~\lngra~~
 \dprod[|s|=k]{\cJ(x,s)} \dprod[|v|<k]{\cJ(s,v)} \Yof{v}
\end{myeq}
\noindent mapping to the copy of \w{\Yof{v}} on the right with index
\w{x\stackrel{g}{\to} s \stackrel{f}{\to} v} by projection of the left hand product
onto the copy of \w{\Yof{v}} indexed by the composite \w{x \stackrel{fg}{\lngra} v}
(followed by \w[).]{\Id\sb{\Yof{v}}}
\end{enumerate}
\end{notn}

\begin{example}\label{egreedy}
Consider the following \Jass[ ] $\cJ$:
\diags{
&& a \ar[r] \ar[dr] & u \ar[r] & s \\
&x \ar[ur] \ar[dr] && v \ar[ur] \ar[dr] \\
&& b \ar[r] \ar[ur] & w \ar[r] & t\\
\deg: & 3 & 2 & 1 & 0
}

where all subdiagrams commute, and the degrees are as indicated. Then
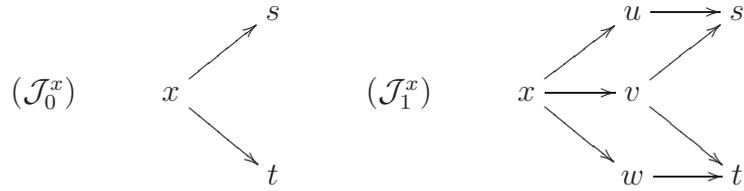
\diagt{
&& s                    &&& u \ar[r] & s    \\
(\Jxk[0])&x \ar[ur] \ar[dr] && (\Jxk[1]) &
x \ar[r] \ar[ur] \ar[dr] & v \ar[ur] \ar[dr] \\
&& t                        &&& w \ar[r] & t
}

\noindent with  \w[,]{\Jxk[2]=\cJ} and \w{\partial \Jxk[0]} is the discrete category
with objects \w[.]{\{s,t\}} Furthermore we have:
\diagt{
& a \ar[r] \ar[dr] & u \ar[r] & s       &&      u \ar[r] & s            \\
(\partial\Jxk[2]) && v \ar[ur] \ar[dr]  &&(\partial\Jxk[1])& v \ar[ur] \ar[dr]  \\
& b \ar[r] \ar[ur] & w \ar[r] & t       &&      w \ar[r] & t~.
}
\end{example}

\vsn\quad

\begin{defn}\label{dreedy}
For a \Jass[ ] $\cJ$ as above and any \w{x \in\cJ} of degree $>k$:
\begin{enumerate}
\renewcommand{\labelenumi}{(\alph{enumi})~}
\item The \emph{comma category} \w{\xcomma=(x\downarrow\partial\Jxk)} has as objects
the morphisms in $\cJ$ from $x$ to objects in \w[,]{\Jul{k}} with maps in \w{\xcomma}
given by commutative triangles in $\cJ$ of the form
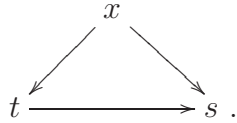
\diagr{
& x \ar[dl] \ar[dr] \\
t \ar[rr] && s ~.
}
\item For any functor \w{Y:\partial \Jxk\to\cE} and \w[,]{k<|x|} we define
the object \w{\Mul[x]{k}(Y)} (functorial in $Y$) to be the limit in $\cE$
$$
\Mxk(Y)~:=~\lim\sb{(x \downarrow \Jxk)} \widehat Y ~,
$$
\noindent where \w{\widehat Y(f:x \to s)=Y(s)} (see \cite[X.3]{MacC}).

We often write \w{\Mxk} for \w{\Mxk(Y)} when $Y$ is clear from the context.
\item For any slightly larger diagram \w[,]{Y: \Jxk \to \cE} there is a canonical map in $\cE$
  defined using the universal property of the limit,
  \w[,]{\mxk{k}(Y):\Yof{x} \to \Mxk(Y)}  and \w{\sigmaul{x}{<k+1}} is
the composite of \w{\mxk{k}} with the forgetful map (inclusion)
\diagr{
\Mxk\ \ar@{^{(}->}[rr]^-{\forget} && \dprod[|t|\leq k]{\cJ(x,t)} \Yof{t}
}
\noindent from the limit to the product, so it is closely related to the Reedy
matching map when \w[.]{k=|x|-1}

Note that \w{\Mul[x]{0}} is simply a product of entries of degree zero,
indexed by the set of maps from $x$ to the discrete category \w[,]{\Jul[x]{0}}
and \w{\mxk{0}=\sigmaul{x}{0}}.
\noindent When $\cE$ is a model category, $Y$ is called \emph{Reedy fibrant} if
each \w{\mxk{|x|-1}(Y)} is a fibration; the special case \w{k=|x|-1} is the
standard Reedy matching construction (cf.\ \cite[Defn. 15.2.3 (2)]{PHirM}).
\end{enumerate}
\end{defn}

\begin{lemma}\label{lraneq}
Given a functor \w{Y:\partial\Jxk\to\cE} as above, an extension to
\w{\overline Y:\smx[k] \to\cE} is (uniquely) determined by a choice of an object
\w[,]{\Ybarof{x}\in\cE} together with a map \w[.]{\Ybarof{x} \to \Mxk(Y)}
\end{lemma}

\begin{proof}
Recall that there is an adjoint pair given by forgetting and the right Kan extension
over \w[.]{\cI^k_x}  The fact that \w{\cI^k_x} is fully faithful implies that the
right Kan extension restricts back to the original functor (hence the term extension).
Moreover, \w{\Mxk(Y)} is the formula for the value of the right Kan extension,
\w[,]{\Ran_{\partial \Jxk}^{\Jxk} (Y)} at the
entry $x$ (see \cite[X.3, Thm 1]{MacC}).

Because of the adjunction, $\overline Y$ extends $Y$ on \w{\partial \Jxk} precisely
when there is a natural transformation
\w{\overline Y \to \Ran_{\partial \Jxk}^{\Jxk} (Y)} restricting to the identity
away from $x$. It is thus completely determined by the entry
\w[.]{\Ybarof{x}\to \Mxk(Y)}
\end{proof}

Embedding the limit \w{\Mxk(Y)} as usual into
\w[,]{\prod_{\cJ(x,u), |u|\leq k}\,\Yof{u}} we see that
there are two kinds of conditions needed for an element in this product to
be in the limit (when $\cE$ is a concrete category):
\begin{enumerate}
\renewcommand{\labelenumi}{(\alph{enumi})~}
\item Those not involving \w{\Yof{s}} with \w[,]{|s|=k} yielding
\w{\Mxk[k-1](Y)} in the lower left corner of \wref[;]{usepsi}
\item Those which do involve \w{\Yof{s}} with \w[,]{|s|=k} where the compatibility
conditions necessarily involve objects in degree $<k$, since all maps in $\cJ$
lower degree.
\end{enumerate}

This implies:

\begin{lemma}\label{lmatchpull}
If $\cJ$ is a \Jass[ ] and  \w[,]{|x|>k>0} a functor \w{Y:\partial\Jxk \to \cE} induces a pullback square:
\mydiagram[\label{usepsi}]{
  \Mxk(Y) \ar@{} [drrr] |<<{\mbox{\large{$\lrcorner$}}} \ar[d] \ar[rrr] &&&
  \dprod[|s|=k]{\cJ(x,s)} \Yof{s} \,
\ar[d]^(0.6){\prod_{\cJ(x,s)} \sigmaul{s}{<k}} \\
\Mxk[k-1](Y) \ar@{^{(}->}[rr]^{\forget} &&
\dprod[|t|<k]{\cJ(x,t)} \Yof{t} \ar[r]^-{\Psi} &
\dprod[|s|=k]{\cJ(x,s)} \dprod[|v|<k]{\cJ(s,v)} \Yof{v}.
}
\end{lemma}

Here \w{\Psi=\Psul{x}{k}} is the generalized diagonal map of \wref[,]{eqgendiag}
and the maps \w{\sigmaul{s}{<k}} on the right (given by \S \ref{nreedy}(d))
all have sources in \w[,]{\partial\Jxk} where \w{\Yof{f}} is defined.

\begin{proof}
Note that the existence of \w{Y} suffices to define each component of the diagram.
In particular, \w{\Yof{f}} is defined for each morphism $f$ in \w[,]{\partial\Jxk} and even forms
part of the definition of the factors of the right vertical, but such maps are not defined for
any \w{g:x \to v} with \w[.]{|v| \geq k}

Denote the pullback of the lower right part of the diagram by \w[.]{\Rxk} We first show that \w{\Rxk}
induces a cone on  \w[,]{(x \downarrow \Jxk)} thus inducing a map \w{\Rxk\to\Mxk}  by the universal
property of the limit: projecting off to the right for targets of degree $k$, or projecting after moving down
followed by the forgetful map for targets of lower degree, yields maps \w{\Ybarof{g}:\Rxk \to \Yof{s}}
for each \w{g:x \to s} in \w[.]{(x \downarrow \Jxk)} We must verify that whenever \w{h=f g} for \w{h:x \to t}
we have a commutative diagram in $\cE$, so that \w[.]{\Ybarof{h}=\Yof{f} \Ybarof{g}}
If the codomain of $g$ has degree less than $k$, the upper right corner is not involved, and commutativity
follows from the fact that the map from \w{\Rxk} factors through \w{\Mxk[k-1](Y)} in the lower left.
On the other hand, if the codomain of $g$ has degree exactly $k$, then projecting off at the chosen pair
\w{(g,f)} in the assumed (commutative) pullback diagram, we see that
\mydiagram{
\Rxk \ar[d]_{\Ybarof{h}} \ar[rr]^{\Ybarof{g}} && \Yof{s} \ar[d]^{\Yof{f}} \\
\Yof{t} \ar[rr]^{=} && \Yof{t}
}
\noindent commutes by the definition of the generalized diagonal $\Psi$, which establishes the cone
condition. Thus, the universal property of the limit yields a unique map \w[.]{\Rxk \to \Mxk}

On the other hand, the forgetful map \w{\forget:\Mxk \hookrightarrow \dprod[|t|\leq k]{\cJ(x,t)} \Yof{t}}
    can be split into factors with \w[,]{|t|=k} and the factors with \w[,]{|t|<k} thereby defining maps to the
    two corners of the pullback which will make the outer diagram commute, by inspection.
    Thus, there is also a map \w{\Mxk \to \Rxk} and the induced cone, as above, is the standard one, so the
    composite is the identity on \w[.]{\Mxk}

 Finally, starting from \w[,]{\Rxk} building the cone as above and then projecting as just discussed
    recovers the same maps \w{\Ybarof{h}} as entries, so this composite is the identity on \w{\Rxk} as well.
\end{proof}

\supsect{\protect{\ref{cgrms}}.B}{Pointed Graded Matching Objects}

Higher homotopy operations have traditionally appeared as obstructions
to vanishing in a pointed context, so we shall need a pointed version of the constructions above.

\begin{defn}\label{dpoint}
When $\cE$ is any category with limits (such as a model category), a
\emph{pointed object} in $\cE$ is one equipped with a map from the final object
(or empty limit), denoted by $\ast$. The most commonly occurring case is
where $\ast$ is a \emph{zero object} (both initial and final in $\cE$).
Similarly, a \emph{pointed map} in $\cE$ is one under $\ast$.  This defines the
pointed category \w{\rE} (which inherits any model category structure on $\cE$ \wh
cf.\ \cite[1.1.8]{HovM}). Note that there is a canonical zero map, also
denoted by $\ast$, between any two objects in \w[.]{\rE}
\end{defn}

\begin{defn}\label{dpindex}
We say that a small category $\cJ$ as in \S \ref{areedy} is
a \emph{pointed indexing category} if the set of morphisms  has a partition
\w{\Mor(\cJ)=\tJ\sqcup \bJ} (and thus \w{\cJ(x,t)=\tJ(x,t) \sqcup \bJ(x,t)}
for each \w[)]{x,t\in\Obj\cJ} such that:
\begin{enumerate}
\renewcommand{\labelenumi}{(\alph{enumi})~}
\item \w{\bJ(x,x)} contains \w{\Id\sb{X}} if and only if $x$ is a zero object
in $\cJ$.
\item The subsets \w{\bJ(x,t)} are absorbing under composition \wh that is,
if $f$ and $g$ are composable and either of $f$ or $g$ lies in $\bJ$, then so does
their composite. Thus $\bJ$ behaves like a (2-sided) ideal and $\tJ$ like the
corresponding cosets.
\end{enumerate}

Given \w{\rE} and a pointed \Jass[]  $\cJ$ \wh that is, a pointed indexing category which is also
  a \Jass \wh a \emph{pointed diagram} in $\rE$ is a functor \w{\Yul{}:\cJ \to \rE} such that
\w{\Yof{g}=\ast} whenever \w[.]{g \in \bJ(x,t)}
\end{defn}

\begin{example}\label{ptchain}
We can make the decreasing poset category
$$
\cJ=[n] = \{n > n-1 > \dots > 0 \}
$$
\noindent pointed by setting \w{\bJ(t,s):=\cJ(t,s)} whenever \w[,]{t-s>1}
so only indecomposable maps lie in $\tJ$. A pointed diagram \w{\cJ \to \rE} is then
simply a chain complex in $\rE$.
\end{example}

\begin{remark}\label{rpoinmat}
Making a diagram commute while also forcing certain maps to be zero is more
restrictive than simply making it commute.  Thus, we would like to construct an analog
of \w{\Mxk} tailored to the pointed case.

Note that in a pointed category \w{\rE} there is a canonical map
\w{\ast \to \cprod_{\bJ(x,t)} \Yof{t}} for any $t$, hence a section
\begin{myeq}\label{eqpointmat}
\Theta:\prod_{\tJ(x,t)} \Yof{t} \to \prod_{\cJ(x,t)} \Yof{t}
\end{myeq}
\noindent of the projection map.
\end{remark}

\begin{defn}\label{dpoinmat}
Given any diagram \w[,]{\Yul{}:\cJ \to \rE} where $\cJ$ is a pointed \Jass, define its
\emph{reduced matching space} (for $x$ and $k$) as the object of \w{\cE} defined
by the pullback:
\diagr{
  \rMxk(\Yul{}) \ar@{} [drr] |<<{\mbox{\large{$\lrcorner$}}} \ar[d]_{\rforget} \ar[rr]^{\iota^x_k}
  && \Mxk(\Yul{}) \ar[d]^\forget  \\
\dprod[|t| \leq k]{\tJ(x,t)} \Yof{t}  \ar[rr]^{\Theta} &&
\dprod[|t| \leq k]{\cJ(x,t)} \Yof{t}
}
\noindent which also determines the maps \w{\iota^x_k} and \w[.]{\rforget}
In effect, we have replaced any factor indexed on a map in $\bJ$ by $\ast$,
like reducing modulo the ideal \w[,]{\bJ}
precisely as one would expect for a pointed diagram.
\end{defn}

We then have the following analogues of Lemmas \ref{lraneq} and \ref{lmatchpull}:

\begin{lemma}\label{lptraneq}
Given a pointed functor \w[,]{Y:\partial \Jxk\to\rE} a pointed extension to
\w{\overline Y:\smx[k]\to\rE} is (uniquely) equivalent to a choice of an object
\w[,]{\Ybarof{x}} together with a morphism in \w[,]{\rE}  \w[.]{\Ybarof{x}\to\rMxk(Y)}
\end{lemma}

\begin{lemma}\label{lptmatchpull}
If \w[,]{|x|>k>0} a pointed functor \w{Y:\partial\Jxk\to\rE} (for $\cJ$ and \w{\rE} as
above) induces a pullback square:
\mytdiag[\label{ptusepsi}]{
\rMxk(Y) \ar@{} [drr] |<<{\mbox{\large{$\lrcorner$}}} \ar[d] \ar[rr] &&
\dprod[|s|=k]{\tJ(x,s)} \Yof{s} \ar[d]^(0.6){\dprod[|s|=k]{\tJ(x,s)} \sigmaul{s}{<k}} \\
\rMxk[k-1](Y) \ar@{^{(}->}[r]^-{\rforget} &
\dprod[|t|<k]{\tJ(x,t)} \Yof{t} \ar[r]^-{\uPsi}
    & \dprod[|s|=k]{\tJ(x,s)} \dprod[|v|<k]{\tJ(s,v)} \Yof{v}
}
\noindent where \w{\sigmaul{s}{<k}} is as in \S \ref{nreedy}, and
\w{\uPsi=\uPsul{x}{k}} is defined by analogy with \wref[.]{eqgendiag}
\end{lemma}

\begin{proof}
Follow the proof of Lemma \ref{lmatchpull}, with $\tJ$ replacing $\cJ$.
The absence of factors indexed in $\bJ$ implies that the structure maps
\w{\Ybarof{h}} from the pullback of \wref{ptusepsi} to the copy of \w{\Yof{s}}
indexed by \w{h:X\to s} is the zero map whenever \w[,]{h \in \bJ}
so the result follows from the absorbing property of $\bJ$.
\end{proof}

From the two lemmas we have:

\begin{cor}\label{cptmatchpull}
Any pointed diagram \w{Y: \Jxk \to\rE} induces a structure map
\w{\rmxk{k}:\Yof{x}\to\rMxk} for each \w[.]{|x|>k>0}
\end{cor}

\begin{defn}\label{dprfib}
If $\cE$ is a model category, and $\cJ$ is a pointed \Jass, a pointed diagram
\w{Y:\cJ\to\rE} is called \emph{pointed Reedy fibrant} if each map
\w{\rmxk{|x|-1}} is a fibration.
\end{defn}

\begin{lemma}\label{lptreedyfib}
If $\cE$ is a model category and $\cJ$ is a pointed \Jass, a pointed diagram
\w{Y:\cJ\to\cE} which is Reedy fibrant in the sense of \S \ref{dreedy} is
also pointed Reedy fibrant.  Moreover, for any pointed Reedy fibrant $Y$, \w{\rMxk[k](Y)} is
 fibrant in $\rE$ for each $k$.
\end{lemma}

\begin{proof}
Let \w[,]{k=|x|-1} and consider a lifting square for \w{\rmxk{k}} with respect
to an acyclic cofibration $\alpha$; extend the diagram to include \w[:]{\mxk{k}}
\diagr{
  C \ar@{ >->}[d]_{\alpha} \ar[rr] && \Yof{x}
  \ar[d]_{\rmxk{k}} \ar@{.>}[drr]^{\mxk{k}} &&\\
D \ar[rr] &&  \rMxk[k] \ar@{.>}[rr]^{\iota^x_k} &&  \Mxk[k] ~.
}
\noindent Note that a lift in the outer, distorted square will serve as a lift
for the inner square, since \w{\iota^x_k} is a base change of another monomorphism,
so is itself monic.

To show that \w{\rMxk[k](Y)} is fibrant in $\rE$ whenever $Y$ is pointed Reedy fibrant, we
adapt the argument of Lemma 15.3.9(2) through Corollary
15.3.12(2) of \cite{PHirM}, as follows:

Given a lifting diagram in $\rE$,
\mydiagram[\label{liftPTReedy}]{
C \ar@{ >->}[d]^{\sim} \ar[r] & \rMxk[n] \ar[d] \\
D \ar[r] \ar@{.>}[ur]^{h} & \ast }
\noindent we construct the dotted lift by induction on \w[.]{0\leq
k<n} For a pointed Reedy fibrant object, we assume the zero
entries are each fibrant, so their product \w{\rMxk[0]} will also
be fibrant.  For the induction step, suppose we have a lift in the
diagram
\mydiagram[\label{liftTwoPTReedy}]{
C \ar@{ >->}[d]^{\sim} \ar[r] & \rMxk[n] \ar[r] & \rMxk[k-1]  \ar[d] \\
D \ar[rr] \ar@{.>}[urr]^{h_{-1}} && \ast }
\noindent Note that the structure for any \w{f:x \to s} with
\w{|s|=k} induces a commutative diagram
\mydiagram[\label{rightPTReedy}]{
\rMxk[k-1] \ar[r] \ar[dr] & \Yof{s} \ar[d] \\
& \rMul[s]{k-1} }
\noindent so in the new lifting diagram:
\mydiagram[\label{liftThreePTReedy}]{
C \ar@{ >->}[d]^{\sim} \ar[r] & \rMxk[n] \ar[r] & \Yof{s} \ar@{->>}[d] \\
D \ar[rr] \ar@{.>}[urr]^{h_f} && \rMul[s]{k-1}
}
combining the previous two, the lift \w{h_f} exists because $Y$
was assumed to be pointed Reedy fibrant. All of these maps
together define \w[.]{h_0:D \to \prod \Yof{s}}

Compatibility with lower degree pieces then implies that \w{h_0}
factors through the limit defining \w{\rMxk[k]} which completes
our induction step, showing that \w{\rMxk[n]} is fibrant in $\rE$.
\end{proof}

\begin{lemma} \label{newPtReedy}
Each pointed diagram $Z$ has a pointed Reedy fibrant replacement
$\bY{}$ which is weakly equivalent to its Reedy fibrant
replacement $Y$ as an unpointed diagram.
\end{lemma}

\begin{proof}
In the following commuting diagram:
\diagr{
  Z(x) \ar[d]_{\alpha} \ar[rr] && \rMxk(Z) \ar[rr] && \rMxk(Y) \ar[d]\\
\Mxk(Z) \ar[rrrr] &&&&  \Mxk(Y)
}
\noindent factor the top horizontal composite as an acyclic
cofibration \w[.]{Z(x) \hra \bY{}(x)} followed by a fibration
\w[.]{\bY{}(x)\epic \rMxk(Y)} A lift in the diagram
\diagr{
Z(x) \ar@ { >->}[d]^{\sim} \ar@ { >->}[rr]^{\sim} && Y(x) \ar@{->>}[d] \\
\bY{}(x) \ar@{-->}[urr]^{\sim} \ar@{->>}[r] & \rMxk(Y) \ar[r] &
\Mxk(Y) }
\noindent will allow us to construct inductively a weak
equivalence between the new diagram $\bY{}$ and the standard Reedy
fibrant replacement $Y$ for $Z$.
\end{proof}

%
%
\section{General Definition of higher order operations}
\label{cgdhho}

{F}rom now on $\cE$ will be a model category, and we assume given a ``homotopy
commutative diagram'' in $\cE$ \wh that is, a functor
\w[,]{\tY:\cJ \to \ho(\cE)} with $\cJ$ as in \S \ref{areedy}.
Our higher homotopy operations will serve as obstructions to \emph{rectification}
of such a $\tY$ \wh that is, lifting it to \w[.]{Y:\cJ \to \cE}

We may assume for simplicity that each \w{\widetilde Y(s)} is both cofibrant
and fibrant, which can always be arranged without altering any homotopy types (see \S \ref{rassfibcof}).

\begin{mysubsection}{The double induction}
\label{sdoubind}
We attempt to construct the rectification $Y$ by a double induction:
\begin{enumerate}
\renewcommand{\labelenumi}{\Roman{enumi}.~}
\item In the outer induction, we assume we have succeeded in finding a functor
\w{\Yul{n}:\Jul{n}\to\cE} ($\Yul{n}$ \ is assumed to be Reedy fibrant), realizing
\w[.]{\tY\rest{\Jul{n}}} In fact, for our induction step it suffices to assume only
the existence of \w{\tYk{n+1}:\Jul{n+1}\to\ho(\cE)} extending \w[.]{\Yul{n}}
\item By the Reedy conditions, lifting \w{\tYk{n+1}} to
\w{\Yul{n+1}:\Jul{n+1}\to\cE} extending \w{\Yul{n}} is equivalent to
extending the latter to \emph{a point-wise extension}
\w{\Yul[x]{n}:\Jxk[n]\to\cE} for each \w{x\in\Obj\cJ}
of degree \w{n+1} separately.

Given such an $x$, the restriction of \w{\Yul{n}} produces a diagram
\w{\Yul{k}:\partial\Jxk\to\cE} for each \w{k \leq n} and the restriction of
\w{\tYxk[]{n+1}} produces a diagram \w[,]{\tYxk{k}:\Jxk \to \ho(\cE)} with the
two remaining compatible.  Thus, for our inner induction hypothesis,
assume a pointwise extension of \w{\Yul{k-1}} at $x$ (agreeing with appropriate
restrictions of both of these) has been chosen, so
\w[.]{\Yul[x]{k-1}: \Jxk[k-1]\to\cE}  Our inner induction step then asks if it is
possible to lift \w{\tYxk{k}} to \w{\Yul[x]{k}:\Jxk \to \cE} strictly extending both
\w{\Yul[x]{k-1}} and \w[,]{\Yul{k}} with the final case of the inner induction
being \w[.]{k=n}
\end{enumerate}

Notice, our inner induction step is equivalent to making coherent choices
for each homotopy class of maps out of $x$ to an object of degree $k$, leaving all
maps not involving $x$ (so those from \w[)]{\Yul{k}} or maps into objects of lower
degree (so those from \w[)]{\Yul[x]{k-1}} unchanged.
By Lemma \ref{bottomY} below, we may start the inner induction with \w{\Yul[x]{0}}
defined by the values on objects of \w[.]{\tYxk{0}} The assumption that \w{\Yul{n}}
is Reedy fibrant implies that \w{\Yul{1}} is Reedy fibrant, too, which will allow us to
use the homotopy pullback property to extend \w{\Yul[x]{0}} to \w[.]{\Yul[x]{1}}
The general step in the inner induction will use Lemma \ref{lmatchpull}:
By assumption, we have a map into the lower left corner of \wref[,]{usepsi} which we
want to extend to a map into the upper left corner still representing
the appropriate class required by \w[.]{\tYxk{k}}

\begin{remark} \label{rassfibcof}
Our induction assumption that the diagram \w{\Yul{n}} is Reedy fibrant implies that
  \w{\Yul{n}(t)} is fibrant in $\cE$ for each \w[,]{t\in\Obj\Jul{n}} and the same will  hold for the
  pullbacks that we consider below (see, e.g., \S \ref{rfibcy}).
  We will assume in addition that in the inner induction, for each \w[,]{x\in\Obj\cJ}
  \w{\Yul[x]{n}(x)} is cofibrant in $\cE$. Together this will ensure that the
  left and right homotopy classes, appearing in various results from the Appendix,
  coincide (cf.\ \cite[1.2.6]{HovM}), and the distinction can thus be disregarded.
\end{remark}

Theorem \ref{unptedThm} then yields an obstruction theory for this step in the inner
induction.
\end{mysubsection}

\begin{lemma}\label{bottomY}
In the setup described in \S \ref{sdoubind} given \w{x\in\Obj\cJ} with \w[:]{|x|>0}
\begin{enumerate}
\renewcommand{\labelenumi}{(\alph{enumi})~}
\item Any choice of representatives for a homotopy commutative
\w{\tYxk{0}:\Jxk[0] \to \ho(\cE)} provides a lift \w[.]{\Yxk[0]:\Jxk[0]\to\cE}
\item Any Reedy fibrant \w{\Yul{1}:\partial\Jxk[1] \to \cE} as above has a pointwise
  extension to a functor \w{\Yxk[1]:\Jxk[1]\to\cE} which lifts \w[.]{\tYxk{1}}
\end{enumerate}
\end{lemma}

\begin{proof}
For (a), note that \w{\Jxk[0]} has no non-trivial compositions by definition\vsm.

For (b), consider the pullback diagram
\mytdiag{
\Yof{x} \ar@{.>}[dr]_{\mxk{1}} \ar@/_2em/[ddr]_{\mxk{0}}
\ar@{-->}@/^2em/[drrr]_{\sigmaul{\widetilde Y}{1}} \\
& \Mxk[1](\Yul{1}) \ar@{} [drr] |<<{\mbox{\large{$\lrcorner$}}} \ar@{->>}[d] \ar[rr] &&
\dprod[|s|=1]{\cJ(x,s)} \Yof{s}
    \ar@{->>}[d]^(0.6){\prod\,\mxk[s]{0}(\Yul{1})}\\
& \Mxk[0](\Yxk[0]) \ar@{=}[r] & {\dprod[|t|=0]{\cJ(x,t)} \Yof{t}}
    \ar[r]^-{\Psi} &  \dprod[|s|=1]{\cJ(x,s)} \dprod[|v|=0]{\cJ(s,v)} \Yof{v} ~,
}
\noindent where the right vertical is a fibration (being a product of fibrations
by the Reedy fibrancy assumption).
This is a special case of \wref{usepsi} where the forgetful (inclusion) map
on the lower left is the identity, since \w{\partial\Jxk[0]} is discrete.

Note that the outer diagram commutes up to homotopy, since it simply compares
composites representing maps in \w{\tYxk{1}} in a somewhat unusual presentation.
By Lemma \ref{lhpp}, we can then alter the dashed map \w{\sigmaul{\widetilde Y}{1}}
within its homotopy class to obtain the dotted map \w{\mxk{1}} into \w[.]{\Mxk[1]}
Equivalently, by Lemma \ref{lraneq} one can find a representative of
\w{\tYxk{1}} extending  to \w{\smx[1]} without altering the restriction to
\w{\partial \Jxk[1]} (although this may not be the original \w[,]{\tYxk{1}}
since we might have altered \w{\sigmaul{\widetilde Y}{1}} within its homotopy class when
  applying Lemma \ref{lhpp}).
\end{proof}

\begin{remark}\label{rzeroone}
Using Lemma \ref{bottomY}, we shall henceforth assume that in the inner induction
we may start with \w[.]{k \geq 1} In order to ensure Reedy fibrancy for \w[,]{k=1}
we factor \w{\mxk{1}:\Yof{x}\to\Mxk[1]} as an acyclic cofibration
\w{\Yof{x}\hra\hY(x)} followed by a fibration \w[.]{\hmxk{1}:\hY(x)\to\Mxk[1]}
We must verify that \w{\hY(x)} and \w{\hmxk{1}} may be chosen in such a way
that the maps to the other objects \w{\tY(s)}  (with \w[)]{|s|>1} have the correct
homotopy type. However, by assumption all such objects \w{\tY(s)} are fibrant, so we
can use the left lifting property for
\diagr{
\Yof{x} \ar@ { >->}[d]_{\sim} \ar[rr]^{\alpha} && \tY(s) \ar@{->>}[d] \\
\hY(x) \ar[rr] \ar@{-->}[urr]^{\widehat{\alpha}} && \ast
}
to ensure that $\alpha$ and $\widehat{\alpha}$ have the same homotopy class.
\end{remark}

In the inner induction on $k$, we build up the diagram under the
fixed \w{x\in\Obj\cJ} by extending \w{\Yxk[k-1]} to objects in degree $k$, using:

\begin{lemma}\label{llowmatch}
Assume \w[.]{|x|>k} Given \w{\Yxk[k-1]:\Jxk[k-1] \to \cE} and \w[,]{|s|=k}
any \w{g \in \cJ(x,s)} induces a map \w[.]{\rho(g):\Yof{x} \to \Mul[s]{k-1}}
\end{lemma}

\begin{proof}
Given $g$, the diagram \w{\Yxk[k-1]} induces a cone on \w[,]{(s \downarrow \Jxk[k-1])}
sending \w{f:s \to v} to the value of \w{\Yxk[k-1]} at the target of \w[.]{f g}
Moreover, given a morphism
\diagr{
& s \ar[dl]_{f} \ar[d]^{f'} \\
v \ar[r]_{h} & u
}
\noindent in \w[,]{(s \downarrow \Jxk[k-1])} precomposition with $g$ yields
\diagr{
& x \ar[dl]_{fg} \ar[d]^{f'g} \\
v \ar[r]_{h} & u
}
\noindent which commutes in $\cJ$ \wh that is, a morphism in \w[.]{\smx[k-1]}
Applying \w{\Yxk[k-1]} yields a commutative diagram in $\cE$, showing
that we have a cone, and thus a map \w{\rho(g)} to the limit.
\end{proof}

\begin{cor}\label{clowmatch}
Combining all maps \w{\rho(g)} of  Lemma \ref{llowmatch}, a functor
\w{\Yxk[k-1]:\Jxk[k-1] \to \cE} induces a natural map
\w[.]{\displaystyle \rho_{k-1}:\Yof{x} \to \dprod[|s|=k]{\cJ(x,s)} \Mul[s]{k-1}}
\end{cor}

\begin{defn}\label{dpbgrid}
A \emph{pullback grid} is a commutative diagram tiled by squares where
each square, hence each rectangle in the diagram, is a pullback.
\end{defn}

Next, we embed the maps  \w{\rho_{k-1}} and \w{\mxk{k-1}} in a pullback grid,
in order to apply Lemma \ref{lmatchpull}:

\begin{lemma}\label{lowpiece}
Assuming \w[,]{|x|>n\geq k\geq 2} any functor \w{\Yxk[k-1]:\Jxk[k-1]\to\cE} induces
a pullback grid defined by the lower horizontal and right vertical maps,
with the natural (dashed) maps into the pullbacks:
\myudiag[\label{eqlowpiece}]{
\Yof{x} \ar@/^1.5em/[drrr]^{\rho_{k-1}} \ar@{-->}[dr]_{\beta_{k-1}}
    \ar@{-->}@/^1em/[drr]_{\eta_{k-1}} \ar@/_2em/[ddr]_{\mxk{k-1}} \\
    & \Nxk[k-1] \ar@{} [dr] |<{\mbox{\large{$\lrcorner$}}} \ar[d] \ar[r]^{q_{k-1}} &
    \Qxk[k-1] \ar@{} [dr] |<{\mbox{\large{$\lrcorner$}}} \ar[d]^{u} \ar[r]^{v}
& \dprod[|s|=k]{\cJ(x,s)} \Mul[s]{k-1} \ar[d]^(0.55){\prod\,\forget} \\
& \Mxk[k-1] \ar@{^{(}->}[r]^-{\forget} &
\dprod[|t|<k]{\cJ(x,t)} \Yof{t} \ar[r]^-{\Psi}
& \dprod[|s|=k]{\cJ(x,s)} \dprod[|v|<k]{\cJ(s,v)} \Yof{v} ~.
}
\end{lemma}

\begin{proof}
To verify commutativity of the outer diagram, note that for each
composable pair \w{x \stackrel{g}{\to} s\stackrel{f}{\to} v} in $\cJ$, the
projection of either composite from \w{\Yof{x}} onto the copy of \w{\Yof{v}}
indexed by \w{(g,f)} (in the lower right corner) is \w[,]{\Yof{fg}} by definition.
\end{proof}

We now set the stage for our obstruction theory by combining all of these pieces in a
single diagram:

\begin{prop}\label{startHHO}
Assuming \w[,]{|x|>n\geq k\geq 2} any functor \w{\Yul{k}:\partial \Jxk \to \cE}
as in \S \ref{sdoubind} induces maps into a pullback grid:
\mydiagram[\label{basicHHO}]{
\Yof{x} \ar@/_2em/[dddr]_{\mxk{k-1}}  \ar@/_1em/[ddr]_(0.7){\beta_{k-1}}
\ar@/_2em/[ddrr]^(0.7){\eta_{k-1}}
\ar@{-->}@/^1em/[drrr]^{\sigmaul{x}{k}:=\sigmaul{x}{k}(\tYxk{k})}
\ar@{.>}[dr]^(0.6){\mxk{k}}
    \ar@{.>}@/^1em/[drr]^(0.7){\alpha_k} \\
    & \Mxk \ar@{} [dr] |<<{\mbox{\large{$\lrcorner$}}} \ar@{->>}[d] \ar[r] &
    \Pxk \ar@{} [dr] |<<{\mbox{\large{$\lrcorner$}}} \ar@{->>}[d]^{p_{k-1}} \ar[r]^-{r_k} &
\dprod[|s|=k]{\cJ(x,s)} \Yof{s}
    \ar@{->>}[d]^{\prod \mxk[s]{k-1}} \\
    & \Nxk[k-1] \ar@{} [dr] |<<{\mbox{\large{$\lrcorner$}}} \ar[d] \ar[r]_{q_{k-1}} &
    \Qxk[k-1] \ar@{} [dr] |<<{\mbox{\large{$\lrcorner$}}} \ar[d] ^{u} \ar[r]^{v} &
\dprod[|s|=k]{\cJ(x,s)} \Mul[s]{k-1} \ar[d] ^{\prod\forget}\\
& \Mxk[k-1] \ar@{^{(}->}[r]^-{\forget}   &
\dprod[|t|<k]{\cJ(x,t)} \Yof{t} \ar[r]^-{\Psi}
& \dprod[|s|=k]{\cJ(x,s)} \dprod[|v|<k]{\cJ(s,v)} \Yof{v} ~.
}
\noindent Here \w{\sigmaul{x}{k}:=\sigmaul{x}{k}(\tYxk{k})} only makes the outermost
diagram commute up to homotopy.

Furthermore, the map \w{\mxk{k}} exists (after altering
\w{\sigmaul{x}{k}} within its homotopy class) if and only if there is
a map \w{\alpha_k} such that \w{p_{k-1} \alpha_k = \eta_{k-1}} and
\w[.]{r_k \alpha_k \sim \sigmaul{x}{k}}
\end{prop}

\begin{proof}
The outer pullback is \w{\Mxk} by Lemma \ref{lmatchpull} and the fact that
\w{\mxk[s]{k-1}} followed by the inclusion ``forget'' is \w{\sigmaul{s}{<k}}
(cf.\ \S \ref{nreedy}).

Note that the lower half of the grid involves only objects of $\cJ$ in degrees
\w[,]{<k} so the fact that \w{\Yul{k}} agrees with \w{\Yxk[k-1]:\Jxk[k-1]\to\cE}
implies that \w{\beta_{k-1}} and \w{\eta_{k-1}} exist, by Lemma \ref{lowpiece}.

The outer diagram commutes up to homotopy because
\w{(\Yul{k})\rest{\partial \Jxk[k-1]}} agrees with \w{\Yxk[k-1]} and lifts
\w[,]{\tYxk{k}} which is homotopy commutative.

Since the upper left square is a pullback, producing a lift of
\w{\beta_{k-1}:\Yof{x}\to\Nxk[k-1]} to \w{\Mxk} is equivalent to choosing a lift
of \w{\eta_{k-1}:\Yof{x}\to\Qxk[k-1]} to \w{\alpha_k:\Yof{x}\to\Pxk} (with
\w[).]{p_{k-1}\circ\alpha_{k}=\eta_{k-1}=q_{k-1}\circ\beta_{k-1}}

The fact that we only alter \w{\sigmaul{x}{k}} within its homotopy class
ensures that \w[,]{r_{k}\circ\alpha_{k}\sim\sigmaul{x}{k}} with the left hand
side serving as the replacement for the right hand side.
\end{proof}

\begin{remark}\label{rfibcy}
The problem here is that even though the two maps from \w{\Yof{x}} into
\w{\prod_{\cJ(x,s)}\prod_{\cJ(s,v)} \Yof{v}} (in the lower right corner of
\wref[)]{basicHHO} agree up to homotopy, this need not hold for the two maps into
\w[,]{\prod_{\cJ(x,s)}\,\Mul[s]{k-1}} the middle term on the right.  Thus we cannot
simply apply Lemma \ref{lhpp} to work with just the upper half of \wref[.]{basicHHO}

In connection with Remark \ref{rassfibcof}, one should note that all three of the objects
  along the right vertical edge
of \wref{basicHHO} are fibrant in $\cE$.  The top and bottom objects are products of entries
we assumed were fibrant.  However, the middle object is a product of the usual Reedy
matching spaces for the factors in the product above, so by \cite[Cor. 15.3.12 (2)]{PHirM},
our assumption of Reedy fibrancy implies these factors are also fibrant.

Lemma \ref{lptreedyfib} implies that this holds in the pointed case, too.
\end{remark}

\begin{mysubsection}{The Total Higher Homotopy Operation}
Following our inner induction hypothesis as in \S \ref{sdoubind}(II),
assume given \w[,]{\tYxk{k}:\Jxk\to\ho(\cE)}
\w{\Yxk[k-1]:\Jxk[k-1]\to\cE} and a Reedy fibrant \w[.]{\Yul{k}:\partial\Jxk\to\cE}

Factor the generalized diagonal map \w{\Psi=\Psul{x}{k}} of \wref{eqgendiag} as a
trivial cofibration \w{\iota:\prod_{\cJ(x,t)}\,\Yof{t}~\xra{\simeq} F^1}
followed by a fibration
\w[.]{\Psi':F^{1}\epic~\prod_{\cJ(x,s)}\,\prod_{\cJ(s,v)}~\Yof{v}}

(If we want a canonical choice of \w[,]{F^{1}} we will use the product of free
path spaces for the non-zero factors appearing in the target and the reduced path space
for each zero factor (see \S \ref{cgrms}.B), with $\iota$ defined by
the constant paths for non-zero factors.)

We then pull back the right vertical maps of \wref{basicHHO} to produce the
following pullback grid, with fibrations indicated as usual by \w[\,\,:]{\epic}
\mydiagram[\label{firstHHO}]{
\Yof{x} \ar@/_2em/[ddr]^{\eta_{k-1}} \ar@/^2em/[drrr]^{\sigmaul{x}{k}:=\sigmaul{x}{k}(\tYxk{k})}
\ar@{-->}@/_4em/[dddrr]^{\varphi} \ar@/^1em/@{-->}[drr]^{\kappa} \ar@{.>}[dr]^{\alpha_k}
\ar@/_3em/[dddr]_{\sigmaul{x}{<k}:=\sigmaul{x}{<k}(\Yxk[k-1])} \\
& \Pxk \ar@{} [dr] |<<{\mbox{\large{$\lrcorner$}}} \ar@{->>}[d]^{p_{k-1}} \ar[r]^{w}  &
F^{3} \ar@{}[dr] |<<{\mbox{\large{$\lrcorner$}}} \ar@{->>}[d]^{\mu} \ar@{->>}[r]^-{r'_k} &
\dprod[|s|=k]{\cJ(x,s)} \Yof{s} \ar@{->>}[d]^(0.5){\prod \mxk[s]{k-1}}
\ar@/^4em/[dd]^{\prod \sigmaul{s}{<k}}\\
&  \Qxk[k-1] \ar@{} [dr] |<<{\mbox{\large{$\lrcorner$}}} \ar[d] ^(0.35){u}\ar[r]^{\gamma} &
F^{2} \ar@{} [dr] |<<{\mbox{\large{$\lrcorner$}}} \ar[d]^{q} \ar@{->>}[r]^(0.35){s} &
\dprod[|s|=k]{\cJ(x,s)} \Mul[s]{k-1} \ar[d]^{\prod\forget} \\
&  \dprod[|t|<k]{\cJ(x,t)} \Yof{t} \ar[r]_(0.6){\sim}^(0.6){\iota}
\ar@/_1.5em/[rr]_{\Psi} &
F^{1} \ar@{->>}[r]^(0.3){\Psi'} &
\dprod[|s|=k]{\cJ(x,s)}\dprod[|v|<k]{\cJ(s,v)} \Yof{v}
}
\noindent where the outermost diagram commutes up to homotopy (and the map
\w{\eta_{k-1}} exists by Lemma \ref{lowpiece}).

In order to construct a lift \w[,]{\Yxk:\Jxk \to \cE} by Proposition \ref{startHHO},
we need to produce the dotted map \w{\alpha_{k}} with
\w{p_{k-1}\circ\alpha_{k}=\eta_{k-1}} and
\w[.]{r_{k}\circ\alpha_{k}=r'_{k}\circ w\circ\alpha_{k}\sim\sigmaul{x}{k}}
The problem is that the large square is a strict pullback, but not a homotopy pullback,
so the outermost diagram commuting up to homotopy is not enough.

However,  the top left square \emph{is} a pullback over a fibration, so by Lemma \ref{lhpp} producing
\w{\alpha_{k}}  is equivalent to finding a map $\kappa$ with
\w{\mu\circ\kappa\sim\gamma\circ\eta_{k-1}} and \w[.]{r'_{k}\circ\kappa\sim\sigmaul{x}{k}}

Moreover, Lemma \ref{lhpp} applies to the right vertical rectangle, which implies that  choosing
$\kappa$ is equivalent to finding a map $\varphi$ in the same homotopy class as the composite
\w[,]{\iota\circ \sigmaul{x}{<k}}  making the outer diagram commute.
Thus, the only question is whether the two composites \w{\Yof{x}\to F^{2}} agree:
that is,  given $\varphi$, with the map $\kappa$ induced by
$\varphi$ (for which necessarily \w[),]{r'_{k}\circ\kappa\sim\sigmaul{x}{k}}
is it true that \w[?]{\mu\circ\kappa\sim\gamma\circ\eta_{k-1}}
\end{mysubsection}

\begin{defn}\label{dhho}
We define the \emph{total higher homotopy operation for $x$} to be the set
\w{\lra{\Yxk[k-1]}} of all homotopy classes of maps \w{\theta:\Yof{x} \to F^{2}} with
\w{\varphi:=q\circ\theta\sim\iota\circ \sigmaul{x}{<k}} and
\w[.]{\Psi' \circ \varphi= (\prod \sigmaul{s}{<k}) \circ \sigmaul{x}{k}}
We say that\w{\lra{\Yxk[k-1]}} \emph{vanishes at} such a
\w{\theta:\Yof{x}\to F^{2}} if also \w[,]{\theta\sim\gamma\circ\eta_{k-1}}
and we say that \w{\lra{\Yxk[k-1]}} \emph{vanishes} if it vanishes at some $\theta$, or equivalently,
if this subset of the homotopy classes contains the specified class \w[.]{[\gamma\circ\eta_{k-1}]}
\end{defn}

\begin{remark}
By Corollary \ref{cthpp2} and the fact that \w{\prod \forget} is a monomorphism,
the homotopy classes \w{[\theta]} making up \w{\lra{\Yxk[k-1]}} are precisely those of the form
\w{[\mu \circ \kappa]} for a $\kappa$ with \w{r'_k \circ \kappa=\sigmaul{x}{k}} and
\w[.]{q \circ \mu \circ\kappa \sim \iota \circ \sigmaul{x}{<k}}
We may apply Corollary \ref{cthpp2} to the right vertical rectangle with horizontal
fibrations, since by assumption the outer diagram commutes up to homotopy.
This implies that the subset \w{\lra{\Yxk[k-1]}} of Definition \ref{dhho} is non-empty:
  i.e., some such $\varphi$ and so some $\kappa$ and in turn some $\theta$, exist. Thus
the total higher homotopy operation \emph{is defined} at this point.
The total higher homotopy operation \emph{vanishes} if there is such a  $\kappa$ with
\w[.]{\mu \circ \kappa \sim \gamma\circ\eta_{k-1}}
\end{remark}

This somewhat incongruous terminology of ``vanishing" is explained by the following.

\begin{prop}\label{obstructk}
Assume given \w{\tY:\cJ \to \ho(\cE)} with $\cJ$  a \Jass, and
\w{x\in\Obj\cJ} with \w[,]{|x|>n\geq k\geq 2} and let
\w[,]{\Yul{k}:\partial\Jxk\to\cE} \w[,]{\Yxk[k-1]} and \w{\tYxk{k}} be as in
\S \ref{sdoubind}. We can then extend \w{\Yul{k}} to \w{\Yxk:\Jxk \to \cE}
if and only if \w{\lra{\Yxk[k-1]}} vanishes.
\end{prop}

\begin{proof}
Note that \w{\prod\forget} is a monomorphism, since the class of monomorphisms is
closed under categorical products and the inclusion of a limit into the underlying product is always a
monomorphism.  Thus, the last statement in Corollary \ref{cthpp2} implies each value $\theta$
of \w{\lra{\Yxk[k-1]}}
satisfies \w{\theta \sim \mu \circ \kappa} for some $\kappa$ with \w{r'_{k}\circ\kappa=\sigmaul{x}{k}}
and \w[.]{q \circ\mu \circ \kappa \sim \iota \circ \sigmaul{x}{<k}}
As a consequence, if we assume
\w{\lra{\Yxk[k-1]}} vanishes at $\theta$, then there is a choice of $\kappa$ which satisfies
\w[.]{\mu \circ\kappa \sim \theta \sim \gamma \circ \eta_{k-1}}
After possibly altering $\kappa$ (and so \w[,]{\mu\circ\kappa} $\varphi$, and
\w[)]{r'_{k} \kappa} within their homotopy classes, by Lemma \ref{lhpp} applied
to the upper left square in \wref{firstHHO} we then have a dotted map \w{\alpha_{k}}
with \w[.]{p_{k-1} \circ \alpha_{k} = \eta_{k-1}}
Replacing \w{\sigmaul{x}{k}} with \w[,]{r'_k \circ \kappa'} we still have the same homotopy
  commutative diagram since \w[.]{\kappa' \sim \kappa}  Moreover, if we disregard the dashed arrows
  $\kappa$ and $\varphi$, the remaining solid diagram commutes on the nose, since
\w[,]{q \circ \mu \circ \kappa' = q \circ \gamma \circ \eta_{k-1}=\iota\circ\sigmaul{x}{<k}}
\w[,]{s \circ \gamma \circ \eta_{k-1} = s \circ \mu \circ \kappa' = \prod \mxk[s]{k-1}\circ (r'_k \circ \kappa')}
and the lower right square commutes by construction.
The upper left pullback square in \wref{basicHHO} then yields \w{\mxk{k}}
and so defines the required extension \w{\Yxk:\Jxk \to \cE} by Lemma \ref{lraneq}.

On the other hand, if \w{\lra{\Yxk[k-1]}} does not vanish, then no choice of
$\varphi$ yields a map $\kappa$ with \w[.]{\mu\circ\kappa \sim\gamma\circ\eta_{k-1}}
Thus \w{\eta_{k-1}} does not lift over \w[,]{p_{k-1}} so no such map \w{\mxk{k}}
exists. Thus there is no extension \w[,]{\Yxk} by Lemma \ref{lraneq}.
\end{proof}

\begin{remark}
As a consequence of Proposition \ref{obstructk}, our total higher homotopy operations are the
obstructions to extending a certain choice of representative of a \wwb{k}truncation of a homotopy
commutative diagram in order to produce a \wwb{k+1}truncated representative.  As in any
obstruction theory, if the obstruction does not vanish at a certain stage, we must backtrack and
reconsider earlier choices, to see whether by altering them we can make the new
  obstruction vanish at the stage in question.

It is natural to ask more generally whether there is any \wwb{k+1}truncated (strict)
representative of the given homotopy commutative diagram.  Rephrasing this in our context,
we ask whether for \emph{any} choice of a \wwb{k}truncated representative our obstruction sets
  contain the particular class which constitutes ``vanishing''. In those cases where one can identify the
  ambient collections of homotopy classes of maps with one another, a positive answer to the more
  general question is equivalent to that particular class lying in the union of our obstruction subsets.
\end{remark}

%
%
\section{Separating Total Operations}
\label{csto}

At this level of generality, we cannot expect Proposition \ref{obstructk} to be
of much help in practice: its purpose is to codify an obstruction theory for
rectifying certain homotopy-commutative diagrams, using the double induction
described in \S \ref{sdoubind}.

We now explain how to factor the right vertical map of \wref{basicHHO} or
\wref{firstHHO} as a composite of (mostly) fibrations with a view to decomposing
the obstruction \w{\lra{\Yxk[k-1]}} into more tractable pieces.
A key tool will be the following

\begin{mysubsection}{The Separation Lemma}
\label{sseparate}
Assume given a solid commutative diagram as follows\vsm:

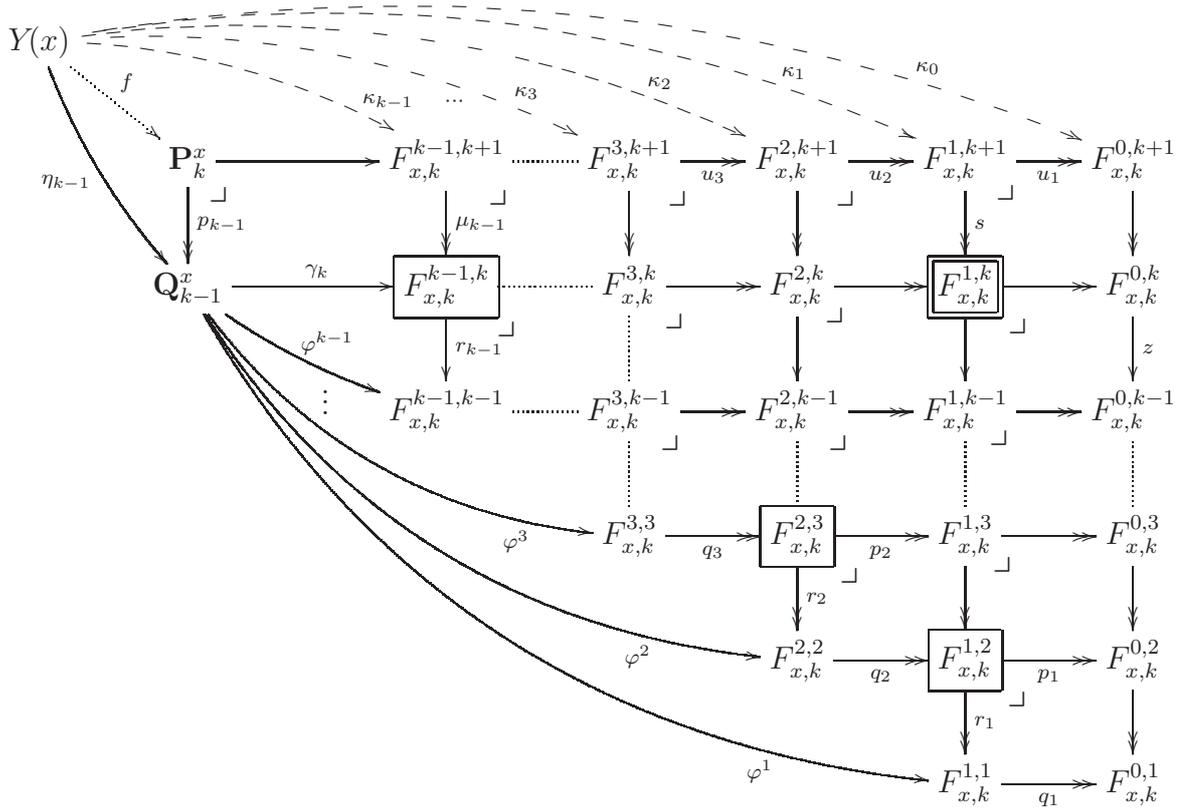
\diagr{
\Yof{x} \ar@{.>}[dr]^{f} \ar@/_/[ddr]_{\eta_{k-1}}
\ar@{-->}@/^2.5em/[drrrrrr]^(0.8){\kappa_1}
\ar@{-->}@/^2em/[drrrrr]^(0.8){\kappa_2}
\ar@{-->}@/^1.5em/[drrrr]^(0.8){\kappa_3}_(0.7){\cdots}
\ar@{-->}@/^1em/[drrr]^(0.8){\kappa_{k-1}}
\ar@{-->}@/^2.75em/[drrrrrrr]^(0.8){\kappa_{0}}
&&&&&& \\
& \Pxk \ar@{}[dr] |<<{\mbox{\large{$\lrcorner$}}} \ar@{->>}[d]^{p_{k-1}} \ar[rr] &&
F_{x,k}^{k-1,k+1} \ar@{}[dr] |<<{\mbox{\large{$\lrcorner$}}} \ar@{->>}[d]^{\mu_{k-1}}
  \ar@{.}[r] & F_{x,k}^{3,k+1} \ar@{}[dr] |<<{\mbox{\large{$\lrcorner$}}} \ar@{->>}[d] \ar@{->>}[r]_{u_3} &
F_{x,k}^{2,k+1} \ar@{}[dr] |<<{\mbox{\large{$\lrcorner$}}} \ar@{->>}[d] \ar@{->>}[r]_{u_2} &
F_{x,k}^{1,k+1} \ar@{}[dr] |<{\mbox{\large{$\lrcorner$}}} \ar@{->>}[d]^{s} \ar@{->>}[r]_{u_1} &
F_{x,k}^{0,k+1} \ar@{->>}[d] \\
& \Qxk[k-1] \ar[rr]^{\gamma_k} \ar@/_/[drr]^{\varphi^{k-1}}
\ar@/_2em/[ddrrr]_(0.8){\varphi^3}^(0.35){\vdots}
\ar@/_3em/[dddrrrr]_(0.8){\varphi^2}  \ar@/_4em/[ddddrrrrr]_(0.8){\varphi^1}
&& *+[F]{F_{x,k}^{k-1,k}} \ar@{}[dr] |<<<{\mbox{\large{$\lrcorner$}}} \ar[d]^{r_{k-1}} \ar@{.}[r] &
F_{x,k}^{3,k} \ar@{}[dr] |<<{\mbox{\large{$\lrcorner$}}} \ar@{.}[d] \ar@{->>}[r]
& F_{x,k}^{2,k} \ar[d] \ar@{}[dr] |<{\mbox{\large{$\lrcorner$}}} \ar@{->>}[r] &
*+[F=]{F_{x,k}^{1,k}} \ar@{}[dr] |<<<{\mbox{\large{$\lrcorner$}}} \ar[d] \ar@{->>}[r] &
  F_{x,k}^{0,k} \ar[d]^{z} \\
  &&& F_{x,k}^{k-1,k-1} \ar@{.}[r] &
  F_{x,k}^{3,k-1} \ar@{}[dr] |<{\mbox{\large{$\lrcorner$}}} \ar@{->>}[r] \ar@{.}[d]
  & F_{x,k}^{2,k-1} \ar@{}[dr] |<{\mbox{\large{$\lrcorner$}}} \ar@{.}[d] \ar@{->>}[r]  &
  F_{x,k}^{1,k-1} \ar@{}[dr] |>>>>>>>>>>>>>>{\mbox{\large{$\lrcorner$}}} \ar@{.}[d] \ar@{->>}[r] &
F_{x,k}^{0,k-1} \ar@{.}[d] \\
&&&& F_{x,k}^{3,3} \ar@{->>}[r]_{q_{3}} &
*+[F]{F_{x,k}^{2,3}} \ar@{}[dr] |<<<{\mbox{\large{$\lrcorner$}}} \ar@{->>}[d]^{r_{2}} \ar@{->>}[r]_{p_2} &
F_{x,k}^{1,3} \ar@{}[dr] |<{\mbox{\large{$\lrcorner$}}} \ar@{->>}[d] \ar@{->>}[r]
    & F_{x,k}^{0,3} \ar@{->>}[d]\\
&&&&& F_{x,k}^{2,2}  \ar@{->>}[r]_{q_{2}}  &
*+[F]{F_{x,k}^{1,2}} \ar@{}[dr] |<<<{\mbox{\large{$\lrcorner$}}} \ar@{->>}[d]^{r_{1}}
    \ar@{->>}[r]_{p_1} & F_{x,k}^{0,2} \ar@{->>}[d]\\
&&&&&& F_{x,k}^{1,1} \ar@{->>}[r]_{q_1} & F_{x,k}^{0,1}
}
\noindent in which:
\begin{itemize}
\item all rectangles are pullbacks,
\item the indicated maps are fibrations,
\item the objects \w{F_{x,k}^{0,1}} and \w{F_{0,k}^{j,k}} are fibrant, and
\item the vertical map $z$ is a monomorphism.
\end{itemize}

Note that as a consequence, all objects in the diagram, other than possibly \w{\Pxk} and
  \w[,]{\Qxk[k-1]} are fibrant, while all vertical maps \w{F_{x,k}^{j,k} \to F_{x,k}^{j,k-1}} are monomorphisms.

Denote the horizontal composite \w{\Qxk[k-1]\to F_{x,k}^{1,k}} by \w{\Gamma_{k-1}}
and the vertical composite \w{F_{x,k}^{j,k+1} \to F_{x,k}^{j,j+1}} by \w[,]{\Phi^{j}}
so \w[,]{\Phi^{k-1}=\mu_{k-1}} and also define \w{\varphi^k} to be the identity on
\w{\Qxk[k-1]} with \w[.]{q_k=\gamma_k}
In addition, let \w{\beta_j} denote the vertical composite \w[.]{F_{x,k}^{j,k+1} \to F_{x,k}^{j,j+2}}

Now assume that we also have a map \w{\kappa_{0}:\Yof{x}\to F_{x,k}^{0,k+1}}
such that \w[.]{\Phi^{0}\circ\kappa_{0}\sim q_1\circ \varphi^{1}\circ\eta_{k-1}}
Then by Lemma \ref{lhpp} applied to the right vertical rectangle
(with horizontal fibrations) there exists \w{\kappa_1} with \w{u_1 \circ \kappa_1 = \kappa_0}
and \w[.]{r_1 \circ \Phi^1 \circ \kappa_1 \sim \varphi^{1}\circ\eta_{k-1}}
We are interested in decomposing the question of whether
\w{s\circ\kappa_{1}\sim\Gamma_{k-1}\circ\eta_{k-1}} into a series of smaller questions.
This question will become important once we demonstrate it to be
an instance of asking for a total higher homotopy operation to vanish.

If it is true that \w[,]{q_{2}\circ\varphi^{2}\circ\eta_{k-1} \sim \Phi^{1}\circ\kappa_{1}}
then Lemma \ref{lhpp} for the next vertical rectangle imply the
existence of the dashed map \w[,]{\kappa_2} such that
\w{u_2 \circ \kappa_2 = \kappa_1}
and \w[.]{r_2 \circ \Phi^2 \circ \kappa_2 \sim \varphi^{2}\circ\eta_{k-1}}
Proceeding in this manner, and assuming the maps into the indicated ``staircase terms"
remain homotopic, even though we are only certain they agree up to homotopy after
applying the relevant \w[,]{r_j}  one produces \w{\kappa_{k-1}} such that
 \w{u_{k-1} \circ \kappa_{k-1} = \kappa_{k-2}}
 and \w[,]{r_{k-1} \circ \mu_{k-1} \circ \kappa_{k-1}=
   r_{k-1} \circ \Phi^{k-1} \circ \kappa_{k-1} \sim \varphi^{k-1}\circ\eta_{k-1}}
since \w[.]{\mu_{k-1}=\Phi^{k-1}}
The final step is then to ask whether
\w[,]{\mu_{k-1} \kappa_{k-1} \sim q_k \circ \varphi^k \circ \eta_{k-1}
  =\gamma_k \circ \eta_{k-1}} and if so, it follows by composing with most
of the rectangle across the top of the diagram that
\w[.]{s\circ\kappa_{1}\sim\Gamma_{k-1}\circ\eta_{k-1}}
In fact, we will be able to characterize when this procedure is possible in terms of
obstructions, which we will
view as ``separated" versions of the total higher homotopy operation corresponding
to the original question.
\end{mysubsection}

\begin{slemma}\label{lseparate}
Given the pullback grid as indicated above along with a choice of \w{\kappa_0}
satisfying \w[,]{\Phi^0 \circ \kappa_0 \sim q_1 \circ \varphi^1 \circ \eta_{k-1}}
there exists the indicated \w{\kappa_1} satisfying \w{u_1 \circ \kappa_1 = \kappa_0} and
\w[.]{r_1 \circ \Phi^1 \circ \kappa_1 \sim \varphi^1} Then \w{\kappa_1} also satisfies
the constraint \w{\Gamma_{k-1}\circ\eta_{k-1}\sim s\circ\kappa_{1}}
if and only if there exists an inductively chosen sequence of maps
\w{\kappa_{j}:\Yof{x} \to F_{x,k}^{j,k+1}} for \w{1\leq j<k} (starting with the
given \w[)]{\kappa_{1}} satisfying
\begin{myeq}\label{eqkappacond}
q_{j+1}\circ\varphi^{j+1}\circ\eta_{k-1}\sim \Phi^{j}\circ\kappa_{j}
\text{ and } \kappa_{j-1}=u_j \circ \kappa_j~.
\end{myeq}
\end{slemma}

The reader should note that with our conventions, in the final case \w[,]{j=k-1}
the conclusion is that \w[.]{\gamma_k \circ \eta_{k-1} \sim \mu_{k-1} \circ \kappa_{k-1}}

\begin{cor}\label{cseparate}
If either of the two equivalent conditions of Lemma \ref{lseparate} holds, then
by changing \w{\kappa_1:\Yof{x}\to F_{x,k}^{1,k+1}} within its homotopy class,
(and so using its image under \w{u_1} to replace \w{\kappa_0} within its homotopy
class as well) but without altering $\Gamma_{k-1}$, we can lift \w{\eta_{k-1}}
to the dotted map \w{f:\Yof{x}\to\Pxk} shown in the diagram.
\end{cor}

\begin{proof}[Proof of Corollary \protect{\ref{cseparate}}]
This follows from Lemma \ref{lhpp}, since the long horizontal rectangle across the
top of the diagram is a pullback over a vertical fibration.
\end{proof}

\begin{remark}\label{rseparate}
  In the case we have in mind,  \w{F_{x,k}^{0,1}} will be a product of objects \w[,]{\Yof{s}} as will
  \w[,]{F_{x,k}^{0,k+1}}  this time with \w[,]{|s|=k} and \w{F_{x,k}^{0,k}} will be the corresponding product
of matching objects \w[,]{\Mul[s]{k-1}} which will be fibrant by \cite[Cor. 15.3.12 (2)]{PHirM}.
Later, we will also have a pointed version, instead relying on pointed Reedy fibrancy
and Lemma \ref{lptreedyfib}.
Note that the second vertical map in
each column of the grid is \emph{not} required to be a fibration, but instead a
monomorphism.  Recall that monomorphisms are closed under base change and
forgetting from a limit to the underlying product is always a monomorphism, so its first
factor in any factorization must also be a monomorphism, hence these conditions will arise
naturally in our cases of interest.
\end{remark}

\begin{proof}[Proof of Lemma \protect{\ref{lseparate}}]
We will repeatedly apply Lemma \ref{lhpp} using a vertical rectangle with horizontal fibrations,
with \w{\kappa_{j-1}} as $p$ and \w{\varphi^j \circ \eta_{k-1}} as $f$, showing \w{\kappa_j} exists
and satisfies
\begin{myeq}\label{bottomcorner}
r_j \circ \Phi^j \circ \kappa_j \sim \varphi^j \circ \eta_{k-1}
\end{myeq}
provided that
\begin{myeq}\label{oneabove}
\Phi^{j-1}\circ\kappa_{j-1}  \sim q_{j}\circ\varphi^{j}\circ\eta_{k-1} ~.
\end{myeq}
Since \w{\kappa_1} exists by the assumption on \w[,]{\kappa_0}
which is really \eqref{oneabove} for \w[,]{j=1}
we begin the induction by assuming \w{\kappa_1} satisfies \eqref{oneabove} for \w[,]{j=2}
in which case \w{\kappa_2} exists and satisfies \eqref{bottomcorner} for \w[.]{j=2}
Now assuming the stricter condition \eqref{oneabove} for \w{j=3} implies the existence
of \w{\kappa_3} satisfying \eqref{bottomcorner} for \w[,]{j=3} and so on.

When our induction constructs \w{\kappa_{k-1}} satisfying
\eqref{bottomcorner} for \w[,]{j=k-1} we assume the stricter condition \eqref{oneabove} for
\w[,]{j=k} which, as noted above, is the statement that
\w[.]{\gamma_k \circ \eta_{k-1} \sim \mu_{k-1} \circ \kappa_{k-1}}
However, then composing with the horizontal rectangle across the top of the
diagram from \w{\mu_{k-1}} to $s$ implies the constraint
\w[.]{\Gamma_{k-1}\circ\eta_{k-1}\sim s\circ\kappa_{1}}

On the other hand, if \w{\kappa_1} satisfies the constraint
\w[,]{\Gamma_{k-1}\circ\eta_{k-1}\sim s\circ\kappa_{1}}
then we proceed by applying Lemma \ref{lhpp} inductively to each square along the top of the diagram
using \w{\kappa_{j-1}} for $p$ and  \w{\gamma_k \circ \eta_{k-1}} followed by the
composite \w{F_{x,k}^{k-1,k} \to F_{x,k}^{j,k}} for $f$, exploiting the horizontal fibrations
in the rectangle.  This yields \w{\kappa_j} satisfying more than \eqref{oneabove}, since the homotopy
relation is satisfied up in \w[,]{F_{x,k}^{j,k}} and this
also implies \eqref{bottomcorner} by construction.
\end{proof}

Given \w[,]{\tYxk{k}:\Jxk\to\ho(\cE)}
\w{\Yxk[k-1]:\Jxk[k-1]\to\cE} and a Reedy fibrant \w{\Yul{k}:\partial\Jxk\to\cE}
as in \S \ref{sdoubind}(II), assume that we can refine diagram \wref{firstHHO}
(used to define \w[,]{\lra{\Yxk[k-1]}} the total higher homotopy operation for $x$)
to a pullback grid as in Lemma \ref{lseparate}. Then
\w{F_{x,k}^{0,k+1}=\prod_{|s|=k}\,\prod_{\cJ(x,s)}\,\Yof{s}}
and \w[,]{F_{x,k}^{0,1}=\prod_{|s|=k}\prod_{|v|<k}\prod_{\cJ(x,s)}\prod_{\cJ(s,v)}\,\Yof{v}}
in conformity with Remark \ref{rseparate}, while one of the two equivalent conditions in
Lemma \ref{lseparate} is the vanishing of the total higher homotopy operation. Recall the vertical
composite \w{F_{x,k}^{j,k+1}\to F_{x,k}^{j,j}} in this diagram is the composite \w[.]{r_j \circ \Phi^{j}}

\begin{defn}\label{dseparate}
If we can produce a pullback grid as in Lemma \ref{lseparate} refining diagram \wref[,]{firstHHO} then
for each \w[,]{1 \leq j < k} the associated \emph{separated higher homotopy operation
for $x$ of order $j+1$}, denoted by \w[,]{\lra{\Yxk[k-1]}^{j+1}} is the set of homotopy
classes of maps \w{\theta:\Yof{x}\to F_{x,k}^{j,j+1}}
such that:
\begin{itemize}
\item if \w[,]{j < k-1} \w{r_j \circ \theta \sim \varphi^{j}\circ \eta_{k-1}} and \w{p_{j}\circ\theta}
equals the composite
\[
Y(x) \stackrel{\kappa_j}{\to} F_{x,k}^{j,k+1} \stackrel{\beta_j}{\to} F_{x,k}^{j,j+2} \text{ , or}
\]
\item if \w[,]{j=k-1} \w{r_{k-1}\circ\theta\sim\varphi^{k-1} \circ \eta_{k-1}} and
  \w{q_{k-1}\circ r_{k-1}\circ\theta=r_{k-2}\circ\mu_{k-1}\circ\kappa_{k-2}} (using the notation of the
  top two rows of vertical arrows in \S \ref{sseparate}).
\end{itemize}

We say that \w{\lra{\Yxk[k-1]}^{j+1}} \emph{vanishes} at \w{\theta:\Yof{x}\to F_{x,k}^{j,j+1}} as above if
\w{\theta\sim q_{j+1} \circ \varphi^{j+1} \circ \eta_{k-1}}
(in the notation of the Lemma), and we say it \emph{vanishes} if it vanishes at some value.
\end{defn}

Note that if we assume
\w{q_{j} \circ \varphi^{j} \circ \eta_{k-1} \sim \Phi^{j-1} \circ \kappa_{j-1}}
then by
Lemma \ref{lhpp}, \w{\kappa_j} exists, while \w{\lra{\Yxk[k-1]}^{j+1}} can then be defined and
by Corollary \ref{cthpp2} each \w{\theta^{j+1}} will satisfy \w[.]{\theta^{j+1} \sim \Phi^{j} \circ \kappa_j}
Thus, the vanishing of some value \w{\theta^{j+1}} becomes equivalent to assuming
\w[.]{q_{j+1} \circ \varphi^{j+1} \circ \eta_{k-1} \sim \Phi^{j} \circ \kappa_{j}}
In other words, the vanishing of \w{\lra{\Yxk[k-1]}^{j+1}}
(\textit{scilicet} at some map $\theta^{j+1}$) is a necessary and sufficient
condition for \w{\lra{\Yxk[k-1]}^{j+2}} to be defined.  (For comments on coherent vanishing,
see Remark \ref{cohvan}).

\begin{remark} \label{cohvan}
  Those familiar with other definitions of higher homotopy operations may have expected
  a stricter, \emph{coherent vanishing} condition in order for a subsequent operation to be defined.
    However, this need not be made explicit in our framework, as it is
  a consequence of compatibility with previous choices.

For example, our version of the ordinary Toda bracket, denoted by \w[,]{\lrau{f,g}{,h}}
is the obstruction to having a given \wwb{2}truncated commuting diagram, satisfying just
\w[,]{f \circ g=\ast} extending to a \wwb{3}truncated diagram simply
by altering $h$ within its homotopy class to satisfy \w[,]{g \circ h=\ast} \emph{without} altering $g$ or $f$.
Each choice of \wwb{2}truncation (of which there is at least one, by Lemma \ref{bottomY})
has an obstruction which is a subset of the homotopy classes of maps
\w[.]{[\Yof{3},\Omega'\Yof{0}]} The usual Toda bracket is the union of these subsets:
  \w[.]{\lra{f,g,h}=\cup \lrau{f,g}{,h}}  Thus, the more general existence question has a positive
  answer (i.e., a vanishing Toda bracket) exactly when, for \emph{some} choice of \wwb{2}truncation,
the obstruction vanishes in our sense.

When defining our long Toda brackets, say \w[,]{\lrau{f,g,h}{,k}}  we will begin by
  building the \wwb{3}truncation only if the ``front" bracket \w{\lrau{f,g}{,h}}
vanishes for some choice of \wwb{2}truncation, and we make an appropriate choice of $h$.  At that point,
we only consider values of the ``back" bracket \w{\lrau{g,h}{,k}} which use the
previously chosen maps $g$ and $h$. Thus asking that our obstruction vanish is
  automatically a kind of coherent vanishing. If it does not vanish, we must
alter our choice of \wwb{3}truncation until we obtain a coherently vanishing
``back" bracket.  Once again, one interpretation of
the traditional long Toda bracket would then be a union \w[,]{\cup\lrau{f,g,h}{,k}} this time indexed over all
possible strict rectifications of \w[,]{\lra{f,g,h}} so all such $3$-truncations.
\end{remark}

\begin{mysubsection}{Applying the Separation Lemma}
\label{sasl}
By Proposition \ref{obstructk}, a necessary and sufficient condition for the inner
induction step in \S \ref{sdoubind} is the vanishing of the total higher homotopy
operation \w{\lra{\Yxk[k-1]}} \wwh that is, by Lemma \ref{lraneq}, the existence
of a suitable map \w{\mxk{k}} in \wref[.]{basicHHO} According to
Proposition \ref{startHHO}, this in turn is equivalent to having a map
$\kappa$ in \wref{firstHHO} satisfying a certain homotopy-commutativity
requirement.

In order to apply Lemma \ref{lseparate}, we need to break up the lower right square
of \wref{firstHHO} into a pullback grid (which then induces a horizontal
decomposition of the upper right square). This will be done by decomposing
the lower right vertical map, which is a product (over \w[,]{\cJ(x,s)} with
\w[)]{|s|=k} of the forgetful maps \w{\Mul[s]{k-1}\to\prod_{\cJ(s,v)}\,\Yof{v}}
(with \w[).]{|v|\leq k-1} The target of this forgetful map can
be further broken up as in \wref{usepsi} to a product over \w{|v|=k-1} and one
over \w[.]{|v|<k-1}
\end{mysubsection}

\begin{example}\label{eglthree}
When \w[,]{|s|=3} we factor the top horizontal arrow in \wref{usepsi}
as a weak equivalence followed by a fibration:
\begin{myeq}\label{eqpoinmat}
\Mul[s]{2}~\stackrel{\sim}{\to}~F_{s,2}^{1,3}~\epic~ \dprod[|v|=2]{\cJ(s,v)}\, \Yof{v}~.
\end{myeq}
\noindent Similarly, we can factor the map in \wref{eqlowpiece} from \w{\Nul[s]{1}}
to the product of lower degree copies of \w{\Yof{t}} to produce a factorization
\begin{myeq}\label{eqothermat}
\Mul[s]{2}~\to~\Nul[s]{1}~\stackrel{\simeq}{\hra}~
G_{s,2}^{1,3}~\epic~\dprod[|t|<2]{\cJ(s,t)}\, \Yof{t}
\end{myeq}
\noindent for the lower degree forgetful map in \wref[.]{usepsi} Together these yield a factorization
of the full forgetful map:
\begin{myeq}\label{eqfforget}
\Mul[s]{2}~\to~F_{s,2}^{1,3}\times G_{s,2}^{1,3}~\epic~\dprod[|v|<3]{\cJ(s,v)}\,\Yof{v}~,
\end{myeq}
\noindent with the second map a fibration and the first necessarily a monomorphism, since
  the composite is a monomorphism as the inclusion of a limit into the underlying product.
Precomposing with structure maps \w{\Yof{s}\epic\Mul[s]{2}} (which are fibrations,
because we assumed our diagram $Y$ was Reedy fibrant) yields
\begin{myeq}\label{eqpfforget}
\dprod[|s|=3]{\cJ(x,s)}\,\Yof{s}~\epic~
\dprod[|s|=3]{\cJ(x,s)}\,\Mul[s]{2}~\to~
\dprod[|s|=3]{\cJ(x,s)}\,(F_{s,2}^{1,3}\times G_{s,2}^{1,3})~\epic~
\dprod[|s|=3]{\cJ(x,s)}\,\dprod[|v|\leq 2]{\cJ(s,v)}\,\Yof{v}~.
\end{myeq}
\noindent This is a refinement of the right column in
\wref[,]{firstHHO} in which all but the second map is a fibration, and that second map is a monomorphism.

Taking \wref{eqpfforget} as the right column in the diagram of
Lemma \ref{lseparate}, we pull it back along the bottom row of \wref{firstHHO}
to get the two right columns of the intended diagram, as shown in \wref[.]{eqegthree}

For the next column, note that the two maps out of \w{\Qxk[2]} in \wref{firstHHO}
induce a map \w[,]{\Qxk[2] \to F_{x,3}^{1,2}} in the notation of \wref[.]{eqegthree}
Factoring this as an acyclic cofibration followed by a fibration:
$$
\Qxk[2]~\stackrel{\simeq}{\hra}~F_{x,3}^{2,2}~\epic~F_{x,3}^{1,2}
$$
and taking pullbacks yields the required pullback grid:
\mydiagram[\label{eqegthree}]{
  \Pxk[3] \ar@{}[dr] |<{\mbox{\large{$\lrcorner$}}} \ar@{->>}[d]^{p_2} \ar[r] &
  F_{x,3}^{2,4} \ar@{}[dr] |<{\mbox{\large{$\lrcorner$}}} \ar@{->>}[d] \ar@{->>}[r] &
F_{x,3}^{1,4} \ar@{}[dr] |<{\mbox{\large{$\lrcorner$}}} \ar@{->>}[d] \ar@{->>}[r] &
\dprod[|s|=3]{\cJ(x,s)} \Yof{s} \ar@{->>}[d]\\
\Qxk[2] \ar[dd] \ar[r] \ar[dr]^{\sim} &
*+[F]{F_{x,3}^{2,3}} \ar@{}[dr] |<<<{\mbox{\large{$\lrcorner$}}} \ar[d] \ar@{->>}[r] &
*+[F=]{F_{x,3}^{1,3}} \ar@{}[dr] |<<<{\mbox{\large{$\lrcorner$}}} \ar[d] \ar@{->>}[r] &
\dprod[|s|=3]{\cJ(x,s)} \Mul[s]{2} \ar[d] \\
& F_{x,3}^{2,2} \ar@{->>}[r] &
*+[F]{F_{x,3}^{1,2}} \ar@{}[dr] |<<<{\mbox{\large{$\lrcorner$}}} \ar@{->>}[d] \ar@{->>}[r]
 & \dprod[|s|=3]{\cJ(x,s)} F_{s,2}^{1,3} \times G_{s,2}^{1,3} \ar@{->>}[d] \\
\dprod[|t|<3]{\cJ(x,t)} \Yof{t} \ar[rr] && F_{x,3}^{1,1} \ar@{->>}[r]
 & \dprod[|s|=3]{\cJ(x,s)}\dprod[|v|<3]{\cJ(s,v)} \Yof{v} ~.
}

Note that \w{F^{1,1}_{x,3}} is the \w{F^{1}} of Definition \ref{dhho}, while
\w{F_{x,3}^{1,3}} is \w{F^{2}} \wwh that is, the target of our total higher operation
$\theta$. Separation Lemma \ref{lseparate} tells us that this operation
vanishes precisely when the following two ``separated'' operations vanish:
\begin{enumerate}
\renewcommand{\labelenumi}{(\alph{enumi})~}
\item The first, landing in \w[,]{F_{x,3}^{1,2}} being defined by the two composite
maps from \w[;]{\Yof{x}}
\item The vanishing of the first yields a second map into \w[,]{F_{x,3}^{2,3}}
where this second map defines the values of the second of the ``separated'' operations,
and the formally defined first map defines the possible vanishing of such operations.
\end{enumerate}
\end{example}

This example is indicative of the general pattern, described by:

\begin{lemma}\label{lprodgrid}
Assume given \w[,]{\tYxk{k}:\Jxk\to\ho(\cE)}
\w{\Yxk[k-1]:\Jxk[k-1]\to\cE} and a Reedy fibrant \w{\Yul{k}:\partial\Jxk\to\cE}
as in \S \ref{sdoubind}(II). If for each\w{\Mul[s]{k-1}} we have
a pullback grid as in Lemma \ref{lseparate}, these induce a pullback grid:
\myqdiag[\label{eqsepgrid}]{
  \Pxk \ar@{}[dr] |<{\mbox{\large{$\lrcorner$}}} \ar@{->>}[d]^{p_{k-1}} \ar[r] &
  F_{x,k}^{k-1,k+1} \ar@{}[dr] |<{\mbox{\large{$\lrcorner$}}} \ar@{->>}[d] \ar@{->>}[r] &
  F_{x,k}^{k-2,k+1} \ar@{}[dr] |<{\mbox{\large{$\lrcorner$}}} \ar@{->>}[d] \ar@{.>>}[r] &
  F_{x,k}^{1,k+1} \ar@{}[dr] |<{\mbox{\large{$\lrcorner$}}} \ar@{->>}[d] \ar@{->>}[r]
    & \dprod[|s|=k]{\cJ(x,s)} \Yof{s} \ar@{->>}[d] \\
\Qxk[k-1] \ar[ddd] \ar[r] \ar@/_/[dr]^{\sim} \ar@/_2em/[ddrr]_{\sim}
& *+[F]{F_{x,k}^{k-1,k}} \ar@{}[dr] |<<<{\mbox{\large{$\lrcorner$}}} \ar[d] \ar@{->>}[r] &
F_{x,k}^{k-2,k} \ar@{}[dr] |<{\mbox{\large{$\lrcorner$}}} \ar[d] \ar@{.>>}[r]
& *+[F=]{F_{x,k}^{1,k}} \ar@{}[dr] |<<<{\mbox{\large{$\lrcorner$}}} \ar[d] \ar@{->>}[r] &
\dprod[|s|=k]{\cJ(x,s)}
    \Mul[s]{k-1} \ar[d] \\
    & F_{x,k}^{k-1,k-1} \ar@{->>}[r] &
    *+[F]{F_{x,k}^{k-2,k-1}} \ar@{}[dr] |<<<{\mbox{\large{$\lrcorner$}}} \ar@{.>>}[d] \ar@{.>>}[r]
    & F_{x,k}^{1,k-1} \ar@{}[dr] |<{\mbox{\large{$\lrcorner$}}} \ar@{.>>}[d] \ar@{->>}[r] &
\dprod[|s|=k]{\cJ(x,s)} F_{s,k-1}^{k-2,k}\times G_{s,k-1}^{k-2,k} \ar@{.>>}[d] \\
&& \ar@{.>>}[r] &
*+[F]{F_{x,k}^{1,2}} \ar@{}[dr] |<<<{\mbox{\large{$\lrcorner$}}} \ar@{->>}[d] \ar@{->>}[r] &
\dprod[|s|=k]{\cJ(x,s)} F_{s,k-1}^{1,k}\times G_{s,k-1}^{1,k} \ar@{->>}[d] \\
 \dprod[|t|<k]{\cJ(x,t)} \Yof{t} \ar[rrr]^{\sim} &&& F_{x,k}^{1,1} \ar@{->>}[r]
    & \dprod[|s|=k]{\cJ(x,s)}\dprod[|v|<k]{\cJ(s,v)} \Yof{v}
}
\noindent suitable for lifting \w{\eta_{k-1}:\Yof{x}\to\Qxk[k-1]} to \w[.]{\Pxk}
\end{lemma}

Note that the two top right slots in \wref{eqsepgrid} are consistent
with Remark \ref{rseparate}.

\begin{proof}
We prove the Lemma by induction on $k$, beginning with \wref{eqegthree} for
\w[.]{k=3} We start with a decomposition
\begin{myeq}\label{eqfacfor}
\Mul[s]{k-1}~\to~F_{s,k-1}^{k-2,k}\epic\dotsc\epic~F_{s,k-1}^{2,k}~\epic~
F_{s,k-1}^{1,k}~\epic~\dprod[|v|=k-1]{\cJ(s,v)} \Yof{v}
\end{myeq}
\noindent of the top map in \wref[,]{usepsi} where all but the first map are
fibrations; this first map is a monomorphism since the composite is such, being
the inclusion of a limit into the underlying product.
This is generated using Step \w{k-1} in the induction,
by precomposing the top row in \wref{eqsepgrid} for \w{k-1} with the map
\w{\Mul[s]{k-1}\to\Pul[s]{k-1}} of \wref[.]{basicHHO}

For \w{\Nul[s]{k-1}\to\Qul[s]{k-1}\to\prod_{|s|=k}\,\Mul[s]{k-1}}
(the middle row of \wref[),]{basicHHO} we pull back the right column of
\wref{eqsepgrid} for \w{k-1} along the generalized diagonal $\Psi$
of \wref{eqgendiag} to obtain a sequence of pullbacks
\mydiagram[\label{eqfgpull}]{
G_{s,k-1}^{j,k} \ar@{}[dr] |<{\mbox{\large{$\lrcorner$}}} \ar@{->>}[d] \ar[r] &
\dprod[|v|=k-1]{\cJ(s,v)} F_{v,k-2}^{j,k-1} \times G_{v,k-2}^{j,k-1} \ar@{->>}[d] \\
\dprod[|t|<k-1]{\cJ(s,t)} \Yof{t} \ar[r]^-{\Psi} &
\dprod[|v|=k-1]{\cJ(s,v)}\dprod[|u|<k-1]{\cJ(v,u)} \Yof{u},
}
\noindent for each \w[,]{1 \leq j \leq k-3} where the right vertical map is a
fibration by the induction assumption.

For \w[,]{j=k-2} we instead factor the composite of the top row in:
\mydiagram[\label{eqlastfact}]{
\Nul[s]{k-2} \ar[rr]^{q_{k-2}} \ar@ { >->}[rrd]^{\simeq}_{i} &&
\Qul[s]{k-2} \ar[rr] &&   G_{s,k-1}^{k-3,k} \\
&& G_{s,k-1}^{k-2,k} \ar@{->>}[rru]_{r} &&
}
\noindent into an acyclic cofibration $i$ followed by a fibration $r$, as shown
(where the top maps are those of \wref{basicHHO} and \wref{eqsepgrid} for
\w[,]{k-1} respectively).
Precomposing this with the map \w{\Mul[s]{k-1} \to \Nxk[k-2]} of \wref{basicHHO}
and then taking products as in Example \ref{eglthree} yields the
desired factorization of the forgetful map:
\begin{myeq}\label{eqotherfacfor}
\Mul[s]{k-1} \to F_{s,k-1}^{k-2,k} \times G_{s,k-1}^{k-2,k}~\dotsc
\epic~F_{s,k-1}^{2,k} \times G_{s,k-1}^{2,k}~\epic~
F_{s,k-1}^{1,k} \times
 G_{s,k-1}^{1,k}~\epic~\dprod[|v|<k]{\cJ(s,v)}\, \Yof{v}.
\end{myeq}

Now factor the next generalized diagonal \w{\Psul{x}{k-1}} as an acyclic
cofibration followed by a fibration
\w[.]{p^{1,1}:F_{x,k}^{1,1}\epic\dprod[|v|<k]{\cJ(s,v)}\,\Yof{v}}
Pulling back the tower \wref{eqotherfacfor} along \w{p^{1,1}} yields the second
column on the right in our new grid \wref[.]{eqsepgrid} The total higher operation
will then land in the twice-boxed pullback object \w[.]{F_{x,k}^{1,k}}

To construct the $j$-th column from the right \wb[,]{j\geq 2} with entries
\w[,]{F_{x,k}^{j+1,\bullet}} factor the previously defined map
\w{\Qxk[k-1]\to F_{x,k}^{j,j+1}} as an acyclic cofibration
\w{\Qxk[k-1]\stackrel{\sim}{\to} F_{x,k}^{j+1,j+1}} followed by a fibration
\w[.]{p:F_{x,k}^{j+1,j+1} \epic F_{x,k}^{j,j+1}} We then pull back the \wwb{j-1}st
column along $p$ to form the $j$-th column of \wref[.]{eqsepgrid}

Note that upon completion of this process, the map \w{\Qxk[k-1] \to F_{x,k}^{k-1,k}}
need not be a fibration, but the vertical maps in the upper left square are
fibrations, by successive base-change from the product of maps
\w[,]{\Yof{s}\epic\Mul[s]{k-1}} each of which is a fibration by Reedy fibrancy of
\w[.]{\Yul{k}}
\end{proof}

\begin{defn}\label{dsepdiag}
The diagram of Lemma \ref{lseparate}, when constructed inductively as in
Lemma \ref{lprodgrid}, will be called a \emph{separation grid} for
\w[.]{\Yul{k}}
\end{defn}

Combining Lemma \ref{lprodgrid} with the Separation Lemma \ref{lseparate}
and Corollary \ref{cseparate}
yields the following refinement of Proposition \ref{obstructk}:

\begin{thm}\label{unptedThm}
Assume given \w[,]{\tYxk{k}:\Jxk\to\ho(\cE)} \w{\Yxk[k-1]:\Jxk[k-1]\to\cE} and a
Reedy fibrant \w{\Yul{k}:\partial\Jxk\to\cE} as in \S \ref{sdoubind}(II) for
\w[.]{|x|>n\geq k\geq 2}
Then our total higher homotopy operation separates into a sequence of $k-1$
    obstructions and the following are equivalent:
\begin{enumerate}
\item A further extension to \w{\Yxk:\Jxk \to \cE} exists;
\item The total operation \w{\lra{\Yxk[k-1]}} vanishes;
\item The associated sequence \w{\lra{\Yxk[k-1]}^{j+1}} \wb{1\leq j<k} of separated
higher homotopy operations of \S \ref{dseparate} vanish
(so in particular each in turn is defined).
\end{enumerate}
\end{thm}

\begin{remark}\label{rgenpict}
The machinery of the separated higher homotopy operations has been formulated to
agree with (long) Toda brackets in pointed cases. We shall deal with
these in Section \ref{cltbmp}, after a more detailed study of the special issues
involving pointed diagrams.  In particular, the role of \w{\Qxk[k-1]} will be
played by a point, so the weak equivalence followed by a fibration
factorizations out of it will be provided by taking
reduced path objects on the target. However, we first present a simple example of
the (less familiar) general unpointed situation before focusing on the details
for the pointed situation.
\end{remark}

%
%
\section{Rigidifying Simplicial Diagrams up to Homotopy}
\label{crsd}

A commonly occurring instance of a homotopy-commutative diagram which needs to be rectified
are restricted (co)simplicial objects, also known as $\Delta$-simplicial objects
(i.e., without (co)degeneracies). Examples appear in \cite[\S 6]{BJTurnR},
\cite[\S 4.1]{BJTurnHA}, \cite[\S 5]{BlaAI}, and implicitly in \cite{MayG,SegCC,PrasmS},
and more. We now show how the double inductive approach described in \S \ref{sdoubind}
applies to such diagrams.

We denote the objects of the simplicial indexing category $\Delta$ by
\w[,]{\bze,\bo,\dotsc,\bn,\dotsc} with the value of \w{Y:\Delta\to\cE}
at $\bn$ thus denoted by \w{\Yof{\bn}} instead of the usual \w[.]{Y_{n}}

\begin{mysubsection}{$1$-Truncated $\Delta$-Simplicial Objects}
\label{sotdso}
We start the outer induction with \w[.]{n=0} Our $1$-truncated diagram in \w{\ho(\cE)}
then consists of a pair of parallel arrows, so we have only the stage \w{k=0} in
the inner induction: this means choosing representatives for each of the two
face maps \w[.]{d_0,d_1:\Yof{\bo} \to \Yof{\bze}} Making this Reedy
fibrant means changing the combined map
\w{(d_0,d_1):\Yof{\bo} \to \Yof{\bze}^{d_0} \times\Yof{\bze}^{d_1}}
into a fibration (i.e., factoring this as
\w{\Yof{\bo}\stackrel{\simeq}{\hra}\Yof{\bo}'\epic \Yof{\bze}\times\Yof{\bze}}
and replacing \w{\Yof{\bo}} by \w[).]{\Yof{\bo}'}
\end{mysubsection}

\begin{mysubsection}{$2$-Truncated $\Delta$-Simplicial Objects}
\label{sttdso}
For \w[,]{n=1} $x$ is $\bt$ and \w{\Yul{1}:\partial \Jul[1]{0} \to \cE} is
the Reedy fibrant diagram just constructed.

To define \w{\Yul[\bt]{0}:\Jul[\bt]{0} \to \cE} at stage \w{k=0} in the inner
induction, pick representatives for each of the full length composites:
in this case, the three maps \w{\Yof{\bt}\to\Yof{\bze}} denoted by \w[,]{d_0 d_1}
\w[,]{d_0 d_2} and \w{d_1 d_2} in canonical form.
This means \w{\Mul[\bt]{0}} is the product of three copies of \w{\Yof{\bze}} indexed
by \w{d_i d_j} \wb[,]{0 \leq i < j \leq 2} and our choice of representatives yields
a single map \w{\mxk[\bt]{0}} into the product.

At stage \w[,]{k=1} we must first choose representatives for the components of
\w{\sigmaul{\bt}{1}(\tYxk[\bt]{1})} \wwh that is, for the maps \w[,]{d_{0}}
\w[,]{d_{1}} and \w[,]{d_{2}:\Yof{\bt}\to\Yof{\bo}} which are all the maps
\w{\bt\to\bo} in $\cJ$).
The generalized diagonal map \w{\Psi=\Psul{\bt}{1}} of \wref{eqgendiag} takes
\w{\Yof{\bze}^{d_i d_j}} \wb{i<j} to the product
\w[,]{\Yof{\bze}^{d_i d_j}\times \Yof{\bze}^{d_{j-1} d_i}} in accordance with the
simplicial identities. Note that the target of \w{\sigmaul{\bt}{1}}
is \w[.]{\prod_{0 \leq j \leq 2}\,\Yof{\bo}^{d_j}}

Thus we have a pair of maps into a pullback diagram:
\mydiagram[\label{eqtwtrunc}]{
\Yof{\bt} \ar@/_1em/[ddr]_{\mxk[\bt]{0}}
\ar@/^1em/@{-->}[drrr]^{\sigmaul{\bt}{1}(\tYxk[\bt]{1})=(d_0,d_1,d_2)}
\ar@{.>}[dr]^{\mxk[\bt]{1}}\\
& \Mul[\bt]{1} \ar@{}[drr] |<<<{\mbox{\large{$\lrcorner$}}} \ar@{->>}[d] \ar[rr] &&
\cprod_{j \leq 2} \Yof{\bo}^{d_j} \ar@{->>}[d] \\
& \Mul[\bt]{0} \ar@{=}[r] & \cprod_{i<j\leq2} \Yof{\bze}^{d_i d_j}  \ar[r]^-{\Psi} &
\cprod_{j \leq 2}\ \cprod_{i \leq 1}\, \Yof{\bze}^{d_i d_j} ~.
}
\noindent where the outer diagram commutes up to homotopy (for any choice of
representatives for \w[,]{d_{0}} \w[,]{d_{1}} and \w[).]{d_{2}}
The dotted map exists by Lemma \ref{lhpp} (after possibly altering the dashed map
within its homotopy class), yielding a full $2$-truncated $\Delta$-simplicial object
(which rectifies  \w[)]{\tYxk[\bt]{1}} by Lemma \ref{lraneq}. Changing
\w{\mxk[\bt]{1}} into a fibration provides us with a Reedy fibrant replacement
\w[.]{\Yul{2}:\partial\Jul{2}\to\cE}
\end{mysubsection}

\begin{mysubsection}{$3$-Truncated $\Delta$-Simplicial Objects}
\label{sthtdso}
At stage \w{n=2} (with \w[),]{x=\bth}  for the first time we are in the situation
of \S \ref{dhho}, somewhat simplified by the fact that we have a single
object $\bn$ in each grading $n$ of \w[.]{\cJ=\Delta} In particular, we will have
no separated operations yet.

In the inner induction, for \w[,]{k=0} we choose
representatives for each full length map in \w{\tYxk[\bth]{2}} to obtain
\w[;]{\Yul[\bth]{0}}  the full length composites are the four maps
\w{d_i d_j d_{\ell}}
with \w[,]{0 \leq i < j < \ell \leq 3} so \w{\Mul[\bth]{0}} is a
product of four copies of \w{\Yof{0}} indexed by these maps, and the generalized
diagonal of \wref{eqgendiag} takes each copy of
\w{\Yof{\bze}^{d_{i}d_{j}d_{\ell}}} to the product
$$
\Yof{\bze}^{d_{i}d_{j}d_{\ell}} \times
\Yof{\bze}^{d_{j-1}d_{i}d_{\ell}}\times
\Yof{\bze}^{d_{\ell-2} d_{i}d_{j}}~.
$$
\noindent We make an initial choice (to be modified below) of
\w{\sigmaul{\bt}{1}(\tYxk[\bt]{1})} (i.e., of each composite
\w{d_{j}d_{\ell}:\bth\to\bo} for \w{0\leq j<\ell\leq 3} within its homotopy class).
Again this yields a pair of maps into a pullback diagram:
\mydiagram[\label{eqtwotr}]{
  \Yof{\bth} \ar@/_1em/[ddr]_{\mxk[3]{0}}
  \ar@/^1em/@{-->}[drrrr]^{\sigmaul{2}{=1}(\tYxk[2]{1})} \ar@{.>}[dr]^{\mxk[3]{1}} \\
& \Mul[\bth]{1} \ar@{}[drrr] |<<<<{\mbox{\large{$\lrcorner$}}} \ar@{->>}[d] \ar[rrr] &&&
\cprod_{j <k \leq 3} \Yof{\bo}^{d_j d_{\ell}} \ar@{->>}[d] \\
& \Mul[\bth]{0} \ar@{=}[r] &
\cprod_{i<j<k\leq3} \Yof{\bze}^{d_i d_j d_{\ell}}  \ar[rr]^-{\Psul{\bth}{1}}
    && \cprod_{j<\ell \leq 3} \cprod_{i \leq 1} \Yof{\bze}^{d_i d_j d_{\ell}} ~.
}
\noindent where the right vertical is a product of fibrations
\w{\Yof{\bo}\to\Mul[\bo]{0}=\prod_{i \leq 1} \Yof{\bze}^{d_i}} (by Reedy fibrancy of
\w[).]{\Yul{\bt}}

Since \w{\tYxk[\bth]{2}} is homotopy commutative, by Lemma \ref{lhpp} we obtain a
dotted map \w{\mxk[\bth]{1}} (after altering the dashed map \wh that is, our choice for
each \w{d_{j}d_{\ell}} \wwh  within its homotopy class). By Lemma \ref{lraneq} this yields
\w[,]{\Yul[\bth]{1}} still representing \w[.]{\tYxk[\bth]{2}}

It is at stage \w{k=2} in the inner induction that we first encounter a possible
obstruction: we must now choose representatives for \w{d_{\ell}:\bth\to\bt}
\wb{0 \leq \ell \leq 3} in the homotopy class given by \w[.]{\tYxk[\bth]{2}}

As in \wref[,]{usepsi} we know that the target of the forgetful map from
\w{\Mul[\bth]{1}} is the product of the lower left and upper right corners
of \wref[.]{eqtwotr} Thus \w{\Psi=\Psul{\bth}{2}} is a product of two maps:
the first taking each factor \w{\Yof{\bo}^{d_j d_{\ell}}}
\wb{0 \leq j < \ell \leq 3} diagonally to a product
\w[,]{\Yof{\bo}^{d_j d_{\ell}}\times \Yof{\bo}^{d_{\ell-1} d_j}} and the second taking
\w{\Yof{\bze}^{d_i d_j d_{\ell}}} \wb{0 \leq i < j < \ell \leq 3} diagonally to the product
\w[.]{\Yof{\bze}^{d_i d_j d_{\ell}}\times\Yof{\bze}^{d_{i}d_{\ell-1} d_{j}}
\times\Yof{\bze}^{d_{j-1} d_{\ell-1} d_{i}}}

As in \S \ref{dhho}, we now factor $\Psi$ as a trivial cofibration to \w{F^1}
followed by a fibration \w[,]{\Psi'} and pull back the product of the
forgetful maps
$$
\Psul{\bth}{2}:\Mul[\bt]{1}~\to~
\prod_{j \leq 2} \Yof{\bo}^{d_j} \times \prod_{i<j\leq2} \Yof{\bze}^{d_i d_j}
$$
\noindent as in \wref[,]{eqtwtrunc} indexed by the first face maps
\w{d_{\ell}:\bth\to\bt} \wb{0 \leq\ell\leq 3} along \w{\Psi'} to obtain a
``potential mapping diagram'' as in \wref[:]{firstHHO}
$$
\xymatrix@R=39pt@C=14pt{
\Yof{\bth} \ar@/_2em/[ddr]^{\eta_{1}}
\ar@{-->}@/^2em/[drrr]^{\sigmaul{\bth}{2}(\tYxk[\bth]{2})=(d_0,d_1,d_2,d_3)}
    \ar@{-->}@/_4em/[dddrr]^{\varphi}
    \ar@/^1em/@{-->}[drr]^{\kappa} \ar@{.>}[dr]^{\alpha_{2}}
        \ar@/_3em/[dddr]_{\sigmaul{\bth}{<2}(\Yul[\bth]{1})} \\
        & \Pul[\bth]{2} \ar@{}[dr] |<{\mbox{\large{$\lrcorner$}}} \ar@{->>}[d]^{p_{1}} \ar[r] &
        F^3 \ar@{}[dr] |<{\mbox{\large{$\lrcorner$}}} \ar@{->>}[d]^{\mu} \ar@{->>}[r]^-{r'_2}
    & \cprod_{\ell \leq 3} \Yof{\bt} \ar@{->>}[d]^{\prod \mxk[s]{k-1}} \\
        &  \Qul[\bth]{1} \ar@{}[dr] |<{\mbox{\large{$\lrcorner$}}} \ar[d] \ar[r]^{\gamma} &
        {F^2} \ar@{}[dr] |<{\mbox{\large{$\lrcorner$}}} \ar[d]^{q} \ar@{->>}[r]^{s} & \cprod_{\ell \leq 3}
       \Mul[\bt]{1} \ar[d] \\
&  \cprod_{j<\ell\leq 3} \Yof{\bo}^{d_j d_{\ell}} \times\hspace*{-3mm}
         \cprod_{i<j<\ell\leq 3}\Yof{\bze}^{d_i d_j d_{\ell}} \ar[r]_-{\sim} &
           F^1 \ar@{->>}[r]^(0.2){\Psi'}
& \cprod_{\ell \leq 3}( \cprod_{j \leq 2} \Yof{\bo}^{d_j d_{\ell}} \times \hspace*{-2mm}
\cprod_{i<j \leq 2} \hspace*{-1mm} \Yof{\bze}^{d_{i}d_{j}d_{\ell}})
}
$$
\noindent Note that as in \S \ref{dhho}, we may choose \w{F^{1}} to be a product of
free path spaces, so we can think of $\varphi$ as a choice of homotopies between
the various decompositions in \w{\Yul{2}} of maps \w{\bth\to\bze} in $\Delta$.

As the right vertical rectangular pullback has horizontal fibrations, we can apply
Lemma \ref{lhpp} and the fact that the original outermost diagram commutes
up to homotopy (because \w{\tYxk[\bth]{2}} is homotopy commutative)
to deduce that there is a map $\varphi$ in the correct homotopy class, yielding
$\kappa$ as indicated.

The question is whether \w[.]{\mu \kappa \sim \gamma \eta_{1}}
By Corollary \ref{cthpp2}, our secondary operation consists precisely of
those \w{[\theta]} satisfying \w[.]{\theta \sim \mu \circ \kappa}  Thus, the question
 is answered in the affirmative precisely when our secondary operation \w{\lra{\Yul[\bth]{2}}} vanishes.
In that case, by Lemma \ref{lhpp} applied to the upper left square, with $\mu$ a fibration, we can
find \w{\kappa' \sim \kappa} satisfying \w[,]{\mu \circ \kappa'=\gamma \circ \eta_1} so
inducing the dotted \w{\alpha_2} by the pullback property.  We then alter the map
labeled \w{(d_0,d_1,d_2,d_3)} within its homotopy class by instead using \w[,]{r'_2 \circ \kappa'}
which will make the entire diagram now commute, since
\[
\prod \mxk[s]{k-1} \circ (r'_2 \circ \kappa')=s \circ \mu \circ \kappa'=s \circ \gamma \circ \eta_{1}
\]
and \w[.]{q \circ \mu \circ \kappa' = q \circ \gamma \circ \eta_{1} =\iota\circ\sigmaul{\bth}{<2}}
Thus, we obtain a full $3$-truncated $\Delta$-simplicial object \w{\Yul{3}}
(if we wish to proceed further, we take a Reedy fibrant replacement).

If \w{\lra{\Yul[\bth]{2}}} does not vanish, then there is no way to extend this
\w{\Yul{2}} to a full $3$-truncated object.
\end{mysubsection}

\begin{remark}\label{robstr}
  As with any obstruction theory,  when \w{\lra{\Yul[\bth]{2}}} does not vanish, we need
  to backtrack, and see if we can get our obstruction to vanish by modifying previous choices.
  We observe that in special cases, given a truncated $\Delta$-simplicial object, there is a
  formal procedure for adding degeneracies to obtain a full (similarly truncated)
  simplicial object (see, e.g., \cite[\S 6]{BlaHH}).
\end{remark}

%
%
\section{Pointed higher operations}
\label{cgho}

Most familiar examples of higher homotopy operations are pointed, so we now
describe the modifications needed in our general setup when the indexing category
$\cJ$, as well as the model category $\cE$, are pointed (see
\S \ref{cgrms}.B). This will also cover ``hybrid'' cases, where certain
composites in the diagram are required to be \emph{zero} in $\cE$, rather than just
null homotopic.

\begin{lemma}\label{rbottomY}
If \w{\rE} is a pointed model category, \w{\tY:\cJ \to\ho(\rE)} a pointed diagram,
and \w{x\in\Obj\cJ} with \w[,]{|x|>0} then
\begin{enumerate}
\renewcommand{\labelenumi}{(\alph{enumi})~}
\item Any choice of a representative \w{\Yxk[0](g)} of \w{\tY(g)} for
every \w{g \in \tJ^{x}_{0}} yields a lifting of \w{\tY\rest{\Jxk[0]}}
to \w[.]{\Yxk[0]:\Jxk[0]\to\rE}
\item Any pointed Reedy fibrant \w{\Yul{1}:\partial \Jxk[1] \to \rE} as in  \S \ref{sdoubind}(II)
has a pointwise extension to a functor \w{\Yxk[1]:\Jxk[1] \to \rE} which lifts \w[.]{\tYxk{1}}
\end{enumerate}
\end{lemma}

\begin{proof}
For (a), note that if \w[,]{g \in \bJ} \w{\Yof{g}} must be the zero map, but
otherwise any choice of lifting will do, since \w{\Jxk[0]} has
no non-trivial compositions.
For  (b), follow the proof of Lemma \ref{bottomY} with $\tJ$ replacing $\cJ$,
using reduced matching spaces and Definition \ref{dprfib} for the fibrancy.
\end{proof}

We also have the following version of Lemma \ref{lowpiece}:

\begin{lemma}\label{ptlowpiece}
Assuming \w[,]{2 \leq k \leq n < |x|} any pointed functor \w{Y:\cJ_n \to\rE}
with a pointed extension to \w{\smx[k-1]} induces a pullback grid with
natural dashed maps :
\myudiag[\label{eqtlowpiece}]{
\Yof{x} \ar@/^1em/[drrr]^{\rho_{k-1}} \ar@{-->}[dr]_{\beta_{k-1}}
\ar@{-->}@/^1em/[drr]_{\eta_{k-1}} \ar@/_2em/[ddr]_{\rmxk{k-1}} \\
& \rNxk[k-1] \ar@{}[dr] |<<{\mbox{\large{$\lrcorner$}}} \ar[d] \ar[r]^{q_{k-1}} &
\rQxk[k-1] \ar@{}[dr] |<<{\mbox{\large{$\lrcorner$}}} \ar[d] \ar[r] &
\dprod[|s|=k]{\tJ(x,s)} \rMul[s]{k-1} \ar[d]^(0.65){\prod_{\tJ(x,s)} \rforget} \\
& \rMxk[k-1] \ar@{^{(}->}[r]^{\rforget} &
\dprod[|t|<k]{\tJ(x,t)} \Yof{t} \ar[r]^(0.4){\uPsi}
 & \dprod[|s|=k]{\tJ(x,s)}\dprod[|v|<k]{\tJ(s,v)} \Yof{v} ~.
}
\end{lemma}

We then deduce the following analogue of Proposition \ref{startHHO} (with a similar
proof):
\begin{prop}\label{startpHHO}
Assuming \w[,]{2 \leq k \leq n < |x|} any pointed functor
\w{\Yul{k}:\partial \Jxk \to \cE} as in \S \ref{sdoubind} induces maps into a
pullback grid:
\mydiagram[\label{basicpHHO}]{
\Yof{x} \ar@/_2em/[dddr]_{\rmxk{k-1}} \ar@{-->}@/^1em/[drrr]^{\sigmaul{x}{=k}(\tYxk{k})}
\ar@/_1em/[ddr]_(0.7){\beta_{k-1}} \ar@/_1em/[ddrr]^(0.25){\eta_{k-1}}
\ar@{.>}@/^1.5em/[dr]^(0.7){\rmxk{k}}    \ar@{.>}@/^1em/[drr]^(0.7){\alpha_k} \\
& \rMxk \ar@{}[dr] |<{\mbox{\large{$\lrcorner$}}} \ar[d] \ar[r] &
\rPxk \ar@{}[dr] |<{\mbox{\large{$\lrcorner$}}} \ar[d]^{p_{k-1}} \ar[r]^{r_k} & \dprod[|s|=k]{\tJ(x,s)}
\Yof{s} \ar[d]^{\prod \rmxk[s]{k-1}} \\
& \rNxk[k-1] \ar@{}[dr] |<<{\mbox{\large{$\lrcorner$}}} \ar[d] \ar[r]_{q_{k-1}} &
\rQxk[k-1] \ar@{}[dr] |<<{\mbox{\large{$\lrcorner$}}} \ar[d] \ar[r] &
\dprod[|s|=k]{\tJ(x,s)} \rMul[s]{k-1} \ar[d]^(0.6){\prod \rforget} \\
& \rMxk[k-1] \ar@{^{(}->}[r]^{\rforget}   & \dprod[|t|<k]{\tJ(x,t)} \Yof{t}
\ar[r]^(0.4){\uPsi} & \dprod[|s|=k]{\tJ(x,s)}\dprod[|v|<k]{\tJ(s,v)} \Yof{v} ~.
}
\noindent Again, the dashed map only makes the outermost diagram commute up to
homotopy.

Furthermore, the dotted map \w{\rmxk{k}} exists (after altering
\w{\sigmaul{x}{=k}(\tYxk{k})} within its homotopy class) if and only if there is a
dotted map \w{\alpha_k} such that \w{p_{k-1} \alpha_k = \eta_{k-1}} and
\w[.]{r_k \alpha_k \simeq \sigmaul{x}{=k}(\tYxk{k})}
\end{prop}

With this at hand, we may modify Definition \ref{dhho} as follows to obtain
a sequence of obstructions to extending pointed diagrams:

\begin{mysubsection}{Total Pointed Higher Homotopy Operations}
Assume given pointed functors \w[,]{\tYxk{k}:\Jxk\to\ho(\rE)}
\w{\Yxk[k-1]:\Jxk[k-1]\to\rE} and a pointed Reedy fibrant \w{\Yul{k}:\partial\Jxk\to\rE}
as in \S \ref{sdoubind}(II). This means each \w{\rmxk[s]{k-1}:\Yof{s}\to\rMul[s]{k-1}}
is a fibration. Factor \w{\uPsi=\uPsul{x}{k}} (see Lemma \ref{lptmatchpull})
as a weak equivalence followed by a fibration \w[,]{\uPsi'} and pull back the right
column of \wref{basicpHHO} along \w{\uPsi'} to obtain the following pullback grid:
\myvdiag[\label{secondHHO}]{
\Yof{x} \ar@/_2em/[ddr]^{\eta_{k-1}} \ar@/^2em/[drrr]^{\sigmaul{x}{k}(\tYxk{k})}
\ar@{-->}@/_4em/[dddrr]^{\varphi} \ar@/^1em/@{-->}[drr]^{\kappa}
      \ar@{.>}[dr]^{\alpha_k} \ar@/_3em/[dddr]_{\sigmaul{x}{<k}(\Yxk[k-1])} \\
      & \rPxk \ar@{}[dr] |<{\mbox{\large{$\lrcorner$}}} \ar@{->>}[d]^{p_{k-1}} \ar[r] &
          {F^3} \ar@{}[dr] |<{\mbox{\large{$\lrcorner$}}} \ar@{->>}[d]^{\mu} \ar@{->>}[r]^-{r'_k} &
\dprod[|s|=k]{\tJ(x,s)} \Yof{s} \ar@{->>}[d]^{\alpha} \\
&  \rQxk[k-1] \ar@{}[dr] |<{\mbox{\large{$\lrcorner$}}} \ar[d] \ar[r]^{\gamma} &
      {F^2} \ar@{}[dr] |<{\mbox{\large{$\lrcorner$}}} \ar[d]^{q} \ar@{->>}[r]^(0.4){s} &
     \dprod[|s|=k]{\tJ(x,s)} \rMul[s]{k-1} \ar[d]^{\beta} \\
&  \dprod[|t|<k]{\tJ(x,t)} \Yof{t} \ar[r]^{\iota}_{\sim} &
F^1 \ar@{->>}[r]^(0.4){\uPsi'}
& \dprod[|s|=k]{\tJ(x,s)}\dprod[|v|<k]{\tJ(s,v)} \Yof{v}
}
\noindent As in \S \ref{dhho}, Lemma \ref{lhpp} allows us to modify $\varphi$ so
as to obtain a map \w{\kappa:\Yof{x} \to F^3} into the pullback.
\end{mysubsection}

\begin{defn}\label{dphho}
We define the \emph{total pointed higher homotopy operation for $x$} to be the
set \w{\lra{\Yxk[k-1]}} of homotopy classes of maps \w{\theta:\Yof{x} \to F^2}
with \w{\uPsi'\circ q\circ\theta=\beta\circ\alpha\circ\sigmaul{x}{k}} with
  \w[,]{q\circ\theta\sim\varphi} where $\varphi$ is defined to be the composite
$$
\Yof{x}~\stackrel{\sigmaul{x}{<k}}{\longrightarrow}~
\dprod[|t|<k]{\tJ(x,t)}~\Yof[k]{t}~\stackrel{\iota}{\longrightarrow}~F^{1} ~.
$$
\noindent We say \w{\lra{\Yxk[k-1]}} \emph{vanishes at}
\w{\theta:\Yof{x} \to F^{2}}  as above if $\theta$ is homotopic to the composite
$$
\Yof{x}~\stackrel{\eta_{k-1}}{\longrightarrow}~\rQxk[k-1]~\stackrel{\gamma}{\to}~F^2 ~,
$$
and that \w{\lra{\Yxk[k-1]}} \emph{vanishes} if it vanishes at some value $\theta$.
\end{defn}

\begin{remark}\label{rphho}
In many cases of interest we will have \w[,]{\rQxk[k-1]\simeq\ast} in which case
the pointed operation \w{\lra{\Yxk[k-1]}} vanishes at $\theta$
precisely when \w[,]{\theta\sim\ast} as one might expect, so the subset vanishes
precisely when it contains the zero class.
\end{remark}

We have chosen our definitions so as to have the following analogue of
Proposition \ref{obstructk}:

\begin{prop}\label{pobstructk}
Assume given pointed functors \w[,]{\tYxk{k}:\Jxk\to\ho(\rE)}
\w{\Yxk[k-1]:\Jxk[k-1]\to\rE} and a pointed Reedy fibrant \w{\Yul{k}:\partial\Jxk\to\rE}
as in \S \ref{sdoubind}(II) for \w[.]{|x|>n\geq k\geq 2}
Then there exists a further pointed extension to \w{\Yxk:\Jxk \to \rE} if and only
if the total higher homotopy operation \w{\lra{\Yxk[k-1]}} vanishes.
\end{prop}

\begin{proof}
Once again, the definition of \w{\lra{\Yxk[k-1]}} together with Corollary \ref{cthpp2}
  implies that each value $\theta$ is homotopic to \w{\mu \circ\kappa} for some $\kappa$
  with \w{r'_k \circ \kappa = \sigmaul{x}{k}} and
  \w[.]{q \circ\mu\circ\kappa \sim \iota \circ \sigmaul{x}{<k}}
Thus the obstruction vanishes at $\theta$ if and only if there exists such a $\kappa$ with
\w[,]{\mu \circ \kappa \sim \gamma \circ \eta_{k-1}} precisely as in the proof of Proposition \ref{obstructk}.
The upper left pullback square in \wref{secondHHO}
then produces the lift into \w[,]{\rP_x^k} or equivalently, a map
\w[,]{\Yof{x}\to\rM_{x}^{k}} yielding the required pointed extension by
Lemma \ref{lptraneq}.

If \w{\lra{\Yxk[k-1]}} does not vanish, then there is no choice of $\varphi$ for
which such a lift exists, and so there is no pointed extension compatible with
the given choices.
\end{proof}

\begin{remark}\label{rpsepo}
Given pointed functors \w[,]{\tYxk{k}:\Jxk\to\ho(\rE)}
\w{\Yxk[k-1]:\Jxk[k-1]\to\rE} and a pointed Reedy fibrant \w{\Yul{k}:\partial\Jxk\to\rE}
as in \S \ref{sdoubind}(II) for \w[,]{|x|>n\geq k\geq 2} we may define
\emph{separated pointed higher homotopy operations \w{\lra{\Yxk[k-1]}^{j+1}} for $x$}
as in Definition \ref{dseparate}, using a refinement of \wref{secondHHO}
constructed \emph{mutatis mutandis} with products over \w{\cJ(x,s)} replaced
everywhere by products over \w[.]{\tJ(x,s)}

Separation Lemma \ref{lseparate} is stated in sufficient generality to apply here,
too, with Remark \ref{rseparate} modified accordingly, yielding the following
variant of Theorem \ref{unptedThm}:
\end{remark}

\begin{thm}\label{ptedThm}
Assume given pointed functors \w[,]{\tYxk{k}:\Jxk\to\ho(\rE)}
\w{\Yxk[k-1]:\Jxk[k-1]\to\rE} and a pointed Reedy fibrant \w{\Yul{k}:\partial\Jxk\to\rE}
as in \S \ref{sdoubind}(II) for \w[.]{|x|>n\geq k\geq 2}
Then the total pointed higher homotopy operation separates into a sequence of \w{k-1}
pointed operations, and the following are equivalent:
\begin{enumerate}
\item A further extension to \w{\Yxk:\Jxk\to\rE} exists;
\item The total pointed operation \w{\lra{\Yxk[k-1]}} vanishes;
\item The associated sequence \w{\lra{\Yxk[k-1]}^{j+1}} \wb{1\leq j<k} of separated
pointed higher homotopy operations of \S \ref{dseparate} vanish
(so in particular each in turn is defined).
\end{enumerate}
\end{thm}

%
%
\section{Long Toda Brackets and Massey Products}
\label{cltbmp}

We are finally in a position to apply our general theory to the two most familiar examples
of higher order operations: (long) Toda brackets and (higher) Massey products.
Since both are cases of the (pointed) higher operations fully described in
Sections \ref{cgdhho}-\ref{csto} and \ref{cgho}, we thought it would be easier for
the reader to consider two specific examples in detail, briefly indicating
what needs to be done for the higher version.

\supsect{\protect{\ref{cltbmp}}.A}{Right justified Toda brackets}

Since the ordinary Toda bracket (of length $3$) was treated in Section \ref{cctb},
we start with the next case, the Toda bracket of length $4$ (the first example of
a \emph{long Toda bracket} in the sense of \cite{GWalkL}).

Thus, if $\rE$ is a pointed model category, assume given a diagram \w{\tY:\cJ\to\ho\rE} of the form
\mydiagram[\label{eqfourtoda}]{
\Yof{4} \ar[r]^{[k]} & \Yof{3} \ar[r]^{[h]} & \Yof{2} \ar[r]^{[g]} &
\Yof{1} \ar[r]^{[f]} & \Yof{0}
}
\noindent with each adjacent composite null-homotopic: that is, a chain complex of
length $4$ in \w[,]{\ho\rE} as in Example \ref{ptchain} (compare
\wref[).]{eqtodadiag}  Without loss of generality, we can assume all objects involved are
both cofibrant and fibrant.

Applying the double induction procedure of \S \ref{sdoubind}, we see that we
must deal with chain complexes of length \w[,]{n\leq 4} as follows:
\begin{enumerate}
\renewcommand{\labelenumi}{(\alph{enumi})~}
\item When \w[,]{n=0} we have no inner induction, and making the result Reedy fibrant
consists of factoring the representative to produce a fibration
\w{f:\Yof{1}\epic\Yof{0}} in the specified class \w[.]{[f]}
\item When \w[,]{n=1} note that \w{\tJ(x,t)} is empty if \w[,]{|x|-|t| > 1}
 for this pointed indexing category, so as a consequence \w{\rMxk= \ast} if
\w[.]{|x|-|k|>1}  Thus \w[,]{\rMul[2]{0}=\ast} so \w{\rMul[2]{1}} is simply the
fiber of $f$.  Since \w[,]{[f]\circ[g]=\ast} by Lemma \ref{lhpp} we can choose a
representative $g$ for \w{[g]} which factors as a fibration
\w{\Yof{2}\epic\rMul[2]{1}=\fib(f)} followed by the inclusion
\w[.]{\fib(f)\hra\Yof{1}}
\item When \w[,]{n=2} again \w[,]{\rMul[3]{0}=\ast=\rMul[3]{1}}
while the case \w{k=2} is just that of our (length $3$) Toda bracket
\w[.]{\lrau{f,g}{,h}}

In this case, the indexing set for products in the right column of \wref{secondHHO}
is the singleton \w[,]{\tJ(3,2)} while the forgetful map in the bottom row of
\wref{basicpHHO} is the identity of the zero object, with $\Psi$ the zero map.

Factoring $\Psi$ as a trivial cofibration $\iota$ followed by a fibration
\w[,]{\Psi'} as in the bottom row of \wref[,]{secondHHO} and pulling back the right
column yields the diagram:
\mydiagram[\label{eqthrtoda}]{
\Yof{3} \ar@{.>}[dr] \ar@{.>}@/^1em/[drrr] \ar@{->}@/^1.5em/[drrrrr]^{h}
\ar@{-->}@/_2.4em/[ddrrr]_(0.3){\theta} \\
& \rPul[3]{2} \ar@{}[dr] |<<{\mbox{\large{$\lrcorner$}}} \ar[d] \ar[rr] &&
F^3 \ar@{}[dr] |<<{\mbox{\large{$\lrcorner$}}} \ar[d] \ar@{->>}[rr] &&
\Yof{2} \ar@{->>}[d] \ar@{->}@/^1.5em/[dd]^(0.35){g} \\
& {\ast} \ar[rr] && F^2 \ar@{}[dr] |<<{\mbox{\large{$\lrcorner$}}} \ar[d] \ar@{->>}[rr] &&
\rMul[2]{1} \ar[d]_{\rforget} \ar@{.>>}[rr] && {\ast} \ar@{.>}[d] \\
& {\ast} \ar[rr]_{\sim}^{\iota} && F^1 \ar@{->>}[rr]^{\Psi'} &&
\Yof{1} \ar@{.>>}[rr]_{f} && \Yof{0} ~.
}
\noindent Thus \w{F^{1}} is a model for the reduced path space on \w[,]{\Yof{1}}
with \w{\Psi'} the path fibration.
However, since $f$ was chosen above to be a fibration, the composite
\w{F^{1}\to\Yof{0}} is a fibration, too, with \w{F^{1}} contractible, so we see that
\w[,]{F^{2}} being the pullback of the dotted rectangle, is a model for
the loop space \w[,]{\Omega \Yof{0}} which we denote by \w[.]{\Omega' \Yof{0}}
Similarly, \w{\rMul[2]{1}} is a model for \w[.]{\fib(f)}

Our total secondary pointed homotopy operation \w{\lra{\Yul[\bth]{1}}}
(cf.\ \S \ref{dphho}) is thus a set of maps \w[,]{\theta:\Yof{3} \to \Omega' \Yof{0}}
and it vanishes when this set contains the zero map (cf.\ Remark \ref{rphho}).
This is our usual Toda bracket \w[,]{\lrau{f,g}{,h}} described in the language
of Section \ref{cgho}.
\item In order for our four-fold Toda bracket \w{\lrau{f,g,h}{,k}} (denoted by
\w{\lra{\Yul[\bfo]{2}}} above) to be defined, \w{\lra{\Yul[\bth]{1}}}
must vanish. This allows us to choose a pointed extension \w{\Yul{3}:\Jul{3}\to\rE}
of \w{\Yul{2}} which realizes \w[.]{\tY\rest{\Jul{3}}} The fact that the
diagram \w{\Yul{3}} has realized $\tY$ through filtration degree $3$ means that
each of the maps $g$ and $h$ factors through the fiber of the previous one, as in
the following solid commutative diagram:
\mydiagram[\label{eqfortoda}]{
\Yof{4} \ar@{.>}[r]^-{k_1} \ar[dr]_{k} & \fib(h_1) \ar[d] \ar[r] & \ast \ar[d] \\
& \Yof{3} \ar@{->>}[r]^-{h_{1}} \ar[dr]_{h }& \fib(g_1) \ar@{->>}[r] \ar[d] &
         \ast \ar[d] \\
& & \Yof{2} \ar@{->>}[r]^-{g_{1}} \ar[dr]_{g} & \fib(f) \ar@{->>}[r] \ar[d]^{g_2} &
          \ast \ar[d] \\
& & & \Yof{1} \ar@{->>}[r]_{f} & \Yof{0} ~.
}
\noindent Making \w{\Yul{3}} pointed Reedy fibrant (\S \ref{dprfib}) just means
ensuring that the maps \w{h_{1}} and \w{g_{1}} are fibrations.
\item At stage \w{n=3} in the outer induction, we attempt to find the dotted lift
\w{k_{1}} in \wref[,]{eqfortoda} after having chosen a suitable representative $h$ for
the given homotopy class \w[,]{[h]} which is possible by the vanishing of the previous obstruction.

Again we have \w[,]{\rMul[4]{0}=\ast} \w[,]{\rMul[4]{1}=\ast} and
\w[,]{\rMul[4]{2}=\ast=\rQul[4]{2}=\rNul[4]{2}} so the only interesting case is
\w{k=3} in the inner induction.

The separation grid of Lemma \ref{lseparate} then takes the form:
\mydiagram[\label{eqsepfourtoda}]{
\Yof{4} \ar@{.>}[dr]^{k_{1}} \ar@{-->}@/^1em/[drr]
\ar@{-->}@/^1.5em/[drrr]^(0.6){\kappa} \ar@{->}@/^2em/[drrrr]^{k} \ar@/_1em/[ddr] \\
& \rPul[4]{3} \ar@{}[dr] |<{\mbox{\large{$\lrcorner$}}} \ar[d] \ar[r] &
F_{4,3}^{2,4} \ar@{}[dr] |<{\mbox{\large{$\lrcorner$}}} \ar[d] \ar@{->>}[r] &
F_{4,3}^{1,4} \ar@{}[dr] |<{\mbox{\large{$\lrcorner$}}} \ar[d] \ar@{->>}[r] & \Yof{3} \ar@{->>}[d]_{h_{1}}
\ar@{->}@/^2em/[ddd]^{h}\\
& {\ast} \ar[dr]_{\sim} \ar[r] & F_{4,3}^{2,3} \ar@{}[dr] |<{\mbox{\large{$\lrcorner$}}} \ar[d] \ar@{->>}[r] &
F_{4,3}^{1,3} \ar@{}[dr] |<{\mbox{\large{$\lrcorner$}}} \ar[d] \ar@{->>}[r] &
\rMul[3]{2} \ar[d] \ar@{.>>}[rr] && {\ast} \ar@{.>}[d]  & {\ast} \ar@{.>}[d]^{\sim} \\
& & F_{4,3}^{2,2} \ar@{->>}[r] &
F_{4,3}^{1,2} \ar@{}[dr] |<{\mbox{\large{$\lrcorner$}}} \ar@{->>}[d] \ar@{->>}[r] &
F_{3,2}^{1,3} \ar@{->>}[d] \ar@{.>>}[rr]
&& F_{3,2}^{1,2} \ar@{}[dr] |<{\mbox{\large{$\lrcorner$}}} \ar@{.>>}[d] \ar@{.>}[r] &
F_{3,2}^{1,1} \ar@{.>>}[d] \\
& {\ast} \ar[r] & {\ast} \ar[r]^{\sim} & F_{4,3}^{1,1} \ar@{->>}[r] &
\Yof{2} \ar@{.>>}[rr]^{g_{1}} \ar@{->}@/_1.3em/[rrr]_{g} &&
\rMul[2]{1} \ar@{.>}[d]^{\ast} \ar@{.>}[r]^(0.5){g_{2}} &  \Yof{1} \ar@{.>>}[d]^(0.4){f} \\
&&&&&& \ast \ar@{.>}[r] & \Yof{0}
}
\noindent where we have extended the pullback grid downwards, and to the right,
to show how it was constructed from the previous case (diagram \wref[)]{eqthrtoda}
using Lemma \ref{lprodgrid}. We have also indicated how (representatives of)
the maps of \wref{eqfortoda} fit in.

As in Step (c) above, we can identify \w{F_{3,2}^{1,2}} as a model for
\w[,]{\Omega \Yof{0}}  and \w{\rMul[2]{1}} as a model for \w[.]{\fib(f)}
Similarly, \w{F_{4,3}^{1,2}} is a model for
\w[,]{\Omega \Yof{1}} using the vertical fibrations in the
rectangle with diagonal corners \w{F_{4,3}^{1,2}} and \w[.]{\Yof{1}}
Likewise \w{F_{4,3}^{1,3}} a model for \w{\Omega \rMul[2]{1}} (using horizontal
fibrations in the larger square beneath it), and \w{F_{4,3}^{2,3}} is a model for
\w{\Omega^2 \Yof{0}} (now using the rectangle with diagonal corners \w{F_{4,3}^{2,3}}
and \w[,]{F_{3,2}^{1,2}} along with the previous identification of the latter).
Similarly, \w{\rMul[3]{2}} is a model for \w{\fib(g_{1})} of \wref[,]{eqfortoda}
while \w{\rPul[4]{3}} is \w{\fib(h_{1})} (which is also the homotopy fiber).
See \wref{eqsepfrtoda} below for the full identification.

Therefore, the final obstruction to having a dotted lift \w{k_{1}} in \wref{eqfortoda}
(or \wref[)]{eqsepfourtoda} is the composite \w[.]{k\circ h_{1}}
\end{enumerate}

Note that there are no factors of type \w{G_{i,j}^{k,\ell}} as in \wref{eqsepgrid}
here, since we can always choose the zero map as our factorization of the
zero map between zero objects.

\begin{remark}\label{rasepop}
Our total pointed tertiary homotopy operation \w{\lra{\Yul[\bfo]{2}}} is a set
of homotopy classes \w[.]{\theta:\Yof{4} \to \Omega \rMul[2]{1}}
However, using Lemma \ref{lseparate}, we can replace it by two
separated higher homotopy operations for $\bfo$, in the sense of \S \ref{dseparate}:
\begin{enumerate}
\renewcommand{\labelenumii}{\roman{enumii}.~}
\item The second order operation
\w[.]{\lra{\Yul[\bfo]{2}}^{2}\subseteq[\Yof{4},\,\Omega \Yof{1}]}
\item If \w{\lra{\Yul[\bfo]{2}}^{2}} vanishes, the third order operation
\w{\lra{\Yul[\bfo]{2}}^{3}\subseteq[\Yof{4},\,\Omega^{2}\Yof{0}]} is defined, and
serves as the final obstruction to lifting $\tY$. By definition, this is
our \emph{four-fold Toda bracket} \w[.]{\lrau{f,g,h}{,k}}
\end{enumerate}
\end{remark}

\begin{lemma}\label{lthreetoda}
Given a pointed Reedy fibrant diagram \w{\Yul{3}} realizing \wref{eqfourtoda}
through filtration $3$, the associated second order separated higher homotopy operation
\w{\lra{\Yul[\bfo]{2}}^{2}} is our usual Toda bracket \w[.]{\lrau{g,h}{,k}}
\end{lemma}

\begin{proof}
Note that \w{F_{3,2}^{1,3}} is a model for the homotopy fiber of
\w{g:\Yof{2} \to \Yof{1}} (which is not-itself a fibration).
Thus, the rectangle with corners \w{F_{3,2}^{1,3}} and \w{\Yof{1}} in
\wref{eqsepfourtoda} is a homotopy invariant version of the
rectangle with corners \w{F^{2}} and \w{\Yof{0}} in \wref[,]{eqthrtoda} used to define
our Toda bracket in Step (c) above \wh this time, applied to the
left $3$ maps in \wref[.]{eqfourtoda} The map corresponding to $\theta$ in
\wref{eqthrtoda} \wh the value of the Toda bracket \wh is the map
\w{\Yof{4}\to F_{4,3}^{1,2}} obtained by composing $\kappa$ with the vertical maps
\w[,]{F_{4,3}^{1,4}\to F_{4,3}^{1,2}} which is indeed the definition of the
value of \w{\lra{\Yul[\bfo]{2}}^{2}} associated to our choices (see
Definition \ref{dseparate}).
\end{proof}

\begin{aside}\label{asepopn}
Note that \emph{if} the dotted forgetful map \w{\rMul[2]{1}\to\Yof{1}}
in \wref{eqsepfourtoda} were a fibration, the horizontal dotted map above it would
be a fibration, too, so right properness would imply that the vertical map
\w{\rMul[3]{2}\to F_{3,2}^{1,3}} would be a weak equivalence.
\end{aside}

\begin{mysubsection}{Length $n$ Toda brackets}
\label{slntb}

The general procedure described in Section \ref{cgho} tells us what needs to be done
for Toda diagrams (chain complexes $\tY$ in \w[):]{\ho\rE}
\mydiagram[\label{eqntoda}]{
\Yof{n} \ar[r]^-{[f_{n}]} & \Yof{n-1} \ar[r]^-{[f_{n-1}]} &\dotsc\ar[r] &
\Yof{3} \ar[r]^{[f_{3}]} & \Yof{2} \ar[r]^{[f_{2}]} & \Yof{1} \ar[r]^{[f_{1}]} & \Yof{0}
}
\noindent of arbitrary length $n$. We sketch the main features of the
general construction, already discernible in the case \w{n=4} described above:

In the double induction of \S \ref{sdoubind}, we can concentrate on the last
stage \wh assuming the vanishing of shorter brackets on the right, which guarantees
the existence of a solid diagram
\mydiagram[\label{eqlenntoda}]{
\Yof{n} \ar@{.>}[r]^-{g_{n}} \ar[dr]_{f_n} &
\fib(g_{n-1}) \ar@{}[dr] |<{\mbox{\large{$\lrcorner$}}} \ar[d]_{g^2_{n}}  \ar[r] & \ast \ar[d] \\
& \Yof{n-1} \ar@{->>}[r]^-{g_{n-1}} \ar@{.}[dr] &
\fib(g_{n-2}) \ar@{}[dr] |<{\mbox{\large{$\lrcorner$}}} \ar@{->>}[r] \ar@{.}[d] & \ast \ar@{.}[d] \\
& & \Yof{3} \ar@{->>}[r]^-{g_{3}} \ar[dr]_{f_{3}} &
\fib(g_{2}) \ar@{}[dr] |<{\mbox{\large{$\lrcorner$}}} \ar@{->>}[r] \ar[d]_{g^2_{3}} & \ast \ar[d] \\
& & & \Yof{2} \ar@{->>}[r]^-{g_{2}} \ar[dr]_{f_{2}} &
\fib(f_1) \ar@{}[dr] |<{\mbox{\large{$\lrcorner$}}} \ar@{->>}[r] \ar[d]_{g^2_{2}} & \ast \ar[d] \\
& & & & \Yof{1} \ar@{->>}[r]_{f_1} & \Yof{0}
}
\noindent analogous to \wref[;]{eqfortoda} our length $n$ Toda bracket,
\w[,]{\lrau{f_1,f_2,\dots f_{n-1}}{,f_n}}
 will be the
final obstruction to finding the dotted map \w{g_{n}} in \wref[,]{eqlenntoda}
perhaps after altering \w{f_{n}} within its homotopy class.

The existence of the fibrations \w{g_{k}} for \w[,]{2\leq k<n} and the fact that
\w{f_{1}} is a fibration, mean that we have a lifting \w{\Yul{n-1}:\Jul{n-1}\to\rE}
of \w[,]{\tY\rest{\Jul{n-1}}} which we have made pointed Reedy fibrant. The
underlining in the notation represents our intention to leave that portion fixed.

The construction of the separation grid for \w{\Yul{n-1}} (\S \ref{dsepdiag})
greatly simplifies, in this case, as we see in comparing \wref{eqthrtoda} to
\wref[:]{eqsepfourtoda} at each step, one writes the previous
separation grid vertically (instead of horizontally) on the right
(after changing the previously chosen \w{g_{n-1}} into a fibration, thus altering
\w{\Yof{n-1}} up to homotopy). We then factor the zero map $\Psi$ and pull back the
leftmost existing column to form a new column to its left.  Factoring the zero map
from \w{\rQxk[k-1]} to the second place from the bottom in this new column and
again pulling back, we note that the intermediate object produced by this
factorization is a reduced path object, so by induction the entry immediately above
it is a loop object (being the pullback over a fibration with upper right and lower left corners
contractible \wh one because it is the reduced path object, and the other by induction).
Moreover, the number of loops increases as we move up and to the left
(see Lemma \ref{loopstair}).

Repeat this step until the new column involves just two maps (so the second
object from the bottom is at the same height as the product of the objects
\w{\rMul[s]{k-1}} on the right). The pullback in the upper left corner is now the
actual fiber of \w[.]{g_{n-1}} To illustrate, we reproduce diagram \wref{eqsepfourtoda}
with the pieces identified up to homotopy:
\myodiag[\label{eqsepfrtoda}]{
\Yof{4} \ar@{.>}[dr]^{k_{1}} \ar@{-->}@/^1em/[drr]
\ar@{-->}@/^1.5em/[drrr] \ar@{->}@/^2em/[drrrr]^{k} \ar@/_1em/[ddr] \\
& \fib(h_{1}) \ar@{}[dr] |<{\mbox{\large{$\lrcorner$}}} \ar[d] \ar[r] &
F_{4,3}^{2,4} \ar@{}[dr] |<{\mbox{\large{$\lrcorner$}}} \ar[d] \ar@{->>}[r] &
F_{4,3}^{1,4} \ar@{}[dr] |<{\mbox{\large{$\lrcorner$}}} \ar[d] \ar@{->>}[r] & \Yof{3} \ar@{->>}[d]_{h_{1}}
\ar@{->}@/^2em/[ddd]^{h}\\
& {\ast} \ar[dr]_{\sim} \ar[r] &
*+[F]{\Omega^{2}\Yof{0}} \ar@{}[drr] |<<<{\mbox{\large{$\lrcorner$}}} \ar[d] \ar@{->>}[r] &
\Omega\fib(f) \ar@{}[dr] |<<{\mbox{\large{$\lrcorner$}}} \ar[d] \ar@{->>}[r] &
\fib(g_{1}) \ar[d] \ar@{.>}[rr] & {\ast} \ar@{.>}[d]  & {\ast} \ar@{.>}[d]^{\sim} \\
& & P\Omega\Yof{1} \ar@{->>}[r] &
*+[F]{\Omega\Yof{1}} \ar@{}[dr] |<<<{\mbox{\large{$\lrcorner$}}} \ar@{->>}[d] \ar@{->>}[r] &
\fib(g) \ar@{->>}[d] \ar@{.>>}[r]
    & \Omega\Yof{0} \ar@{.>>}[d] \ar@{.>}[r] & P\Yof{1} \ar@{.>>}[d] \\
&& \ast \ar[r]^{\sim} & P\Yof{2}  \ar@{->>}[r] &
\Yof{2} \ar@{.>>}[r]^{g_{1}} \ar@{->}@/_1.3em/[rr]_{g} &
\fib(f) \ar@{.>}[r]^(0.5){g_{2}} &  \Yof{1}
}
\noindent Note that while not all the pullbacks in the grid can be easily identified,
the targets of the separated operations (boxed) are iterated loop spaces on
the original objects of \wref[,]{eqntoda} as one would expect for long Toda brackets.
This last obstruction, consisting of a subset of the homotopy classes of maps into the
top left iterated loop space, then represents our length $n$ Toda bracket,
\w[,]{\lrau{f_1,f_2,\dots f_{n-1}}{,f_n}} with the lower separated higher homotopy operations corresponding
to the vanishing of the lower obstructions necessary in order to define it (together with those already
assumed to vanish in order to build the current commuting diagram).
\end{mysubsection}

\supsect{\protect{\ref{cltbmp}}.B}{Massey Products as a Hybrid Case}

The classical Massey product (cf.\ \cite{MassN}) is defined for three cohomology
classes of the same space $X$  \w{[\alpha],[\beta],[\gamma]\in H^{\ast}(X;R)} for some ring $R$,
equipped with null homotopies \w{F:\mu(\alpha,\beta) \sim 0} and
\w{G:\mu(\beta,\gamma) \sim 0} for the two products. Like a Toda bracket,
the Massey product serves as the obstruction to simultaneously making both products
strictly zero (see \cite[\S 4]{BBGondH}).

This situation may be described by the pointed indexing category $\cJ$:
\myvdiag[\label{eqmasseyd}]{
& & g \ar[d] \ar[ddl] \ar[ddr] \ar@{-->}[ddrr] \ar@{-->}[ddll] \ar@/^14em/@{-->}[ddd]\\
& & f \ar[dll] \ar[dl] \ar[dr] \ar[drr] \ar@{-->}@/_0.2em/[dd] \ar@/^0.2em/[dd]\\
b \ar@{-->}[drr] & c \ar[dr] \ar[l] && d \ar[dl] \ar[r] & e \ar@{-->}[dll]  \\
& & a
}
\noindent Here the dashed maps are in $\bJ$ and the others are in $\tJ$.
The inner diamond commutes (with the solid composite) and the outer diamond commutes
(with the dashed composite).

The corresponding pointed diagram \w{\tY:\cJ\to\ho\Ta} has products of
Eilenberg-Mac~Lane spaces \w{K_{i}:=K(R,i)} in all but the top slot:
\mywdiag[\label{eqmassey}]{
& & \Yof{g} \ar[d]|{(\alpha,\beta,\gamma)} \ar[ddl]
  \ar[ddr] \ar@{-->}[ddrr]^{\mu(\alpha,\beta)} \ar@{-->}[ddll]_{\mu(\beta,\gamma)} \\
  & & K_r \times K_s \times K_t \ar[dll] \ar[dl]^{(\pi_{1},\mu)} \ar[dr]_{(\mu,\pi_{2})}
\ar[drr] \ar@{-->}@/_0.2em/[dd] \ar@/^0.2em/[dd]^{\mu}\\
\ast\times K_{s+t} \ar@{-->}@/_0.9em/[drr] &
K_r \times K_{s+t} \ar[dr]^{\mu} \ar[l]^{\pi_{2}} &&
K_{r+s} \times K_t \ar[dl]_{\mu} \ar[r]_{\pi_{1}} &
K_{r+s}\times\ast \ar@{-->}@/^0.9em/[dll] \\
& & K_{r+s+t}
}
\noindent where the central diamond represents associativity of the cup product maps
$\mu$; \w{\pi_{1}} and \w{\pi_{2}} are the two projections; and we have omitted the
zero map from top to bottom that appears in \wref{eqmasseyd} in the interest of clarity.

Choose a strictly associative model of the Eilenberg-Mac~Lane $\Omega$-spectrum in
question (cf.\ \cite{RobO}), with strictly pointed multiplication, so in particular at each level
\w{K_{r}} is a simplicial  (or topological) abelian group. We can then make all of
\wref{eqmassey} below \w{\Yof{g}} (involving only the cup product maps) strictly
commutative. Our Massey product will be the total pointed higher homotopy operation
\w{\lra{\Yul[g]{1}}} (for \w[).]{n=k=2}

From \S \ref{dpoinmat} we see that if we let
\w[,]{\bK:=K_{r}\times K_{s+t}\times K_{r+s}\times K_{t}} then
\w{\rMul[f]{1}} is the pullback of the two multiplication maps
\w[,]{K_{r}\times K_{s+t}\to K_{r+s+t}\leftarrow K_{r+s}\times K_{t}}
with a natural inclusion (forgetful map)  \w[.]{i_{1}:\rMul[f]{1}\to\bK} The pullback grid of
\wref{secondHHO} then takes the form:
\myqdiag[\label{eqmasssHHO}]{
\rPul[g]{2} \ar@{}[dr] |<<<{\mbox{\large{$\lrcorner$}}} \ar@{->>}[d]_{p_{k-1}} \ar[r] &
F^3 \ar@{}[dr] |<<<{\mbox{\large{$\lrcorner$}}} \ar@{->>}[d] \ar@{->>}[r]^-{r'_k} &
{K_{r}\times K_{s}\times K_{t}} \ar@{->>}[d] \\
\rQul[g]{1} \ar@{}[dr] |<<<{\mbox{\large{$\lrcorner$}}} \ar[d] \ar[r]^{\gamma} &
     {F^2} \ar@{}[dr] |<<<{\mbox{\large{$\lrcorner$}}} \ar[d]^{q} \ar@{->>}[r]^(0.4){s} &
    \rMul[f]{1} \ar[d]^{(\pi_{2}i_{1},\pi_{3}i_{1},i_{1},\mu i_{1})} \\
\bK \ar[r]^(0.17){\sim} & PK_{r+s}\times PK_{s+t}\times\bK\times PK_{r+s+t}
\ar@{->>}[r]^(0.53){\uPsi'}
& K_{r+s}\times K_{s+t}\times\bK\times K_{r+s+t}
}
\noindent Thus a point in \w{F^{2}} is given by
\w{(U,V,x,u,v,z,W)\in PK_{r+s}\times PK_{s+t}\times\rMul[f]{1}\times PK_{r+s+t}}
with \w[,]{U:u\sim\ast} \w[,]{V:v\sim\ast} and \w[.]{W:xu=vz\sim\ast}
We thus have a natural map \w{\lambda:F^{2}\to\Omega K_{r+s+t}\times\Omega K_{r+s+t}}
sending \w{(U,V,x,u,v,z,W)} to \w[.]{(xU-W,Vz-W)} Postcomposition with the difference map
\w{d:\Omega K_{r+s+t}\times\Omega K_{r+s+t}\to\Omega K_{r+s+t}} yields \w[.]{(xU-Vz)}

Now \w{\Yof{g}} maps into the top right corner of \wref{eqmasssHHO} by
(a lift of) \w[,]{(\alpha,\beta,\gamma)} and thereby on to \w[,]{\rMul[f]{1}}
and into the bottom middle term by
$$
\varphi~:=~\lra{F,\,G,\,\alpha,\,\mu(\beta,\gamma),\,\mu(\alpha,\beta),\,\gamma,\,L}~,
$$
\noindent with $L$ some nullhomotopy of \w[.]{\mu(\alpha,\beta,\gamma)} Together these two maps
induce the map \w{\theta:\Yof{g}\to F^{2}} of \S \ref{dphho}.

Postcomposing $\theta$ with \w{d\circ\lambda} gives the usual Massey product
$$
\lra{\alpha,\beta,\gamma}\in[\Yof{g},\,\Omega K_{r+s+t}]=H^{r+s+t-1}(\Yof{g};R) ~.
$$
The two factors of \w{\lambda\circ\theta} merely give the usual indeterminacy for the Massey product,
as we can see by choosing \w{L:=\mu(F,\gamma)} or \w[.]{L:=\mu(\alpha,G)}

\begin{remark}\label{rmassey}
    An alternative definition of the usual (higher) Massey products, more in line with that given for the Toda
    bracket, appears in \cite[\S 4.1]{BBGondH}.
\end{remark}

%
%
\section{Fully reduced diagrams}
\label{cfrd}

Ultimately, we would like to develop an ``algebra of higher order operations,''
along the lines of Toda's original juggling lemmas (see \cite[\S 1]{TodC}).
As a first step in this direction, we consider a special type of pointed diagram,
which most closely resembles the long Toda diagram of \wref[.]{eqntoda}

The most useful property of the separated higher operations associated to Toda diagrams
is that we can often identify their targets \w{F_{x,k}^{j,j+1}} as loop spaces (as we saw
in \wref[).]{eqsepfrtoda}

It turns out the property of the pointed indexing category $\cJ$ needed for this
to happen is the following:

\begin{defn}\label{dfred}
A pointed indexing category $\cJ$ as in \S \ref{dpindex} is called
\emph{fully reduced} if any morphism decreasing degree by at least $2$ lies in $\bJ$.
\end{defn}

\begin{remark} \label{ptfullred}
If $\cJ$ is fully reduced, for \w{|x| \geq k+1} we have
\w{\prod_{\tJ(x,t), |t| < k} \Yof{t}=\ast} and so \w{\rMxk[k-1]=\ast} (cf.\ \S \ref{dpoinmat})
as well.  We deduce that \w[,]{\rNxk[k-1]=\ast=\rQxk[k-1]} too (cf.\
\wref[),]{eqtlowpiece} since both are fibers of a product of monomorphisms, by Lemma
\ref{ptlowpiece} (under mild assumptions on \w[).]{\rE}

Furthermore, the map \w{\rforget} of \S \ref{dpoinmat} factors through
\w[,]{\prod_{\tJ(s,t), |t|=|s|-1}\,\Yof{t}} so no factors of type \w{G_{x,k}^{k+1,j}}
(cf.\ \wref[)]{eqfgpull} are needed when constructing the separation grid
\wref[.]{eqsepgrid} This also implies that \w{F_{x,k}^{j,j}} is contractible for
\w[,]{j < k} which is the key ingredient for identifying the targets of
the separated operations as loop spaces.
\end{remark}

Our key decomposition result is the following.

\begin{lemma} \label{loopstair}
If $\cJ$ is a fully reduced pointed indexing category and \w[,]{n \geq k \geq j \geq 2}
 we have:
$$
F_{x,k}^{j-1,j}
\sim  \prod_{\begin{aligned} (&f_{k-j}, \dots, f_k) \\
f_{k-j} \circ &\dots \circ f_k:x \to v \end{aligned}} \Omega^{j-1} \Yof{v}
$$
\noindent in \wref[,]{eqsepgrid} where each \w{f_i} is a non-identity map in $\tJ$,
with target of degree $i$.
\end{lemma}

\begin{proof}
We prove this by induction on $k$ (for fixed $n$ and $x$), as in Lemma \ref{lprodgrid}.
In each case, we combine two pullbacks over fibrations, one of which has fiber
identified at an earlier stage, with two corners contractible; the upper left corner
(source) is then homotopy equivalent to the loop space on the lower right corner,
(see Step (e) of \S \ref{cltbmp}.A).

For \w[,]{2 = j < k} we use the basic pullback rectangle
\mysdiag{
F^{1,3}_{s,2} \ar@{}[drr] |<{\mbox{\large{$\lrcorner$}}} \ar@{->>}[d] \ar@{->>}[rr] &&
\dprod[|u|=1]{\tJ(s,u)} \Yof{u} \ar@{->>}[d] \\
F^{1,2}_{s,2} \ar@{}[drr] |<{\mbox{\large{$\lrcorner$}}} \ar[d] \ar@{->>}[rr] &&
\dprod[|u|=1]{\tJ(s,u)} \rMul[u]{0} \ar[d] \\
F^{1,1}_{s,2} \ar@{->>}[rr] && \dprod[|u|=1]{\tJ(s,u)}
\dprod[|v|=0]{\tJ(u,v)} \Yof{v}
}
\noindent to construct the pullback rectangle
\mysdiag{
F^{1,2}_{x,3} \ar@{}[dr] |<{\mbox{\large{$\lrcorner$}}} \ar@{->>}[d] \ar@{->>}[r] &
\dprod[|s|=3]{\tJ(x,s)} F^{1,3}_{s,2} \ar@{}[dr] |<{\mbox{\large{$\lrcorner$}}} \ar[r] \ar@{->>}[d]
    & \dprod[|s|=3]{\tJ(x,s)} F^{1,1}_{s,2} \ar@{->>}[d] \\
F^{1,1}_{x,3} \ar@{->>}[r] & \dprod[|s|=3]{\tJ(x,s)} \Yof{u} \ar[r] &
\dprod[|s|=3]{\tJ(x,s)}\dprod[|u|=1]{\tJ(s,u)} \dprod[|v|=0]{\tJ(u,v)} \Yof{v}
}
\noindent where the vertical maps are fibrations, and both
\w{\dprod[|s|=3]{\tJ(x,s)} F^{1,1}_{s,2}} and \w{F^{1,1}_{x,3}} contractible, as in Remark \ref{ptfullred}.

For \w[,]{2 < j < k} we similarly use the pullback rectangle
\mysdiag{
  F_{x,k}^{j-1,j} \ar@{}[dr] |<{\mbox{\large{$\lrcorner$}}} \ar@{->>}[d] \ar@{->>}[r] &
  \dprod[|s|=k]{\tJ(x,s)} F_{s,k-1}^{j-1,k} \ar@{}[dr] |<{\mbox{\large{$\lrcorner$}}}
    \ar@{->>}[d] \ar[r] & \dprod[|s|=k]{\tJ(x,s)} F_{s,k-1}^{j-1,j-1} \ar@{->>}[d] \\
F_{x,k}^{j-1,j-1} \ar@{->>}[r] & \dprod[|s|=k]{\tJ(x,s)} F_{s,k-1}^{j-2,k}
    \ar[r] & \dprod[|s|=k]{\tJ(x,s)} F_{s,k-1}^{j-2,j-1}
}
\noindent in which the vertical maps are fibrations,
together with the fact that \w{F_{x,k}^{j-1,j-1}} and each
\w{F_{s,k-1}^{j-1,j-1}} are contractible, to prove the claim by induction on $j$
(since loops commute with products).

For \w[,]{2 \leq j=k} recall that when \w{|s|=2} the first non-trivial case (with
\w[)]{k-1=1} involves the first pullback diagram
\mydiagram{
\rMul[s]{1} \ar@{}[drr] |<{\mbox{\large{$\lrcorner$}}} \ar@{->>}[d] \ar[rr]^(0.4){\rforget} &&
\dprod[|u|=1]{\tJ(s,u)} \Yof{u} \ar@{->>}[d] \\
\ast \ar[rr] &&  \dprod[|u|=1]{\tJ(s,u)}\ \dprod[|v|=0]{\tJ(u,v)} \Yof{v}
}
\noindent For \w{2<j=k} we have the second pullback diagram
\mydiagram{
  \rMul[s]{k-1} \ar@{}[dr] |<{\mbox{\large{$\lrcorner$}}} \ar@{->>}[d] \ar[r] &
  \rPul[s]{k-1} \ar@{}[dr] |<{\mbox{\large{$\lrcorner$}}} \ar@{->>}[d] \ar[r] &
F_{s,k-1}^{k-2,k} \ar@{->>}[d] \\
\ast=\rNul[s]{k-1} \ar[r] & \ast=\rQul[s]{k-1} \ar[r] &  F_{s,k-1}^{k-2,k-1}
}
\noindent and combining (products of) either type into
\mydiagram{
  F_{x,k}^{k-1,k} \ar@{}[dr] |<{\mbox{\large{$\lrcorner$}}} \ar@{->>}[r] \ar[d] &
  \dprod[|s|=k]{\tJ(x,s)} \rMul[s]{k-1} \ar@{}[dr] |<{\mbox{\large{$\lrcorner$}}}
  \ar@{->>}[r] \ar[d] & \ast \ar[d] \\
F_{x,k}^{k-1,k-1} \ar@{->>}[r] & \dprod[|s|=k]{\tJ(x,s)} F_{s,k-1}^{k-2,k} \ar@{->>}[r]
    & \dprod[|s|=k]{\tJ(x,s)} F_{s,k-1}^{k-2,k-1}
}
yields a pullback with horizontal fibrations and with \w{F_{x,k}^{k-1,k-1}} (and of course $\ast$)
contractible, so the result (with \w[)]{2 \leq j=k} also follows by induction.
\end{proof}

With these conventions, each factor in the product
\w{\Yof{x} \to \Omega^{j-1} \Yof{v}} is a $j$-ary Toda bracket by construction,
and vanishing of the product is equivalent to vanishing of each factor.

\begin{thm}
In the fully reduced case, all higher operations decompose into a sequence of
Toda brackets of order no greater than the degree of the first target object in
the string.
\end{thm}

\appendix
%
%
\section{Background Material}\label{abm}

We collect here a number of basic facts about model categories needed in this paper
and one non-standard lemma included for ease of reference elsewhere. We refer
the reader to \cite[\S\S 7.1-7.3]{PHirM} for the basics on model categories and homotopy
assumed for this appendix.

\begin{notn}
Given two maps \w[,]{f, g: X \to Y} we write \w{f\sim^r g}
if the maps are right homotopic, and \w{f\siml g} if the maps are left homotopic.
\end{notn}

\begin{lemma}[Homotopy Lifting Property]\label{lhlp}
Suppose we have the solid diagram with $q$ a fibration and $T$ cofibrant:
\mydiagram[\label{eqlhlp}]{
T \ar@{>}[dr]_{\psi} \ar[r]^{f} & Y \ar@{->>}[d]^{q} \\
 & Z
}
 \noindent Then there is a homotopy \w{\psi \siml q\circ f} if and only if there is a
 map \w{f': T \to Y} with a homotopy \w{f' \siml f} such that \w[.]{\psi = q\circ f'}

 Dually, if $Z$ is fibrant and $f$ is a cofibration then there is a homotopy \w{\psi \sim^r q\circ f}
 precisely when there is a  map \w{q': Y \to Z} with a homotopy \w{q' \sim^r q} such that
 \w[.]{\psi = q'\circ f}
\end{lemma}

\begin{proof}
Assume $q$ is a fibration. Let
$$
T\amalg T \stackrel{i_1\sqcup i_2}{\longrightarrow} \Cyl(T) \stackrel{p}{\longrightarrow} T
$$
be a factorization of the fold map \w{T\amalg T \stackrel{1_T\amalg 1_T}{\longrightarrow} T}
such that \w{i_1\sqcup i_2} is a cofibration and $p$ is a weak equivalence.
Cofibrancy of $T$ implies
\w{i_1:T\to\Cyl(T)} is an acyclic cofibration by \cite[7.3.7]{PHirM}.
Given a homotopy \w{H: \Cyl(T) \to Z} with
\w{H\circ i_{1}= q\circ f} and \w[,]{H\circ i_{2}=\psi} we
may use the left lifting property in
\begin{myeq}
\xymatrix@R=25pt{
T \ar@{ >->}[d]_{\simeq}^{i_{1}} \ar[r]^{f} &
         Y \ar@{->>}[d]^{q} \\
\Cyl(T) \ar@{.>}[ru]^{\hat{H}} \ar[r]^{H}      & Z
}
\end{myeq}
\noindent to factor $H$ as \w[,]{q\circ\hat{H}} and set
\w[.]{f':=\hat{H}\circ i_{2}}
If $f$ is instead a cofibration, use the dual argument.
\end{proof}

\begin{lemma}[Homotopy Pullback Property]\label{lhpp}
Suppose we have the following solid diagram where the square is a pullback, $T$
is cofibrant, and the two vertical maps are fibrations.
\mydiagram[\label{eqlhpp}]{
T \ar@{.>}[dr]^{g} \ar@/_1em/[ddr]_{p} \ar@/^1em/@{.>}[drr]^{f} \\
& W \ar@{}[dr] |<{\mbox{\large{$\lrcorner$}}} \ar@{->>}[d]_{r} \ar[r]^{j} & Y \ar@{->>}[d]^{q} \\
& X \ar[r]_{i} & Z
}
\noindent Then there is a dotted map \w{f:T \to Y} with a homotopy  \w{q\circ f \siml ip}
precisely when there is a dotted map \w{g:T \to W} with a homotopy \w{j\circ g \siml f}
and \w[.]{r\circ g=p}
\end{lemma}

\begin{proof}
Suppose there is a homotopy \w[.]{q\circ f \siml i\circ p}
Since $T$ is cofibrant and $q$ is a fibration,  the
Homotopy Lifting Property (with \w[)]{\psi = i\circ p} produces \w{f': T \to Y} homotopic to $f$,
such that \w[.]{q\circ f' = i\circ p} Since the square is a pullback, there is a map \w{g: T \to W} such that
\w{j\circ g = f'} and \w[.]{r\circ g = p} Since \w[,]{f\siml f'} we conclude that \w[.]{f \siml j\circ g}
\end{proof}

\begin{cor}\label{cnullpull}
  If $X$ is cofibrant, \w{k:X\to Y} is any pointed map,  and \w{h:Y \to Z} is a pointed fibration,
  then the composite \w{h\circ k:X \to Z} is null-homotopic if and only if there exists some
\w[,]{k':X\to Y}  left homotopic to $k$, which factors through \w[.]{\fib(h)}
\end{cor}

\begin{lemma}[Homotopy Ladder Property]\label{lthpp2}
 Suppose we are given the following diagram in which both squares are (strict) pullbacks,
 $T$ is cofibrant, the indicated horizontal maps are fibrations, and the outer diagram commutes
 up to homotopy:
\mydiagram[\label{eqthpp}]{
  T \ar@{-->}[dr]^{\kappa} \ar@/_1em/[dddr]_{\varphi}
  \ar@/_1em/@{.>}[ddr]^{\theta} \ar@{>}@/^1em/[drr]^{\sigma} \\
& U \ar@{}[dr] |<{\mbox{\large{$\lrcorner$}}} \ar@{->>}[r]^{r} \ar[d]_{\Phi} & V \ar[d]^{t} \\
& W \ar@{}[dr] |<{\mbox{\large{$\lrcorner$}}} \ar@{->>}[r]^{p} \ar[d]_{q} & X \ar[d]^{s} \\
& Y \ar@{->>}[r]_{u} & Z ~.
}
\noindent Consider the following three statements:
\begin{enumerate}
\item[(1)] There is a map \w{\kappa: T \to U} such that \w[,]{\sigma = r\circ \kappa} and there are
(left) homotopies \w{\theta \siml \Phi\circ \kappa} and \w[.]{\varphi \siml q\circ \Phi\circ \kappa}
\item[(2)] \w[,]{\varphi \siml q\circ \theta} and there is a map \w{\theta': T \to W} homotopic to
  $\theta$ such that \w[.]{p\circ \theta' = t\circ \sigma}
\item[(3)] There is a map \w{\theta': T \to W} homotopic to $\theta$ such that $\varphi$ is homotopic to
 \w{\varphi' := q\circ \theta'} and \w[.]{u\circ \varphi' = s\circ t\circ \sigma}
\end{enumerate}
\noindent Then \w[.]{(1) \Leftrightarrow (2) \Rightarrow (3)} Furthermore, if $s$ is a
monomorphism, then (1), (2), and (3) are all equivalent.
\end{lemma}

\begin{proof}
\noindent \w[:]{(1) \Longrightarrow (2)} Since \w[,]{\theta \siml \Phi\circ \kappa} it follows that
\w[.]{\varphi \siml q\circ \Phi \circ \kappa \siml q\circ \theta}
Since \w[,]{p\circ \theta \siml p\circ \Phi\circ \kappa = t\circ \sigma} applying the Homotopy
Lifting Property (with \w{q = p} and \w[),]{f = \theta} to \w{\psi = t\circ \sigma} there exists
\w{\theta' \siml \theta} with \w[.]{p\circ \theta'=t\circ \sigma}

\noindent \w[:]{(2) \Rightarrow (1)} Let \w{\theta' \siml \theta} with \w[,]{p\circ \theta' = t\circ \sigma}
and let \w[.]{\varphi' := q\circ \theta'} Then
$$
u\circ \varphi' = u\circ q\circ \theta' = s\circ p\circ \theta' = s\circ t\circ \sigma
$$
\noindent Since the outside rectangle is a pullback, there exists $\kappa: T \to U$ such that
\w[.]{\theta' = \Phi\circ \kappa$ and $\sigma = r\circ \kappa} Thus
\w[.]{\theta \siml \theta' = \Phi\circ \kappa}
Also, \w[.]{\varphi \siml q\circ \theta \siml q\circ \Phi \circ \kappa}

\noindent \w[:]{(2) \Rightarrow (3)} Given \w{\theta' \siml \theta} such that
\w[,]{p\circ \theta' = t\circ \sigma} set \w[;]{\varphi' := q\circ \theta'}  then
\w[.]{\varphi \siml q\circ \theta \siml q\circ \theta' = \varphi'}
Also, from the squares commuting
$$
u\circ \varphi' = u\circ q\circ \theta' = s\circ p\circ \theta' = s\circ t\circ \sigma
$$

Finally, we assume that \w{s: X \to Z} is a monomorphism. We show that
\w[.]{(3) \Rightarrow (2)} From the squares commuting, we have
$$
s\circ t\circ \sigma = u\circ \varphi' = u\circ q\circ \theta' = s\circ p\circ \theta'
$$
\noindent Thus \w[,]{t\circ \sigma = p\circ \theta'} because $s$ is a monomorphism,
and \w{\varphi \siml q\circ \theta} as above.
\end{proof}

\begin{cor}\label{cthpp2}
In \wref{eqthpp} assume again that the squares are pullbacks, $T$ is cofibrant,
and the horizontal maps are fibrations. Assume further that
\w[.]{u\circ \varphi \siml s\circ t\circ \sigma} Then we have the following:
\begin{enumerate}
\item There exists a map \w{\kappa: T \to U} such that \w{\sigma = r\circ\kappa}
  and \w[.]{\varphi \siml q\circ \Phi\circ \kappa}
\item There exists a map \w{\theta: T \to W} such that
\w{\varphi \siml q\circ \theta} and \w[.]{p\circ \theta = t\circ \sigma}
\end{enumerate}
Moreover, if $s$ is additionally a monomorphism then there is a homotopy
\w[.]{\theta \siml\Phi\circ \kappa}
\end{cor}

\begin{proof}
For (1), since \w[,]{u\circ \varphi \siml s\circ t\circ \sigma} by the Homotopy Pullback Property,
there is a map \w{\varphi': T \to U} homotopic to $\varphi$
such that \w[.]{u\circ\varphi' = s\circ t\circ \sigma} Since the outer rectangle is a pullback, there is a
map \w{\kappa: T \to U} such that \w{\varphi' = q\circ \Phi\circ \kappa} and \w[.]{\sigma = r\circ \kappa}
Thus \w[.]{\varphi \siml q\circ \Phi\circ \kappa}

For (2), we have \w[.]{u\circ \varphi \siml s\circ t\circ \sigma}
Again, by the Homotopy Pullback Property,
there is a map \w{\varphi' \siml \varphi} such that \w[,]{u\circ \varphi' = s\circ t\circ \sigma}
so since the bottom square is a pullback, there is a map \w{\theta: T \to W} with
\w{t\circ \sigma = p\circ \theta} and \w[,]{\varphi' = q\circ \theta} and so
  \w[.]{\varphi \siml q\circ \theta}

Finally, \w[,]{u\circ \varphi' = u\circ q\circ \theta = s\circ p\circ \theta = s\circ t\circ \sigma}
so if $s$ is a monomorphism, we may conclude from Lemma \ref{lthpp2} that
\w[.]{\theta \siml \Phi\circ \kappa}
\end{proof}

We have the duals of Lemma \ref{lhpp}, Corollary \ref{cnullpull},
Lemma \ref{lthpp2} and Corollary \ref{cthpp2}:

\begin{lemma}\label{dlhpp}
Suppose the following square is a pushout, $V$ is fibrant, and the two horizontal maps are
cofibrations:
\mydiagram[\label{eqdlhpp}]{
  W \ar@{}[dr] |>{\mbox{\large{$\ulcorner$}}}  \ar@ { >->}[r]^{i} \ar[d]^{\alpha} &
  Y \ar[d]^{\beta} \ar@/^1em/@{.>}[ddr]^{f} & \\
X \ar@ { >->}[r]_{j} \ar@/_1em/[drr]_{p} & Z \ar@{.>}[dr]^{g} & \\
& & V ~.
}
\noindent Then there is a dotted map \w{f:Y \to V} with a homotopy \w{p\circ \alpha \simr f\circ i}
precisely when there is a dotted map \w{g:Z \to V} with a homotopy \w{g\circ\beta\simr f} and  \w[.]{g\circ j = p}
\end{lemma}

\begin{cor}\label{cdnullpull}
  If \w{k:X \to Y} is a pointed cofibration and \w{h:Y \to Z} is any pointed map with Z fibrant,
  the composite \w{h\circ k:X \to Z} is null-homotopic if and only if there exists a
  map \w[,]{h':Y \to Z} right homotopic to $h$, which factors through \w[.]{\cofib{k}}
\end{cor}

\begin{lemma}\label{lthpp3}
Suppose we are given the following diagram in which both squares are (strict) pushouts,
$T$ is fibrant, the indicated horizontal maps are cofibrations, and the outer diagram commutes
up to homotopy:
\mydiagram[\label{eqthpp3}]{
  U \ar@{}[dr] |>{\mbox{\large{$\ulcorner$}}} \ar@ { >->}[r]^{r} \ar[d]^{\Phi} &
  V \ar[d]^{t} \ar@/^1em/[dddr]^{\varphi} & \\
  W \ar@{}[dr] |>{\mbox{\large{$\ulcorner$}}} \ar@ { >->}[r]^{p} \ar[d]^{q} &
  X \ar[d]^{s}  \ar@/^1em/@{.>}[ddr]_{\theta} &\\
Y \ar@ { >->}[r]_{u} \ar@{>}@/_1em/[drr]_{\sigma} & Z \ar@{-->}[dr]^{\kappa}  & \\
  & & T ~.
}
\noindent Consider the following three statements:
\begin{enumerate}
\item[(1)] There exists a map \w{\kappa: Z \to T} such that \w{\sigma = \kappa\circ u} and there are
(right) homotopies \w{\theta \simr  \kappa\circ s} and \w[.]{\varphi \simr \kappa\circ s\circ t}
\item[(2)] \w{\varphi \simr \theta\circ t} and there is a map \w[,]{\theta': X \to T} homotopic to
  $\theta$, such that \w[.]{\theta'\circ p = \sigma\circ q}
\item[(3)] There is a map \w{\theta': X \to T} homotopic to $\theta$ such that $\varphi$ is homotopic
  to \w[,]{\varphi' :=\theta'\circ t} and \w[.]{\varphi'\circ r = \sigma\circ q\circ \Phi}
\end{enumerate}
Then \w[.]{(1) \Leftrightarrow (2) \Rightarrow (3)} Furthermore, if $\Phi$ is
an epimorphism, then (1), (2), and (3) are all equivalent.
\end{lemma}

\begin{cor}\label{cthpp3}
In \wref[,]{eqthpp3} assume again that the squares are pushouts, $T$ is fibrant,
and the horizontal maps are cofibrations. Assume further that
\w[.]{\varphi\circ r \simr \sigma\circ q\circ \Phi} Then we have the following:
\begin{enumerate}
\item There exists a map \w{\kappa: Z \to T}
such that \w{\sigma = \kappa\circ u} and \w[.]{\varphi \simr q\circ \kappa\circ s\circ t}
\item There exists a map \w{\theta: X \to T} such that
\w{\varphi \simr \theta\circ t} and \w[.]{\theta\circ p = \sigma\circ q}
\end{enumerate}
Moreover, if $\Phi$ is additionally an epimorphism then there is a homotopy
\w[.]{\theta\simr\kappa\circ s}
\end{cor}

We define the \emph{reduced path object} \w{PW} associated to a pointed object $W$ by the
pullback
\mydiagram[\label{rpodef}]{
  PW \ar@{}[dr] |<{\mbox{\large{$\lrcorner$}}}  \ar@{->>}[d]_{p_W} \ar[r]^-{j} &
  \Path(W) \ar@{->>}[d]^{p_1\times p_2} \\
W \ar[r]_-{1_W\times 0} & W\times W
}
\begin{lemma}\label{rpoprop}
If $W$ is fibrant, then \w{PW} is weakly contractible. Furthermore, if \w{f: X \to W} is pointed,
then $f$ is right null-homotopic precisely when $f$ factors as
\w[.]{X \to PW \stackrel{p_W}{\to} W} 
%
\end{lemma}
\begin{proof}
First, the diagram \ref{rpodef} can be expanded to the pullback
\mydiagram[\label{rpowc}]{
PW \ar@{}[dr] |<{\mbox{\large{$\lrcorner$}}} \ar@{->>}[d] \ar[r]^-{j} & \Path(W)
\ar@{->>}[d]^{\proj_2\circ (p_1\times p_2)} \\
\ast \ar[r] & W
}
Since $W$ is fibrant, the right hand vertical map is a trivial fibration, by \cite[7.3.7]{PHirM}. Hence
the left hand vertical map is a trivial fibration, by \cite[7.2.12]{PHirM}. Thus \w{PW} is weakly
contractible.

If \w{f: X \to W} is null-homotopic, there is a map \w{H : X \to \Path(W)} with
\w{p_1\circ H = f} and \w[.]{p_2\circ H = 0} From the first factorization, and the pullback property
of \wref[,]{rpodef} there is a map \w{\phi: X \to PW} such that \w[.]{f = p_W\circ \phi}
%
\end{proof}

We similarly define the \emph{reduced cone} \w{CX} on a pointed object $X$ by the pushout
\mydiagram[\label{rpcdef}]{
  X\amalg X \ar@{}[dr] |>{\mbox{\large{$\ulcorner$}}} \ar@{ >->} [d] _{i_1\amalg i_2} \ar[r]^{1_X\amalg 0} &
  X  \ar@{ >->}[d]^{i_X} \\
\Cyl(X) \ar[r]_{} & CX
}
\begin{lemma}\label{rcoprop}
  If $X$ is cofibrant, then \w{CX} is weakly contractible. Furthermore, if \w{f: X \to W} is a pointed map,
  then $f$ is left null-homotopic precisely when $f$ factors as
\w[.]{X \stackrel{i_X}{\rightarrow} CX \to W}
%
\end{lemma}

\begin{lemma}\label{pathpull}
  Let $X$ be cofibrant and both $Z$ and $W$ fibrant. If the composite \w{g\circ h\circ k}
  is right null-homotopic, then the shorter composite \w{h\circ k} is also right null-homotopic if
  and only if there is a null homotopy $\phi$ of \w{g\circ h\circ k} such that the solid
  commutative diagram
\mydiagram[\label{eqpathpull}]{
X \ar[rrrr]^{k} \ar@{.>}[drr]^{\psi} \ar@/_2em/[ddrrr]_{\phi}
 & & & & Y \ar[d]^{h} \\
& & PF_g \ar@{.>>}[r] ^{p_{F_g}} &
F_g \ar @{} [dr] |>>>>>>>>>>>{\mbox{\large{$\lrcorner$}}} \ar@{.>}[d]^j \ar@{.>>}[r]^i &
Z\ \ar[d]^{g} \\
& & & PW \ar@{->>}[r]_{p_W} & W
}
    extends to the full diagram above, with $\psi$ a null homotopy for \w{h\circ k}
    and \w{F_g} the pullback of $g$ along \w[.]{p_W}
\end{lemma}

\begin{proof}
  Suppose the composite \w{g\circ h\circ k} is null-homotopic.
  Then Lemma \ref{rpoprop} gives a factorization \w{g\circ h\circ k = p_W\circ \phi} in
  \wref[.]{eqpathpull}  Since \w{p_W}
is a fibration, so is $i$. If \w{h\circ k} is also null-homotopic then this composite
factors as \w[,]{h\circ k = p_Z\circ \kappa} for some \w[.]{\kappa : X \to PZ}
Now factor $\kappa$ as \w[,]{X \stackrel{\kappa'}{\to} V \stackrel{q}{\to} PZ}
with \w{\kappa'} a cofibration and $q$ a trivial fibration. Since $X$ is cofibrant
and \w{PZ} is weakly contractible by Lemma \ref{rpoprop}, \w{\ast \to V}
is a trivial cofibration. Therefore, \w{p_Z \circ q}
lifts to a map \w{\eta:V\to PF_g}
with \w[.]{i\circ p_{F_g}\circ \eta = p_Z\circ q} Setting
\w{\psi:= \eta\circ \kappa'} makes the whole diagram commute.
\end{proof}

The dual version is:

\begin{lemma}\label{dpathpull}
Let $Y$ be fibrant and both $Z$ and $W$ cofibrant.
Suppose the composite $k\circ h\circ g$ is known to be left null-homotopic.  Then the shorter
composite \w{k\circ h} is also left null-homotopic if and only if for some null homotopy $\phi$
of \w[,]{k\circ h\circ g} the solid commutative diagram
\mydiagram[\label{eqdpathpull}]{
  W \ar@{}[dr] |>{\mbox{\large{$\ulcorner$}}} \ar@ { >->}[r]^{i_W} \ar[d]_{g} &
  CW \ar@{.>}[d] \ar@/^2em/[ddrrr]^{\phi} &&& \\
Z \ar@{ >.>}[r] \ar[d]_{h} & \mapcone{g} \ar@{ >.>}[r] &
C\mapcone{g} \ar@{.>}[rrd]^{\psi} && \\
X \ar[rrrr]^{k} &&&& Y
}
\noindent extends to the full diagram above, with $\psi$ giving a null homotopy for \w{k\circ h}
    and \w{\mapcone{g}} the pushout of $g$ along \w[.]{i_W}
\end{lemma}

%
%
\section{Indeterminacy}\label{abind}

For most higher homotopy operations, one cannot expect a closed formula for the indeterminacy
of operations of the type provided by \cite[Lemma 1.1]{TodC} for the classical (secondary) Toda bracket.
This is because tertiary and higher operations depend on choices made for the vanishing of the
lower order operations, and the amount of choice remaining might vary for different sets of earlier
choices.

However, if we take these earlier choices as given, within the inductive framework described
here the only remaining source of indeterminacy is in the choice of the specific map
\w{\varphi'} which makes the outer diagram in \wref{eqthpp} commute on the nose, and how
that choice affects the resulting lift \w[.]{\theta'}
Note that the homotopy class \w{[\varphi']=[\varphi]}
is then fixed, as is the actual map \w[.]{u\circ\varphi'=s\circ t\circ \sigma:T\to Z}
To help keep track of all this, in this appendix $\varphi$ will denote our initial choice of the map
with the induced lift $\theta$, while \w{\varphi'} will denote some other choice, with induced lift
\w[.]{\theta'} We now investigate how changing $\varphi$ to \w{\varphi'} changes $\theta$ to
\w[,]{\theta'} as maps \w[:]{T \to W}

Given  $\varphi$, a choice of \w{\varphi'} such that \w{u \circ \varphi=u \circ \varphi'}
corresponds uniquely to a map into the pullback
\mydiagram[\label{eqleftpb}]{
T \ar@/_1em/[ddr]_{\varphi} \ar@/^1em/@{-->}[drr]^{\varphi'} \ar@{-->}[dr]\\
& \Yu \ar @{} [dr] |>>>>>>>>>>>{\mbox{\large{$\lrcorner$}}} \ar[d]_{u'} \ar[r] & Y \ar[d]^{u} \\
& Y \ar[r]^{u} & Z
}
\noindent while a choice of such a map \w{\varphi'} equipped with
a (right) homotopy \w{H:\varphi \simr \varphi'} corresponds to a map into the pullback
\mydiagram[\label{eqrightpb}]{
T \ar@/_1em/[ddr]_{\varphi} \ar@/^1em/@{-->}[drr]^{H} \ar@{-->}[dr]\\
& \oYu \ar @{} [dr] |>>>>>>>>>>>>>{\mbox{\large{$\lrcorner$}}} \ar[d]_{\overline u'} \ar[r] &
\Path(Y) \ar[d]^{(1 \times u) \circ m} \\
& Y \ar[r]^-{1\top u} & Y \times Z
}
\noindent where \w{Y\xra{i_{y}}\Path(Y)\xepic{m} Y\times Y} is a path factorization as in \wref[.]{rpodef}
In fact, taking a further pullback
\mydiagram[\label{oneForIndet}]{
  \oWpu \ar @{} [drr] |<<{\mbox{\large{$\lrcorner$}}} \ar@{->>}[d]_{\ovp'} \ar[rr] &&
  \oYu \ar@{->>}[d]^{\overline u'} \\
W \ar[rr]^{q} && Y
}
\noindent we find that the image of the left vertical map \w{\ovp'} is essentially the indeterminacy (see
Corollary \ref{cindet} below).

Note that there is a canonical choice of induced map \w{\psi:T \to \Yu} in \wref[,]{eqleftpb} corresponding
to \w[,]{\varphi'=\varphi} and a similar canonical choice of induced map \w{\overline \psi:T \to \oYu}
in \wref[,]{eqrightpb} corresponding to the canonical self-homotopy \w{H\sb{\varphi}} of $\varphi$ (namely,
the composite \w[),]{T\xra{\varphi} Y \xra{i_{y}}\Path(Y)} which will be used below.

Given a map \w[,]{u:Y \to Z}  consider the following pullback grid:
\mydiagram[\label{IndPull}]{
  \oYu \ar@{} [drr] |<<{\mbox{\large{$\lrcorner$}}} \ar@{-->>}[d]^{\overline u}
  \ar@{->>}@/_2em/[dd]_{\overline u'} \ar[rr] && \Path(Y) \ar@{->>}[d]_{m} && \\
  \Yu \ar @{} [drr] |<<{\mbox{\large{$\lrcorner$}}} \ar@{->>}[d]^{u'} \ar[rr] &&
  Y \times Y \ar@{} [drr] |<<{\mbox{\large{$\lrcorner$}}} \ar@{->>}[d]_{1 \times u} \ar[rr]^-{\proj_2} &&
  Y \ar@{->>}[d]^{u} \\
Y \ar[rr]^-{1 \top u} && Y \times Z  \ar[rr]^-{\proj_2} && Z
}

\begin{notn}\label{defn.lifts}
Assume given four maps \w[,]{u:Y \to Z} \w[,]{\varphi:T \to Y} \w[,]{v:B \to Y} and \w[.]{\rho:A \to Y}
\begin{enumerate}
\renewcommand{\labelenumi}{(\alph{enumi})~}
\item  The pointed set
  \w[,]{\{\varphi':T \to Y ~\vert~u \circ \varphi'= u\circ \varphi \}} based at $\varphi$ itself,
will be denoted by \w[.]{\varu{\varphi}}
\item The pointed set
  \w{\{H:T \to\Path(Y)~\vert~H:\varphi \simr \varphi', u \circ \varphi'= u\circ \varphi \}}
  of (right) homotopies, based at 
  \w[,]{H\sb{\varphi}} 
  will be denoted by \w[.]{\ovaru{\varphi}}
\item  The set \w{\{\sigma:A \to B ~\vert~ v \circ \sigma = \rho\}}  of lifts of $\rho$ with respect to $v$
  will be denoted by \w[.]{\lift{v}{\rho}}
\end{enumerate}
In accordance with Remark \ref{rassfibcof}, we can disregard the distinction between the left
homotopies appearing in the first half of Appendix \ref{abm} and the right homotopies we have here.
\end{notn}

\begin{remark}\label{rvarlift}
  From the pullback properties of the constructions above we see that there are natural bijections of
  pointed sets
\w{\varu{\varphi} \cong \lift{u'}{\varphi}} and \w[,]{\ovaru{\varphi} \cong \lift{\overline u'}{\varphi}}
where \w{\lift{u'}{\varphi}} is based at $\psi$ and \w{\lift{\overline u'}{\varphi}} is based at
$\overline \psi$.
\end{remark}

We then have:

\begin{lemma}
Given \w{\varphi=q \circ \theta:T \to Y} with \w[,]{p \circ \theta=t \circ \sigma}
there is a natural bijection of sets \w[,]{\ovaru{\varphi} \cong \lift{\ovp'}{\theta}}
where \w[.]{\ovp':=p' \circ \ovp}
\end{lemma}

\begin{proof}
We may expand \wref{IndPull} into:
\mydiagram[\label{IndPullb}]{
  & \oWpu \ar[dl] \ar@{-->>}|(.5){\hole}[dd]_(0.65){\ovp}
  \ar@/^1.4em/|(.26){\hole}|(.75){\hole}|(.82){\hole}[dddd]^(0.37){\ovp'} \ar[rr] &&
    \Prel \ar[dl] \ar@{->>}[dd]|(.8){\hole} \\
 \oYu \ar@{-->>}[dd]^(0.3){\overline u} \ar[rr] && \Path(Y) \ar@{->>}[dd] \\
& \Wp \ar[dl]^{q'} \ar@{->>}|(.5){\hole}|(.6){\hole}[dd]_(0.3){p'} \ar[rr] \ar@/^1.7em/[rrrr]^(0.5){p''} &&
Y \times W \ar[dl]^{1\times q} \ar@{->>}|(.5){\hole}|(.6){\hole}[dd]^(0.25){1\times p} \ar[rr]_(0.6){\proj_2}
    && W \ar[dl]^{q} \ar@{->>}[dd]^(0.5){p} \\
\Yu \ar@{->>}[dd]^(0.3){u'} \ar[rr] \ar@/_1em/[rrrr]_(0.4){u''}
&& Y \times Y \ar@{->>}|(.12){\hole}[dd]^(0.3){1 \times u} \ar[rr]^(0.75){\proj_2} &&
Y \ar@{->>}[dd]^(0.3){u} \\
& W \ar[dl]^{q} \ar[rr]|(.5){\hole}^(0.3){q \top p} &&
Y \times X \ar[dl]^{1 \times s} \ar[rr]|(.57){\hole}^(0.76){\proj_2} && X\ar[dl]^{s} \\
Y \ar[rr]^-{1 \top u} && Y \times Z \ar[rr]^(0.6){\proj_2} && Z
}
\noindent Since the rightmost face is a pullback (by assumption), as are both the front and left
long rectangular vertical
faces (by construction), the lower leftmost face, and hence the upper leftmost face, are pullbacks, too.
We define \w{\Prel} by making the upper rightmost face a pullback, so that the back upper
vertical face is, too.

We think of \w{\varphi:T \to Y} as mapping to the front lower left $Y$, and \w{\theta:T\to W} to the
back lower left $W$, with \w{\varphi':T \to Y} mapping to the front right $Y$, and \w{\theta':T\to W} to the
back right $W$. Since \w[,]{u \circ \varphi'=u \circ \varphi} the lower pullback rectangle in
\wref{IndPull} implies that \w{(\varphi,\varphi')} induce a map \w{F:T\to\Yu}  and thus
\w[.]{\widehat{F}:T\to\Wp}
Since also \w{u \circ \varphi= s \circ p \circ \theta=s \circ t \circ \sigma} and
a right homotopy \w{H:\varphi \simr \varphi'} is a map
\w{H:T\to\Path(Y)} which, together with \w[,]{\varphi \top \theta':T\to Y\times W} induces
\w[,]{\widehat{H}:T\to\Prel} together with $\widehat{F}$, these induce a lift of $\theta$ along
\w[.]{\ovp'} Conversely, any lift \w{\widehat{\theta}:T\to\oWpu} of $\theta$ along
\w{\ovp'} yields $\widehat{H}$, and thus $H$, by projecting along the structure maps of the
top pullback square.
\end{proof}

\begin{remark}
When \w[,]{Y \sim \ast} we have \w[,]{\Path(Y) \stackrel{\sim}{\to} Y \times Y} so
\w{\oYu \stackrel{\sim}{\to} \Yu \simeq \Omega Z} and \w[.]{\oWpu \simeq W \times \Omega Z}
In this case a map \w{T \to \oWpu} thus corresponds up to homotopy, to a choice of map $\theta$,
together with  a homotopy class in \w{[T ,\,\Omega Z]} (adjoint to the indeterminacy construction
  of \cite[\S 1]{SpanS}). Note that each of the vertical faces in \wref{IndPullb} is a pullback over a
  fibration, so they are homotopy-meaningful.
\end{remark}

The indeterminacy of our operations is then described by the following.

\begin{cor}\label{cindet}
  Given \w{\varphi=q \circ \theta:T \to Y} (also satisfying \w[)]{p \circ \theta=t \circ \sigma}
  in \wref[,]{eqthpp}  the indeterminacy in our operation produced by varying $\varphi$ lies in the
  image of \w[,]{\ovp''_{\#}:[T,\oWpu] \to [T,W]} where \w[.]{\ovp''=p'' \circ \overline{p}}

  In fact, we can restrict to the fiber of
  \w{\ovp'_{\#}$ over $[\theta]}  (the subset consisting of those homotopy classes containing
  an element of \w[).]{\lift{\ovp'}{\theta}}
\end{cor}

\begin{proof}
  In \wref{IndPullb} each choice of a lifting \w{\theta'} of \w{\varphi'\sim\varphi} has the form
  \w{p'' \circ \ovp \circ \rho} for some \w[.]{\rho:T \to \oWpu} Thus \w[,]{\ovp''_{\#} [\rho]=[\theta']}
  as required. By restricting to those $\rho$ with
\w[,]{[\ovp' \circ \rho]=\ovp'_{\#} [\rho]=[\theta]} we can apply Lemma \ref{lhlp}
to produce a different representative \w{[\rho']=[\rho]} with \w[,]{\ovp' \circ \rho' = \theta}
producing the improved \w[.]{\theta'}
\end{proof}

\end{document}